\documentclass[11pt]{article}

\usepackage{xypic}
\usepackage{amsgen}
\usepackage{amsmath}
\usepackage{amstext}
\usepackage{amsbsy}
\usepackage{amsopn}
\usepackage{amsfonts}
\usepackage{amssymb}
\usepackage{eepic}
\usepackage[pdftex]{graphicx}
\usepackage{epsf}
\usepackage{pstricks}
\xyoption{all}

\def\Box{\square}
\def\edge{\relbar\joinrel\relbar}
\def\mapright#1{\smash{\mathop{\longrightarrow}\limits^{#1}}}
\def\tra#1{\smash{\mathop{\mid\kern
-1pt\joinrel\relbar\joinrel\relbar}\limits^{*}_{#1}}}
\def\longtra#1{\smash{\mathop{\mid\kern
-1pt\joinrel\relbar\joinrel\relbar\joinrel\relbar}\limits^{*}_{#1}}}
\def\vlongtra#1{\smash{\mathop{\mid\kern-1pt\joinrel\relbar\joinrel\relbar\joinrel\relbar\joinrel\relbar}\limits^{*}_{#1}}}
\def\vvlongtra#1{\smash{\mathop{\mid\kern
-1pt\joinrel\relbar\joinrel\relbar\joinrel\relbar\joinrel\relbar\joinrel\relbar}\limits^{*}_{#1}}}
\def\vvvlongtra#1{\smash{\mathop{\mid\kern
-1pt\joinrel\relbar\joinrel\relbar\joinrel\relbar\joinrel\relbar\joinrel\relbar\joinrel\relbar}\limits^{*}_{#1}}}
\def\etra#1{\smash{\mathop{\mid\kern
-1pt\joinrel\relbar\joinrel\relbar}\limits_{#1}}}

\def\A{{\cal{A}}}

\def\iff{\Leftrightarrow}
\def\Rw{\Rightarrow}
\def\oo{\overline}

\def\wt{\widetilde}
\def\wh{\widehat}
\def\dw{\!\downarrow}
\def\B{{\cal{B}}}

\def\L{{\cal{L}}} %{\mathrel{{\mathcal L}}}
\def\M{{\cal{M}}}
\newcommand{\N}{{\rm I}\kern-2pt {\rm N}}
\def\SB{\mathbb{SB}}
\def\BB{\mathbb{B}}

\def\im{\mbox{Im}\,}

\def\BR{\mbox{BR}}

\def\deg{\mbox{deg}}

\def\cl{\mbox{Cl}\,}
\def\clos{\mbox{Cl}}
\def\diam{\mbox{diam}\,}
\def\lat{\mbox{Lat}\,}
\def\latti{\mbox{Lat}}

\def\FLg{\mbox{FLg}}
\def\gth{\mbox{gth}\,}
\def\lin{\mbox{Lin}\,}
\def\poi{\mbox{Pt}\,}

\def\rk{\mbox{rk}\,}

\def\geo{\mbox{Geo}\,}
\def\matro{\mbox{Mat}\,}
\def\matroid{\mbox{Mat}}
\def\flats{\mbox{Fl}\,}
\def\isflats{\mbox{ISFl}\,}
\def\fisflats{\mbox{FISFl}\,}

\def\per{\mbox{Per}\,}
\def\ker{\mbox{Ker}\,}

\def\max{\mbox{max}}
\def\maxdeg{\mbox{maxdeg}\,}
\def\mindeg{\mbox{mindeg}\,}

\def\min{\mbox{min}}

\def\het{\mbox{ht}\,}
\def\hei{\mbox{ht}}

\def\P{{\cal{P}}}

\def\RR{\mathbb{R}}

\def\G{{\cal{G}}}

\def\Y{{\cal{Y}}}
\def\Zcal{{\cal{Z}}}

\def\p{\varphi}

\def\inv{^{-1}}

\def\bi{\begin{itemize}}
\def\ei{\end{itemize}}
\def\beq{\begin{equation}}
\def\eeq{\end{equation}}

\newtheorem{T}{Theorem}[section]
\newcommand{\bt}{\begin{T}}
\newcommand{\et}{\end{T}}
\newcommand{\ftd}{$\square$\end{T}}

\newtheorem{Proposition}[T]{Proposition}
\newcommand{\bp}{\begin{Proposition}}
\newcommand{\ep}{\end{Proposition}}
\newcommand{\fpd}{$\square$\end{Proposition}}

\newtheorem{Definition}[T]{Definition}
\newcommand{\bd}{\begin{Definition}}
\newcommand{\ed}{\end{Definition}}

\newtheorem{Lemma}[T]{Lemma}
\newcommand{\bl}{\begin{Lemma}}
\newcommand{\el}{\end{Lemma}}
\newcommand{\fld}{$\square$\end{Lemma}}

\newtheorem{Corol}[T]{Corollary}
\newcommand{\bc}{\begin{Corol}}
\newcommand{\ec}{\end{Corol}}
\newcommand{\fcd}{$\square$\end{Corol}}

\newtheorem{Result}[T]{Result}
\newcommand{\br}{\begin{Result}}
\newcommand{\er}{\end{Result}}
\newcommand{\frd}{$\square$\end{Result}}

\newtheorem{Remark}[T]{Remark}
\newcommand{\brem}{\begin{Remark}}
\newcommand{\erem}{\end{Remark}}
\newcommand{\fremd}{$\square$\end{Remark}}

\newtheorem{Example}[T]{Example}
\newcommand{\be}{\begin{Example}}
\newcommand{\ee}{\end{Example}}

\newtheorem{Problem}[T]{Problem}
\newcommand{\bq}{\begin{Problem}}
\newcommand{\eq}{\end{Problem}}

\newcommand{\proof}
   {\par\medbreak\noindent{\bf Proof}.\enspace}

\newcommand{\qed}{%\hfill
$\Box$
\par\bigbreak}

\newlength{\lengtha} 
\newlength{\lengthb} 
\normalbaselines
\def\abstract#1{\par\bigskip
\begingroup\small
\baselineskip=12truept
\begin{center}ABSTRACT\end{center}
\par\medskip\par\noindent
\null\hfill\hbox{\vbox{\hsize=5truein\noindent#1}}
\hfill\null\par\endgroup\par}

\title{Matroids, hereditary collections and simplicial complexes having
  boolean representations\footnote{This is a preliminary report and
    the final paper may be amended, extended and perhaps have additional
    authors.}} 
\author{{\bf John Rhodes}\\ 
$ $\\ {\em Department of Mathematics, University of California, Berkeley,}\\ 
{\em California 94720, U.S.A.}\\
{\em email:} rhodes@math.berkeley.edu, BlvdBastille@aol.com\\
$ $\\
{\bf Pedro V. Silva}\\ $ $\\ {\em Centro de
Matem\'{a}tica, Faculdade de Ci\^{e}ncias, Universidade do
Porto,}\\ {\em R. Campo Alegre 687, 4169-007 Porto, Portugal}\\
{\em email:} pvsilva@fc.up.pt} \date{\today}

\begin{document}
\maketitle

\begin{center}\small
2010 Mathematics Subject Classification: 05E45, 05B35, 05B20, 06B15, 15B34

\medskip

Keywords: hereditary collection, simplicial complex, matroid,
representation, boolean matrix, lattice
\end{center}

\abstract{Inspired by the work of Izakhian and Rhodes, a theory of
representation of hereditary collections by boolean matrices is
developed. This corresponds to representation  by finite
$\vee$-generated lattices. The lattice of flats, defined for
hereditary collections, lattices and matrices, plays a central role in
the theory. The representations constitute a lattice and the minimal
and strictly join irreducible elements are studied, as well as
various closure operators.  
} 
 
\tableofcontents

\section{Introduction}

The background and prehistory for this paper goes something like the
following. In 2006 Zur Izakhian \cite{Izh} defined the notion of
independence for columns (rows) of a matrix with coefficients in a
supertropical semiring. Restricting this concept to the superboolean
semiring $\SB$ (see Subsection \ref{sssm}), and then to the subset of boolean
matrices (equals matrices with coefficients 0 and 1), we obtain the
notion of independence of columns (rows) of a boolean matrix. This
notion has several equivalent formulations (see Subsection \ref{sssm} of this
paper and references there), one involving permanent, another being
the following: if $M$ is an $m \times n$ boolean matrix, then a subset
$J$ of columns of $M$ is {\em independent} if and only if there exists
a subset $I$ of rows of $M$ with $|I| = |J| = k$ and the $k \times k$
submatrix $M[I,J]$ can be put into upper triangular form (1's on the
diagonal, 0's strictly above it, and 0's or 1's below it) by
independently permuting the rows and columns of $M[I,J]$. 

This is the notion of independence for columns of a boolean matrix we
will use in this paper. 
In 2008 the first author suggested that this
idea would have application in many branches of Mathematics and
especially in Combinatorial Mathematics. In this paper we apply it to
hereditary collections (also known as abstract simplicial complexes).
For other applications of this
notion to lattices, posets and matroids by Izhakian and the first author, see
\cite{IR1,IR2,IR3}. For applications to finite graphs by the present authors,
see \cite{RSil}.

If $M$ is an $m \times n$ boolean matrix with column space $C$,
then the set $\cal{H}$ of independent subsets of $C$ satisfies the
following axioms (see \cite{IR1,IR2}): 
\bi
\item[(H)] $\cal{H}$ is nonempty and closed under taking subsets
  (making it a {\em hereditary collection});
\item[(PR)] for all nonempty $J, \{ p \} \in {\cal{H}}$, there exists
  some $x \in J$ such that $( J \setminus \{x\}) \cup  \{ p \} \in
  {\cal{H}}$ (the {\em point replacement} property).
\ei
Hereditary collections arising from some boolean matrix $M$ as above
are said to be {\em boolean representable}. 
The elementary properties of such boolean representable
collections were considered in \cite{IR1,IR2,IR3} and it was shown in
\cite{IR2}  
that all matroids have boolean representations.

We describe now the structure and contents of this paper. 

In Section 2,
we present the basic results we need to deal with lattices, superboolean
matrices and hereditary collections. Note that all lattices are finite
in this paper, but many results admit extensions to arbitrary
lattices. 

In Section 3, we establish a bijection between boolean matrices and
$\vee$-generated lattices. Moving the idea of c-independent columns of
a boolean matrix,
via the bijection, over to lattices, we obtain the new idea (to us) that $X
\subseteq L$ ($L$ a finite lattice) is c-independent if and only if
there exists an ordering $X = \{ x_1, \ldots, x_k \}$ $(|X| = k)$ such
that
$$B < x_1 < (x_1 \vee x_2) < \ldots < (x_1 \vee \ldots \vee x_k).$$
Given an $m \times |E|$ boolean matrix $M$, $\flats M$ is
the closure under all intersections of those subsets of $E$ where the
rows of $M$ are zero. Equivalently, for a $\vee$-generated lattice
$(L,E)$ and assuming that the bottom element $B$ is not in $E$, we
have $\flats (L,E) = \{ \ell\dw \cap E \mid \ell \in L \}$. 

As mentioned before, we intend to consider hereditary collections
$(E,H)$ given by a boolean 
matrix $M$ of size $n \times |E|$. Equivalently, $(E,H)$ can be
described through a finite lattice $(L,E)$, $\vee$-generated by $E$,
with $H$ being the set of c-independent subsets of $L \subseteq 2^E$.  

First properties of c-independence are proved in Section 4, where the
key result is Proposition \ref{cindch}. Thus we can consider $(E,H)$
having boolean representations, or equivalently, lattice
representations, in their own right. By the main theorem of
\cite{IR2}, this includes all matroids.

A central thesis or viewpoint is that, perhaps, boolean
representations should replace matroids as the main object of study
in present day matroid theory. More on this at the end of this
Introduction.

In the central Section 5, we start by introducing the concept of flat
(or closed set) of an arbitrary hereditary collection 
$(E,H)$ and the lattice $\flats (E,H)$. This is done by generalizing
one of the formulae in matroid 
theory (they are not all equivalent in the general case, see
(\ref{circuits1}) and the paragraph following it). Then $(E,H)$ is
boolean representable if and only if, considering the transversals of
the partition of successive differences for some chain of $\flats
(E,H)$, equals $H$. See Proposition \ref{cindch} and Section 5.

In Section 6 we use $\vee$-maps to define a natural ordering on all
boolean (or lattice) representations of a 
boolean representable $(E,H)$. This leads to considering minimal boolean
representations of $(E,H)$, and also to the $\vee$ operator which
corresponds to ``stacking''  the matrices of the boolean representations.

Even for matroids, the minimal representation is a new idea (to us)
and it is important to get all  the minimal and sji (strictly join
irreducible) representations in the classical case.

The connections between $(E,H)$ and its lattice representations exist
at all levels. In Section 7 we relate the closure operator
induced by a hereditary collection 
with the closure operator induced by each of its representations.

In Section 8 we do a few examples. This includes computing all the
minimal and sji representations of the Fano matroid $(E,H)$ defined by
taking $E = \{ 1, \ldots,7 \}$ and $H$ equal to all subsets of $E$
with at most 3 elements except $125, 137, 146, 236, 247, 345, 567$.
 
Given integers $2 \leq a < b$, let $U_{a,b} = (E,H)$ be the uniform (simple)
matroid defined by $E = \{ 1, \ldots,b \}$ and $H = \{ X \subseteq E :
|X| \leq a \}$. We also compute all the minimal and sji
representations of $U_{3,b}$ for $b \geq 5$.

Several other aspects of the theory, intersting enough but not
required for the central core of results, are gathered in Sections 9
and 10.< 

To end this Introduction, we would like to outline why, perhaps,
boolean representable hereditary 
collections should replace matroids. 
\begin{enumerate}
\item
All matroids have boolean representations (first proved in
\cite{IR2}, an alternative proof is supplied here in Theorem
\ref{boom}). The proof follows easily by using the lattice of flats of
the matroid, but also smaller lattices can, in general, provide
representations. Calculating the minimal lattices representing the
matroid is a new important question for matroid theory. Also all the
representations of a matroid are endowed with an operation of join
through stacking, so a representation theory ({\em \`a la} ring theory)
begins. Thus the boolean representation theory, even for matroids, is
much richer than the field matrix representation theory of matroids. 
\item
The classical matroid closure operator extends to boolean
representable hereditary collections (see Section 5).
\item
Strong maps are replaced by $\vee$-maps.
\item
Importantly, {\em a geometry}, like for matroids, is attached to a
boolean representable $(E,H)$, see Subsection \ref{sgeom}. Thus
boolean representable hereditary 
collections are ``not too far'' from matroids, since geometry controls both.
\item
The Tutte idea that ``theorems for graphs can be extended to
matroids'' can be extended to boolean representable hereditary 
collections.
\item
Applications: in near future papers, we plan to consider Coxeter matroids
and Bruhat orders \cite{BGW, LS, Rea}.
The methods here provide a missing ingredient in \cite{BGW}, namely the
definition of boolean representable. 
See future papers.
\end{enumerate}

\section{Preliminaries}

\subsection{Lattices}
\label{sslat}

A poset $(P,\leq)$ is called a lattice if, for all $p,q \in P$, there
exist
$$\begin{array}{l}
(p\vee q) = \min\{ x \in P \mid x \geq p,q \},\\
(p\wedge q) = \max\{ x \in P \mid x \leq p,q \}.
\end{array}$$
For the various aspects of lattice theory, the reader is referred to
\cite{GHK, Gra, RS}.

If only the first (respectively the second) of the above conditions is
satisfied, we talk of a $\vee$-{\em semilattice} (respectively
$\wedge$-{\em semilattice}). We assume also that every
$\vee$-semilattice (respectively $\wedge$-semilattice) has a minimum
(respectively maximum) element.

All the lattices in this paper are finite, and we just write $L$
instead of $(L,\leq)$ most of the time. 
If $L$ is a finite lattice, it is immediate that $L$ has a
maximum (or {\em top}) element, which we denote by $T$, 
% (sometimes by 1), 
and a minimum (or {\em bottom}) element, which we denote by $B$.
% (sometimes by 0). 

We say that $E \subseteq L$ is a $\vee$-{\em generating set} of $L$ if
$L = \{ \vee X \mid X \subseteq E\}$. Note that, whenever convenient,
we may assume that 
$B \notin E$ since $B = \vee \emptyset$. 
Following \cite[Chapters 6,8,9]{RS}, we say that $\p:L \to L'$ is a
$\vee$-map if 
$(\vee X)\p = \vee (X\p)$ for every $X \subseteq E$.
We denote by FL the category of {\em
  finite lattices} together with $\vee$-maps. 

We define also the category FLg by taking objects of the form $(L,E)$,
where $L \in$ FL and $E \subseteq L \setminus \{ B \} $ is a
$\vee$-generating set of $L$. The arrows
$\varphi:(L,E) \to (L',E')$ are $\vee$-maps satisfying $E\p \subseteq
E' \cup \{ B \}$.  

We recall that an element $x$ of a finite lattice $L$ is said to be
{\em strictly meet 
  irreducible} (smi) if $x = (y \wedge z)$ implies $y = x$ or $z =
x$. This is equivalent to saying that $x$ is covered by at most one
element of $L$. Similarly, $x$ is {\em strictly join 
  irreducible} (sji) if $x = (y \vee z)$ implies $y = x$ or $z =
x$. This is equivalent to saying that $x$ covers at most one
element of $L$. See \cite[Subsection 3.3]{IR3} for further details. It
is immediate that the sji elements of $L$ constitute the (unique)
minimum $\vee$-generating set of $L$. 

In the well-known {\em boolean semiring} $\BB = \{ 0,1 \}$ , addition and
multiplication are described respectively by  
$$\begin{tabular}{l|ll}
+&0&1\\
\hline
0&0&1\\
1&1&1\\
\end{tabular}
\hspace{1.5cm}
\begin{tabular}{l|ll}
$\cdot$&0&1\\
\hline
0&0&0\\
1&0&1\\
\end{tabular}$$

We denote by $\M_n(\BB)$ the set of all $n \times n$ matrices with
entries in $\BB$.
The standard boolean matrix representation of a poset $(P,\leq)$ is a $(P
\times P)$-matrix $S(P)$ defined by
$$S(P)_{x,y} = \left\{
\begin{array}{ll}
1&\mbox{ if }x \leq y\\
0&\mbox{ otherwise}
\end{array}
\right.$$
for all $x,y \in P$. If $(L,E) \in \FLg$, then we denote by
$S(L,E)$ the restriction of $S(L)$ to $E \times L$.
For reasons which will become clear later, we
prefer the alternative matrix representation $M(L,E) = ((S(L,E))^c)^t$,
where $M^c$ (for a
boolean matrix $M$) denotes $M$ with 
0 and 1 interchanged, and $M^t$ is just the transposed matrix of
$M$. Thus, for all $\ell \in L$ and $e \in E$, we have 
$$M(L,E)_{\ell,e} = 0 \iff e \leq \ell.$$
Note that $(S(P,\leq))^t = S(P,\geq)$ for every poset
$(P,\leq)$. Moreover, $(M^c)^t = (M^t)^c$ for every boolean matrix
$M$.

The following result collects some of the properties of the boolean
matrices $M(L,E)$. We shall see later that these properties
characterize actually all such matrices.

\bp
\label{propma}
Let $(L,E) \in$ {\rm\FLg} and let $M = M(L,E)$. Then:
\bi
\item[(i)] the rows of $M$ are all distinct;
\item[(ii)] the columns of $M$ are all distinct;
\item[(iii)] $M$ contains a row with all entries equal to 0;
\item[(iv)] $M$ contains a row with all entries equal to 1; 
\item[(v)] the rows of $M$ are closed under addition in $\B^{|E|}$. 
\ei
\ep

\proof
(i) Write $M = (m_{\ell e})$. Since $\ell = \vee\{ e
\in E \mid e \leq \ell\} = \vee\{ e \in E 
\mid m_{\ell,e} = 0 \}$ for every $\ell \in L$, the rows of $M$ are all distinct.

(ii) and (iii) Immediate. 

(iv) Since $B \notin E$, we have $m_{0,e} = 1$ for every $e \in E$. 

(v) Let $k,\ell \in L$. It suffices to show that $m_{k\wedge \ell,e} = m_{k,e} +
m_{\ell,e}$ in $\BB$ for every $e \in E$. This follows from the
equivalence
$$\begin{array}{lll}
m_{k\wedge \ell,e} = 0&\iff&e \leq k\wedge \ell \iff (e \leq k
\mbox{ and } e \leq \ell) \iff ( m_{k,e} = 0 \mbox{ and } m_{\ell,e} =
0)\\
&\iff&m_{k,e} + m_{\ell,e} =
0.
\end{array}$$
\qed

\subsection{Superboolean matrices}
\label{sssm}

Following \cite{IR1,IR2,IR3}, we may view boolean matrices as matrices over
the {\em superboolean semiring} $\SB = \{ 0,1,1^{\nu} \}$, where
addition and multiplication are described respectively by 
$$\begin{tabular}{l|lll}
+&0&1&$1^{\nu}$\\
\hline
0&0&1&$1^{\nu}$\\
1&1&$1^{\nu}$&$1^{\nu}$\\
$1^{\nu}$&$1^{\nu}$&$1^{\nu}$&$1^{\nu}$
\end{tabular}
\hspace{1.5cm}
\begin{tabular}{l|lll}
$\cdot$&0&1&$1^{\nu}$\\
\hline
0&0&0&0\\
1&0&1&$1^{\nu}$\\
$1^{\nu}$&0&$1^{\nu}$&$1^{\nu}$
\end{tabular}$$
We denote by $\M_n(\SB)$ the set of all $n \times n$ matrices with
entries in $\SB$.  Note that $\M_n(\BB)$ is {\em not}
a subsemiring of  $\M_n(\SB)$ since $1 + 1  = 1^{\nu}$.

Next we present definitions of independency and rank appropriate in
the context of superboolean matrices, introduced in \cite{Izh} (see
also \cite{IR1}).  

We say that vectors $C_1,\ldots,C_m \in \SB^n$ are {\em dependent} if
$\lambda_1C_1+ \ldots \lambda_mC_m \in \{ 0,1^{\nu} \}$ for some
$\lambda_1,\ldots, \lambda_m \in \{0,1\}$ not all zero. Otherwise,
they are said to be {\em independent}.   

Let $S_n$ denote the symmetric group on $\hat{n} = \{ 1,\ldots,n \}$. The {\em
  permanent} of a matrix $M = (m_{ij}) \in \M_n(\SB)$ (a positive version of the
determinant) is defined by 
$$\per M = \displaystyle\sum_{\sigma \in S_n} \prod_{i=1}^n
m_{i,i\sigma}.$$
Recall that addition and
multiplication take place in the semiring $\SB$ defined above.

Given $I,J \subseteq \hat{n}$, we denote by $M[I,J]$ the submatrix of
$M$ with entries $m_{ij}$ $(i\in I, j \in J)$. In particular,
$M[\hat{n},j]$ denotes the $j$th column vector of $M$ for each $j \in \hat{n}$.
% The original formulation of the following result involves lower
% triangular matrices, but we can of course interchange upper and lower
% triangular matrices through permutations of rows and columns.
 
\bp
\label{nons}
{\rm \cite[Th. 2.10]{Izh}, \cite[Lemma 3.2]{IR1}}
The following conditions are equivalent for every $M \in \M_n(\SB)$:
\bi
\item[(i)] the column vectors $M[\hat{n},1],\ldots, M[\hat{n},n]$ are
  independent; 
\item[(ii)] ${\rm Per}\,M = 1$;
\item[(iii)] $M$ can be transformed into some lower triangular matrix
  of the form 
\beq
\label{nons1}
\left(
\begin{matrix}
1&&0&&0&&\ldots&&0\\
?&&1&&0&&\ldots&&0\\
?&&?&&1&&\ldots&&0\\
\vdots&&\vdots&&\vdots&&\ddots&&\vdots\\
?&&?&&?&&\ldots&&1
\end{matrix}
\right)
\eeq
by permuting rows and permuting columns independently.
\ei
\ep

A square matrix satisfying the above (equivalent) conditions is said to be
{\em nonsingular}. 

Given (equipotent) $I,J \subseteq \hat{n}$, we say that $I$ is a {\em
  witness} for $J$ in $M$ if $M[I,J]$ is nonsingular. 

\bp
\label{indwit}
{\rm \cite[Th. 3.11]{Izh}}
The following conditions are equivalent for every $m \times n$
superboolean matrix $M$ and every $J \subseteq \hat{n}$: 
\bi
\item[(i)] the column vectors $M[\hat{n},j]$ $(j \in J)$ are independent;
\item[(ii)] $J$ has a witness in $M$.
\ei
\ep

\bp
\label{altrank}
{\rm \cite[Th. 3.11]{Izh}}
The following are equal for a given  $m \times n$
superboolean matrix $M$: 
\bi
\item[(i)] the maximum number of independent column vectors in $M$;
\item[(ii)] the maximum number of independent row vectors in $M$;
\item[(iii)] the maximum size of a subset $J \subseteq \hat{n}$ having
  a witness in $M$;  
\item[(iv)] the maximum size of a nonsingular submatrix of $M$.
\ei
\ep

The {\em rank} of a superboolean matrix $M$, denoted by $\rk M$, is
the number described above. A row of $M$ with $n$ entries
is called an $n$-{\em marker}
if it has one entry 1 and all the remaining entries are 0. The
following remark follows from Proposition \ref{nons}: 

\bc
\label{mark}
{\rm \cite[Cor. 3.4]{IR1}}
If $M \in \M_n(\SB)$ is nonsingular, then it has an $n$-marker.
\ec

\subsection{Hereditary collections}
\label{matroi}

Let $E$ be a set and let $H \subseteq 2^E$. We say that $(E,H)$ is a {\em
  hereditary collection} if $H$ is nonempty and closed under taking
subsets. Hereditary collections are also known as (abstract) {\em
  simplicial complexes} (see \cite{Oxl,Wac}).

We say 
that $X \subseteq E$ is {\em independent} if $X \in H$. A maximal
independent subset of $E$ is called a {\em basis}.
Given $k \in \N$, we call $X \subseteq E$ a $k$-{\em subset} of $E$ if $|X| =
k$. we write $P_k(E) = \{ X \subseteq E : |X| \leq k \}$.

The
hereditary collection $(E,H)$ is said to be a {\em matroid} if the following
condition (the {\em 
  exchange property}) holds:
\bi
\item[(EP)] For all $I,J \in H$ with $|I| = |J|+1$, there exists some
  $i \in I\setminus J$ such that $J \cup \{ i \} \in H$.
\ei
Note that this implies that all basis in a matroid have the same size.

There are many other equivalent definitions of matroid. For details,
the reader is referred to \cite{Oxl}.

\section{Matrices versus lattices}
\label{mvsl}

We establish in this section correspondences between boolean matrices
and $\vee$-generated lattices, adapting results from \cite{IR2}.

Let $E$ be a finite set. Given $\Zcal \subseteq 2^E$,
it is easy to see that
$$\widehat{\Zcal} = \{ \cap S \mid S \subseteq \Zcal \}$$
is the $\wedge$-subsemilattice of $(2^E, \subseteq)$ generated by
$\Zcal$. Note that $\cap \Zcal = \min \widehat{\Zcal}$, and also $E =
\cap \emptyset = \max \widehat{\Zcal}$. In fact,
$(\widehat{\Zcal},\subseteq)$ is 
itself a lattice with the {\em determined join}
$$(P \vee Q) = \cap \{ X \in \Zcal \mid P\cup Q \subseteq X\}.$$
However,  $(\widehat{\Zcal},\subseteq)$ is not in general a sublattice of
$(2^E, \subseteq)$ since the determined join $P \vee Q$ (in
$(\widehat{\Zcal},\subseteq)$) needs 
not to coincide with $P \cup Q$ (see \cite{GHK,RS}). 

Let $M = (m_{ij})$ be an $m \times n$
boolean matrix 
and let $E = \hat{n}$ denote the set of columns of
$M$. We may assume that all the rows of $M$ are distinct. 
For $i \in
\hat{m}$, write $Z_i = \{ j \in \hat{n} \mid m_{ij} = 0\}$ and define
$$\Zcal(M) = \{ Z_1, \ldots, Z_m\} \subseteq 2^E.$$
The {\em lattice of flats} of $M$ is then the lattice $\flats M =
(\widehat{\Zcal(M)},\subseteq)$ (with the determined join). 

Now assume that $M$ has no zero columns. This is equivalent to saying
that $\emptyset \in \flats M$. 
For $j \in \hat{n}$, define also $Y_j = \cap\{ Z_i \mid m_{ij} = 0 \}$
and let
$$\Y(M) = \{ Y_1, \ldots, Y_n \} \subseteq \flats M.$$
Note that $Y_j = \cap\{ Z_i \mid j \in Z_i \}$ and so $j \in Y_j$ for every $j$.

\bl
\label{gener}
Let $M = (m_{ij})$ be an $m \times n$
boolean matrix without zero columns. 
Then ({\rm Fl}$\, M, \Y(M)) \in$
{\rm FLg}. 
\el

\proof
First note that $Y_j$ can never be the bottom element $\emptyset$
since $j \in Y_j$. Hence it suffices to show that 
\beq
\label{gener1}
Z_{i_1} \cap \ldots \cap Z_{i_k} = \vee\{ Y_j \mid j \in Z_{i_1} \cap
\ldots \cap Z_{i_k} \}. 
\eeq
holds for all $i_1, \ldots, i_k \in \hat{m}$.

Indeed, take $j \in Z_{i_1} \cap \ldots \cap Z_{i_k}$. On the one
hand, we have $m_{i_1j}
= \ldots = m_{i_kj} = 0$ and so $Y_j \subseteq Z_{i_1} \cap \ldots
\cap Z_{i_k}$. Thus $\vee\{ Y_j \mid j \in Z_{i_1} \cap
\ldots \cap Z_{i_k} \} \subseteq Z_{i_1} \cap \ldots \cap Z_{i_k}$.

On the other hand, since $j \in
Y_j$ for every $j$, we get 
$$Z_{i_1} \cap \ldots \cap Z_{i_k} \subseteq 
\cup\{ Y_j \mid j \in Z_{i_1} \cap
\ldots \cap Z_{i_k} \} \subseteq
\vee\{ Y_j \mid j \in Z_{i_1} \cap
\ldots \cap Z_{i_k} \}$$
and so (\ref{gener1}) holds as required.
\qed

Hence $M \mapsto (\flats M,\Y(M))$ defines an operator from the set of
boolean matrices without zero columns
into $\FLg$. 

We can relate
this operator with the 
matrix representation defined in Subsection \ref{sslat}.
Given a lattice $L$ and $\ell \in L$, let $\ell\dw = \{ x \in L \mid x
\leq \ell \}$. We start with
the following remark:

\bl
\label{zl}
Let $(L,E) \in$ {\rm FLg} and let $M = M(L,E) = (m_{\ell e})$. Then
$Z_{\ell} = \ell\dw \cap E$ for every $\ell \in L$.
\el

\proof
Indeed,
$$Z_{\ell} = \{ e \in E \mid  m_{\ell e} = 0 \} = \{ e \in E \mid e
\leq \ell \} = \ell\dw \cap E.$$
\qed

Next we establish that the lattice of flats of the matrix
representation of a lattice gives back the original lattice:

\bp
\label{ltom}
Let $(L,E) \in$ {\rm FLg} and let $M = M(L,E) = (m_{\ell e})$. Then
({\rm Fl}$\, M,\Y(M)) 
\cong (L,E)$.
\ep

\proof
Let $\p:L \to \flats M$ be defined by $\ell\p = Z_{\ell}$. 
Since $E$ is a $\vee$-generating set of $L$, it follows from Lemma
\ref{zl} that
\beq
\label{ltom1}
Z_k \subseteq Z_{\ell} \iff k \leq \ell
\eeq
holds for all
$k,\ell \in L$.
Thus $\p$ is a
poset embedding. On the other hand, $e \leq (k \wedge \ell)$ if and
only if $e \leq k$ and $e \leq \ell$, hence
$Z_k \cap Z_{\ell} = Z_{k \wedge \ell}$ for all $k,\ell \in L$. This
immediately generalizes to
\beq
\label{ltom2}
Z_{\ell_1} \cap \ldots \cap Z_{\ell_n} = Z_{\ell_1 \wedge \ldots
  \wedge \ell_n}
\eeq
for all $\ell_1, \ldots, \ell_n \in L$, hence
$\p$ is surjective. Thus $\p$ is an isomorphism of posets and
therefore of lattices.

It remains to show that $\Y(M) = \{ Z_e \mid e \in E\}$. It suffices
to prove that $Y_e = Z_e$ for every $e \in E$. Indeed,
$$Y_e = \cap\{ Z_{\ell} \mid m_{\ell e} = 0 \} = \cap\{ Z_{\ell} \mid
e \leq \ell \} = Z_e$$ and we are done.
\qed 

We shall refer to $\flats (L,E) = \flats M(L,E)$ as the {\em lattice of flats}
of $(L,E) \in \FLg$. 

Given matrices $M$ and $M'$, we say that $M$ and $M'$ are {\em
  congruent} and write $M \cong M'$ if $M'$ can be
obtained from $M$ by 
permuting rows and permuting columns (independently!).
Given a boolean matrix $M$ without zero columns, we write
$M^{\nu} = M(\flats M, \Y(M))$. 
In view of Proposition \ref{propma}, it is not true that $M^{\nu} \cong M$
in general. However, we can get a correspondence by focusing our
attention on the set $\M$ of all boolean matrices satisfying
conditions (i)-(v) of Proposition \ref{propma}: 

\bp
\label{mtol}
Let $M \in \M$. Then $M^{\nu} \cong M$.
\ep

\proof
Assume that $M = (m_{ij})$ is an $m \times n$ matrix in $\M$. Since
$M$ satisfies conditions (iii) and (v) of Proposition \ref{propma}, we
have $\flats M = \{ Z_1, \ldots, Z_m\}$.  Since
$M$ satisfies condition (i) of Proposition \ref{propma}, these
elements are all distinct. 
Note that, since
$M$ satisfies condition (iv) of Proposition \ref{propma}, has no 
zero columns and so $(\flats M, \Y(M)) \in \FLg$ by Lemma \ref{gener}. 

Therefore $M^{\nu} = (m'_{Z_iY_j})$ is also a boolean matrix with $m$
rows. To complete the proof, it suffices to show that  
$m'_{Z_iY_j} = m_{ij}$ for all $i \in \hat{m}$ and $j \in \hat{n}$. In
view of condition (ii) of Proposition \ref{propma}, $M^{\nu}$ is then
an $m \times n$ matrix in $\M$ and we shall be done. 

Indeed, $m'_{Z_iY_j} = 0$ if and only if $Y_j \subseteq Z_i$. Since $j
\in Y_j$, this implies $j \in Z_i$. Conversely, $j \in Z_i$ implies
$Y_j \subseteq Z_i$ and so
$$m'_{Z_iY_j} = 0 \iff Y_j \subseteq Z_i \iff j \in Z_i \iff m_{ij} =
0.$$
Therefore $m'_{Z_iY_j} = m_{ij}$ and so $M^{\nu} \cong M$.
\qed

Now it is easy to establish a correspondence between the set
$\FLg/\!\cong$ of isomorphism classes of $\FLg$ and the set
$\M/\!\cong$ of congruence classes of $\M$:

\bc
\label{corre}
The mappings $\M \to$ {\rm FLg}$: M \mapsto (${\rm Fl}$\, M, \Y(M))$ and
{\rm FLg} $\to \M: (L,E) \mapsto M(L,E)$ induce mutually inverse bijections
between $\M/\!\cong$ and {\rm FLg}$/\!\cong$.
\ec

\proof
It follows easily from the definitions that the above operators induce
mappings between $\M/\!\cong$ and $\FLg/\!\cong$. These mappings are
mutually inverse by Propositions \ref{ltom} and \ref{mtol}.
\qed

\be
Let $M$ be the matrix
$$\left(
\begin{matrix}
1&&0&&1&&0&&1\\
1&&0&&0&&1&&1\\
1&&1&&0&&0&&0\\
\end{matrix}
\right)$$
Omitting brackets and commas, and identifying the elements
$Y_1,\ldots, Y_5$ of $\Y(M)$,
the lattice of flats {\rm Fl}$\, M$ can
be represented as
$$
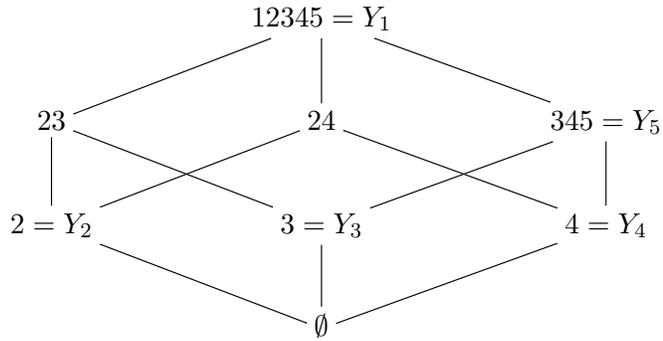
\xymatrix{
&& {12345 = Y_1} \ar@{-}[dll] \ar@{-}[d] \ar@{-}[drr] && \\
23 \ar@{-}[drr] \ar@{-}[d] && 24  \ar@{-}[dll] \ar@{-}[drr] && {345 =
  Y_5} \ar@{-}[dll] \ar@{-}[d] \\
{2 = Y_2} \ar@{-}[drr] && {3 = Y_3} \ar@{-}[d] && {4 = Y_4} \ar@{-}[dll] \\
&& \emptyset &&
}$$
% Note that $\Y(M)$ consists here of 12345, 2, 3, 4 and 345. 
Finally, $M^{\nu}$
is the matrix 
$$\left(
\begin{matrix}
1&&0&&1&&0&&1\\
1&&0&&0&&1&&1\\
1&&1&&0&&0&&0\\
1&&0&&1&&1&&1\\
1&&1&&0&&1&&1\\
1&&1&&1&&0&&1\\
0&&0&&0&&0&&0\\
1&&1&&1&&1&&1
\end{matrix}
\right)$$
\ee

The above example illustrates a simple remark: if all the columns of
$M$ are distinct and nonzero, if all its rows are distinct, then
$M^{\nu}$ can be obtained from $M$ by adding a zero row 
and any new rows obtained by adding the rows of $M$ in $\B^{|E|}$.

\section{c-independence and subset closure}
\label{cisl}

From now on, if we mention a lattice $L$ without specifying a
$\vee$-generating set $E$, it is assumed that $E = L \setminus\{
B\}$. We also assume that $B \neq T$ in $L$. Its {\em height}, denoted
by $\het L$, is the maximum length of a chain in $\L$. In view of the
matrix representation $M(L)$, we say that $\ell_1,\ldots, \ell_k \in
L$ are {\em c-independent} if the the column vectors of $M(L)$
corresponding to $\ell_1,\ldots, \ell_k$ are independent (over $\SB$).
Note that, if $(L,E) \in
\FLg$ and $X
\subseteq E$, this is equivalent to saying that the column vectors of
$M(L,E)$ defined by the 
elements of $X$ are independent (over $\SB$).

The next result
generalizes Theorem 3.6 of \cite{IR2}:

\bp
\label{rtoh}
Let $(L,E) \in$ {\rm FLg}. Then {\rm rk}$\, M(L,E) =$ {\rm rk}$\,
M(L) =$ {\rm ht}$\,L$. 
\ep

\proof
The second equality follows from \cite[Theorem 3.6]{IR2} because the
omitted column corresponding to 0 contains only zeros and is therefore
irrelevant to the computation of the c-rank.

Since $M(L,E)$ is a submatrix of $M(L)$, we have $\rk M(L,E) \leq \rk
M(L)$. To prove the opposite inequality, it suffices to show that 
\beq
\label{rtoh1}
\{ x\vee y\} \cup Z \mbox{ c-independent} \hspace{.7cm} \Rw
\hspace{.7cm} \{ x \} \cup Z \mbox{ or } \{ y\} \cup Z
 \mbox{ c-independent}.
\eeq
Indeed, if (\ref{rtoh1}) holds, we can start with a c-independent
subset $\{ \ell_1,\ldots, \ell_k\}$ of $L$ and by successive
application of (\ref{rtoh1}) replace it by a c-independent
subset of $E$ with the same number of elements.

Assume that $\{ x\vee y, z_1,\ldots,z_k\}$ is c-independent. By
Proposition \ref{indwit}, and permuting columns and rows if necessary,
we may assume that $M(L)$ has a submatrix of
the form (\ref{nons1}), where the columns correspond to
$z_1,\ldots,z_j,x\vee y,z_{j+1},\ldots,z_k$ $(j \in \{ 0, \ldots,
k\})$ and the rows correspond to $\ell_1,\ldots, \ell_j,\ell,
\ell_{j+1},\ldots,\ell_k$. Hence $(x\vee y) \leq \ell_1,\ldots, \ell_j$
and $(x\vee y) \not\leq \ell$. It follows that $x \not\leq \ell$ or $y
\not\leq \ell$. On the other hand, we get $x,y \leq (x\vee y) \leq
\ell_1,\ldots, \ell_j$ and so $\ell_1,\ldots, \ell_j,\ell,
\ell_{j+1},\ldots,\ell_k$ is a witness for at least one of the sets
$\{ x, z_1,\ldots,z_k\}$, $\{ y, z_1,\ldots,z_k\}$. Therefore
(\ref{rtoh1}) holds as required.
\qed

We can use the lattice of flats $\flats (L,E)$ to
define a closure operator (see Subsection \ref{101} of the
Appendix) in the lattice $(2^E,\subseteq)$: given $X
\subseteq E$, let 
$$\clos_LX = \cap\{ Z \in \flats (L,E) \mid X \subseteq Z\}.$$ 
% If the context is clear, we may use the simplified notation $\oo{X}
% = \clos_LX$.

Recalling the notation from Section \ref{mvsl} and Lemma \ref{zl}, it
is easy to see that  
\beq
\label{clz}
\clos_L X = Z_{\vee X} = (\vee X)\dw \cap E.
\eeq
Indeed, we have $X \subseteq
Z_{\vee X} \in \flats (L,E)$, and the equivalence
$$X \subseteq Z_{\ell} \hspace{.3cm}\iff \hspace{.3cm} \forall x\in X
\; x \leq \ell \hspace{.3cm} \iff\hspace{.3cm} \vee X
\leq \ell \hspace{.3cm} \iff \hspace{.3cm}  Z_{\vee X} \subseteq
Z_{\ell}$$ follows from (\ref{ltom1}), hence (\ref{clz}) holds. 

Note that $X \in \flats
(L,E)$ if and only if $\clos_L X = X$, and $\clos_L$ is indeed a
closure operator in the lattice $(2^E,\subseteq)$. 

We say that $X = \{ x_1,\ldots,x_k \} \subseteq E$ is a {\em 
  transversal}  of the partition of the successive differences for the
chain $Y_0 \supset \ldots \supset Y_k$ in $\flats(L,E)$
if $x_i \in Y_{i-1} \setminus Y_i$ for $i =
1,\ldots,k$. A subset of a transversal is a {\em 
partial transversal}.

By adapting the proofs of \cite[Lemmas 3.4 and 3.5]{IR2}, we can prove
the following:

\bp
\label{cindch}
Let $(L,E) \in$ {\rm FLg} and $X \subseteq
E$. Then the following conditions are
equivalent: 
\bi
\item[(i)] $X$ is c-independent;
\item[(ii)] $X$ admits an enumeration $x_1,\ldots , x_k$ such that
\beq
\label{cindch1}
(x_1 \vee \ldots \vee x_k) > (x_2 \vee \ldots \vee x_k) > \ldots >
(x_{k-1} \vee x_k) > x_k;
\eeq  
\item[(iii)] $X$ admits an enumeration $x_1,\ldots , x_k$ such that
$${\rm Cl}_L(x_1,\ldots, x_k) \supset {\rm Cl}_L(x_2,\ldots, x_k) \supset \ldots
\supset {\rm Cl}_L(x_k);$$
\item[(iv)] $X$ admits an enumeration $x_1,\ldots , x_k$ such that
$$x_i \notin {\rm Cl}_L(x_{i+1},\ldots, x_k)\hspace{1cm}(i = 1,\ldots,k-1);$$
\item[(v)] $X$ is a transversal of the partition of successive
  differences for some chain of {\rm Fl}$(L,E)$; 
\item[(vi)] $X$ is a partial transversal of the partition of successive
  differences for some maximal chain of {\rm Fl}$(L,E)$.
\ei
\ep

\proof
(i) $\Rw$ (ii). If $X$ is c-independent, then 
$M(L)$ admits a submatrix of the form (\ref{nons1}), with the columns
labelled, say, by $x_1, \ldots, x_k$. It is a simple exercise to show
that (\ref{cindch1}) holds.

(ii) $\Rw$ (iii). By (\ref{ltom1}) and (\ref{clz}).

(iii) $\Rw$ (iv). If $x_i \in \clos_L(x_{i+1},\ldots, x_k)$, then
$\clos_L(x_{i},\ldots, x_k) = \clos_L(x_{i+1},\ldots, x_k)$.

(iv) $\Rw$ (ii). Clearly, $(x_{i} \vee \ldots \vee x_k) \geq (x_{i+1}
\vee \ldots \vee x_k)$, and equality would imply $\clos_L(x_{i},\ldots,
x_k) = \clos_L(x_{i+1},\ldots, x_k)$ by (\ref{clz}).

(ii) $\Rw$ (i). If (\ref{cindch1}) holds, we build a nonsingular
submatrix of $M(L)$ of the form (\ref{nons1}) by taking rows labelled
by $\ell_1, \ldots, \ell_k \in L$, where $\ell_i = (x_{i+1} \vee \ldots
\vee x_k)$ $(i = 1,\ldots,k)$.

(iv) $\iff$ (v) $\iff$ (vi). Immediate.
\qed

It is easy to characterize c-independence for small numbers of vectors:

\bp
\label{onetwo}
Let $(L,E) \in$ {\rm FLg} and let $X \subseteq E$.
\bi
\item[(i)] If $|X| \leq 2$, then $X$ is c-independent.
\item[(ii)] If $X$ is c-independent and $\vee X < 1$, then $X \cup \{
  p \}$ is c-independent for some $p \in E \setminus X$.
\ei
\ep

\proof
(i) The case $|X| \leq 1$ is immediate, hence we may assume that $X = \{
x_1,x_2 \}$ and $x_1 \not\leq x_2$. Then $(x_1 \vee x_2) > x_2$ and so
$X$ is c-independent by Proposition \ref{cindch}.

(ii) Since $1 = \vee E$, there exists some $p \in E$ such
that $((\vee X) \vee p) > (\vee X)$ and so $X \cup \{ e \}$ is a
c-independent subset of $E$ 
by Proposition \ref{cindch} (using the characterization in (ii)).
\qed

We discuss now the c-independence of 3-subsets.

\bp
\label{threeind}
Let $(L,E) \in$ {\rm FLg} and let $X$ be a 3-subset of $E$. Then the
following conditions are equivalent: 
\bi
\item[(i)] $X$ is c-independent;
\item[(ii)] $X$ admits an enumeration $x_1, x_2, x_3$ such that
$$(x_1 \vee x_2 \vee x_3) > (x_2 \vee x_3) > x_3;$$
\item[(iii)] $X$ admits an enumeration $x_1, x_2, x_3$ such that $x_1
  \notin$ {\rm Cl}$_L(x_2,x_3)$;
\item[(iv)] $X$ is contained in some c-independent 4-subset of $E$ or
  there exists some $x \in X$ such that $\vee (X\setminus\{ x\}) < \vee
  X = 1$.
\ei
\ep

\proof
(i) $\iff$ (ii). By Proposition \ref{cindch}.

(i) $\Rw$ (iii). By Proposition \ref{cindch}.

(iii) $\Rw$ (i). Since $x_2 \neq x_3$, we may assume that $x_2 \not\leq
x_3$. Hence $x_2 \notin Z_{x_3}$ and so $x_2 \notin
\clos_L(x_3)$ in view of (\ref{clz}). By Proposition
\ref{cindch}, $X$ is c-independent. 

(ii) $\Rw$ (iv). Assume that $(x_1 \vee x_2 \vee x_3) > (x_2 \vee x_3)
> x_3$ with $X = \{ x_1,x_2,x_3\}$. The case $(x_1 \vee x_2 \vee x_3)
= 1$ is immediate and the case $(x_1 \vee x_2 \vee x_3) < 1$ follows
from Proposition \ref{onetwo}(ii).

(iv) $\Rw$ (i). Clearly, c-independent sets are closed under
inclusion, hence we may assume that $X = \{ x,y,z\}$ and $(y \vee z) < \vee
  X = 1$. Since we may assume
that $y \not\leq z$, then $z < (z \vee y)$ and so $X$ is c-independent
by Proposition \ref{cindch}. 
\qed

\section{Representation of hereditary collections}

Let $(E,H)$ be a hereditary collection. We say that $X \subseteq E$ is
{\em closed} (or a {\em flat}) if
$$\forall I \in H \cap 2^X\; \forall p \in E \setminus X \hspace{.7cm}
I \cup \{ p \} \in H.$$
The set of all flats of $(E,H)$ is denoted by $\flats (E,H)$. 

An alternative characterization is provided through the notion of {\em
  circuit}: $C \subseteq E$ is said to be a circuit of $(E,H)$ if
$C \notin H$ but all proper subsets of $C$ are in $H$. 

\bp
\label{flatcirc}
Let $(E,H)$ be a hereditary collection and let $X \subseteq E$. Then
the following conditions are equivalent:
\bi
\item[(i)] $X$ is closed;
\item[(ii)] if $p \in C \subseteq X \cup \{ p \}$ for some circuit
  $C$, then $p \in X$.
\ei
\ep

\proof
(i) $\Rw$ (ii). Suppose that there exist a circuit $C$ and $p \in C
\subseteq X \cup \{ p \}$ such that $p \notin X$. Then $C = I \cup \{ p 
\}$ for some $I \subseteq X$. It follows that $I \in H \cap 2^X$
and $p \in E \setminus X$, however $I \cup \{ p \} \notin
H$. Therefore $X$ is not closed.

(ii) $\Rw$ (i). Suppose that $X$ is not closed. Then there exist $I
\in H \cap 2^X$ and $p \in E \setminus X$ such that $I \cup \{ p \} \notin
H$. Let $I_0 \subseteq I$ be minimal for the property $I_0 \cup \{ p \} \notin
H$. Since $I_0 \in H$ due to $I_0 \subseteq I \in H$, it follows that
$I_0 \cup \{ p \}$ is a circuit by minimality of $I_0$. Thus condition
(ii) fails for $C = 
I_0 \cup \{ p \}$ and we are done.
\qed

Note that condition (ii) is the standard characterization of flats for
matroids. 

The following result summarizes some straightforward properties of
$\flats (E,H)$. We say that $(E,H)$ is {\em simple} if $P_2(E) \subseteq H$. A
1-subset of $E$ is also called a {\em point}.  

\bp
\label{trivia}
Let $(E,H)$ be a hereditary collection.
\bi
\item[(i)] If $Y \subseteq$ {\rm Fl}$\,(E,H)$, then $\cap Y \in$ {\rm
    Fl}$\,(E,H)$. 
\item[(ii)] If $P_k(E) \subseteq H$ with $k
  \geq 1$, then $P_{k-1}(E) \subseteq$ {\rm Fl}$\,(E,H)$. 
\item[(iii)] If $(E,H)$ is simple, then the points of $E$ are closed.
\ei
\ep

\proof
(i) We have $E = \cap \emptyset \in \flats(E,H)$ trivially, hence it
suffices to show that $X_1,X_2 \in \flats(E,H)$ implies $X_1\cap X_2 \in
\flats(E,H)$.  

Let $I \in H \cap 2^{X_1\cap X_2}$ and $p \in E \setminus (X_1\cap
X_2)$. Then $p \notin X_1$ or $p \notin X_2$. Hence $X_1 \cap X_2
\subseteq X_j$ and $p \in E \setminus X_j$ for some $j \in \hat{2}$. Since
$X_j \in \flats(E,H)$, we get $I \cup \{ p \} \in H$ and so $X_1\cap X_2 \in
\flats(E,H)$.

(ii) Immediate.

(iii) By part (ii).
\qed

Similarly to Section \ref{cisl}, we  can use the lattice of flats
$\flats (E,H)$ to define a closure operator in $(2^E,\subseteq)$:
given $X \subseteq E$, let 
$$\cl X = \cap\{ Z \in \flats (E,H) \mid X \subseteq Z\}.$$ 
Note that $X \subseteq E$ is closed if and only if $\cl X = X$.
We can also make the following remark:

\bp
\label{newtrivia}
Let $(E,H)$ be a hereditary collection and let $X \subseteq E$ be a
basis. Then {\rm Cl}$\, X = E$.
\ep

\proof
Suppose that $p \in E \setminus \cl X$. Since $X \subseteq \cl X$
and $\cl X$ is closed, we get $X \cup \{ p \} \in H$, contradicting
$X$ being a basis. Thus $\cl X = E$.
\qed

In the matroid case (see Proposition \ref{flatcirc}), the concept of
circuit allows a more constructive perspective of the closure:

\bp
\label{circuits}
{\rm \cite[Proposition 1.4.10(ii)]{Oxl}}
Let $(E,H)$ be a matroid and let $X \subseteq E$. Then
\beq
\label{circuits1}
{\rm Cl}\, X = X \cup \{ p \in E \setminus X \mid p \in C \subseteq
X \cup \{ p \} \mbox{ for some circuit } C\}.
\eeq
\ep

In the general case, it is still true that all such elements $p$ must
be in the closure, but they may not be sufficient. We may then have to
iterate the construction. Eventually, iteration will gives us all the
elements of $\cl X$.

The subsets of independent column vectors of a given (superboolean)
matrix, which include the empty subset and are closed for subsets,
constitute an important example of a hereditary collection. 
In fact, {\em every} hereditary subset is given this way, see \cite{IR1}.
On the
other hand, a {\em boolean representation} of a hereditary collection
$(E,H)$ is a boolean matrix $M$ with column space $E$ such that a subset
$X \subseteq E$ of column vectors of $M$ is independent (over $\SB$) if and
only if $X \in H$. Obviously, we can always assume that the rows in
such a matrix are distinct: the representation is then said to be
{\em reduced}. Note also that by permuting rows in a reduced
representation of $(E,H)$ we get an alternative reduced
representation of $(E,H)$. The number of rows in a boolean
representation $M$ of $(E,H)$ is said to be its {\em degree} and is
denoted by $\deg M$. We
denote by $\mindeg(E,H)$ the minimum degree of a boolean
representation of $(E,H)$.

We remark also that if $(E,H)$ admits a boolean  representation, then
$(E,H)$ satisfies (PR) (see \cite[Theorem 5.3]{IR1}). 

Given an $R \times E$ boolean matrix $M = (m_{re})$ and $r \in R$, we
recall the notation $Z_r = \{ e \in E \mid m_{re} = 0 \}$ introduced
in Section \ref{mvsl}.

\bp
\label{zeroes}
Let $(E,H)$ be a hereditary collection and let $M$ be an $R \times E$
boolean matrix. If $M$ is a boolean representation of $(E,H)$, then
$Z_r \in$ {\rm Fl}$\,(E,H)$ for every $r \in R$.
\ep

\proof
Let $r \in R$ and $J \subseteq Z_r$. Suppose that $J \in H$ and $p
\in E\setminus Z_r$. Since $M = (m_{re})$ is a boolean representation of $(E,H)$,
then the column vectors $M[R,j]$ $(j \in J)$ are independent and so
there exists some $I \subseteq R$ such that $M[I,J]$ is of the form
(\ref{nons1}), for a suitable ordering of $I$ and $J$. Since $J
\subseteq Z_r$, the row vector $M[r,J]$ contains only zeros. On the
other hand, since $p \notin Z_r$, we have $m_{rp} = 1$ and so $M[I
\cup \{ r\}, J \cup \{ p\}]$ is also of the form
(\ref{nons1}) and therefore nonsingular. Thus $J \cup \{ p\}$ defines
an independent set of columns of $M$. Since $M$ is a boolean
representation of $(E,H)$, it follows that $J \cup \{ p\} \in H$ and
so $Z_r \in \flats(E,H)$.
\qed

Let $(E,H)$ be a hereditary collection.
In view of Proposition \ref{trivia}(i), we can view
$(\flats(E,H),\subseteq)$ as a lattice with $(X \wedge Y) = X\cap Y$ and
the determined join $(X \vee Y) = \cap\{ Z \in \flats(E,H) \mid X\cup Y
\subseteq Z \} = 
\clos(X\cup Y)$. If
$(E,H)$ is simple, by identifying $e \in E$ with
$\{ e \}$ and $E$ with $\{ \{e\} \mid e \in E\}$
we may write $(\flats(E,H),E) \in \FLg$. Indeed, for every $X \in
\flats(E,H)$, we have $X = \vee\{ e \mid e \in X\}$. Recalling the
representations of $\vee$-generated lattices from Subsection \ref{sslat}, we can
prove the following: 

\bl
\label{indlh}
Let $(E,H)$ be a simple hereditary collection and let $X \subseteq E$
be c-independent with respect to $M(${\rm Fl}$\,(E,H),E)$.
Then $X \in H$.
\el

\proof
We use induction on $|X|$. Since $(E,H)$ is simple, the case $|X| \leq
1$ is trivial. Hence we assume
that $|X| = m > 1$ and the claim holds for $|X| = m-1$.

By Proposition \ref{indwit}, $X$ has a witness $P$ in $M = M(\flats
(E,H),E)$. We may assume that $X = \{ e_1,\ldots,e_m\}$, $P = \{
P_1,\ldots,P_m\}$ and $M[P,X]$
is of the form (\ref{nons1}), with the rows (respectively the columns)
ordered by $P_1,\ldots,P_m$ (respectively $e_1,\ldots,e_m$). The first
row yields $e_1 \notin P_1$ and $e_2,\ldots,e_m \in P_1$.

Now, since $e_2,\ldots,e_m$ is c-independent, it
follows from the induction hypothesis that $\{ e_2,\ldots,e_m\} \in H$.
Together with $\{ e_2,\ldots,e_m\} \subseteq P_1 \in \flats
(E,H)$ and $e_1 \notin P_1$, this yields $X = \{ e_2,\ldots,e_m\} \cup
\{ e_1 \} \in H$ as required.
\qed

Given matrices $M_1$ and $M_2$ with the same number of columns, we
define $M_1\oplus_b M_2$ to be the matrix obtained by concatenating the
matrices $M_1$ and $M_2$ by
$$\left(
\begin{matrix}
M_1\\
M_2
\end{matrix}
\right)$$
and removing repeated rows (leaving only the first occurrence from top
to bottom, say). We refer to this matrix as $M_1$ {\em stacked over} $M_2$.

\bp
\label{stack}
Let $(E,H)$ be a simple hereditary collection. 
\bi
\item[(i)] If $M_1$ and $M_2$ are
reduced boolean representations of $(E,H)$, so is $M_1\oplus_b M_2$.
\item[(ii)] If $M$ is a reduced boolean representation of $(E,H)$ and
  we add/erase a row which is the sum of other rows in $\BB^{|E|}$, we
    get a matrix $M'$ which is also a reduced boolean representation of $(E,H)$.
\ei
\ep

\proof
(i)
Since $M_1$ and $M_2$ have both space of columns $E$, the matrix
$M = M_1\oplus_b M_2$ is well-defined and has no repeated rows by
definition. Let $R$ be the row space of $M$ and let $X
\subseteq E$. We show that 
\beq
\label{stack1}
X \mbox{ is c-independent with respect to } M \; \iff \; X \in H
\eeq
by induction on $|X|$. The case $|X| = 0$ being trivial, assume that
$|X| > 0$ and (\ref{stack1}) holds for smaller values of $|X|$. 

Suppose that $X$ is c-independent with respect to $M$. By permuting
rows of $M_1\oplus_b M_2$ if necessary, and 
using the appropriate ordering of $E$, we may say that there exists
some $P \subseteq R$ such that $B[P,X]$ is of the form
(\ref{nons1}). Let $p_1$ (respectively $x_1$) denote the first element
of $P$ (respectively $X$) for these orderings, so $M[P\setminus\{
p_1\},X\setminus\{ x_1\}]$ is the submatrix of $M[P,X]$ obtained by
deleting the first row and the first column.
Since reduced boolean representations are closed under
permuting rows, we may assume without loss of generality that the row
$M[p_1,E]$ came from the matrix $M_1$. On the other hand, since the
column vectors $M[R,x]\; (x \in X\setminus\{ x_1\})$ are independent,
it follows from the induction hypothesis that $X\setminus\{ x_1\} \in
H$ and so (since $M_1$ is a boolean representation of $(E,H)$)
$X\setminus\{ x_1\}$ is c-independent with respect to $M_1$. Hence
$M_1$ has a singular submatrix of the form 
$M_1[P',X\setminus\{ x_1\}]$. Now $M_1[P'\cup \{ p_1\},X]$ is still a
nonsingular matrix because the unique nonzero entry in the row
$M_1[p_1,X]$ is $M_1[p_1,x_1]$. Hence $X$ is c-independent with
respect to $M_1$ and so $X \in H$. 

Conversely, if $X \in H$, then  $X$ is c-independent with
respect to $M_1$ and so  $X$ is c-independent with
respect to $M$ as well. Thus (\ref{stack1}) holds and so $M$ is a reduced
boolean representation of $(E,H)$ as claimed.

(ii) The claim is obvious when we add a row, so consider the case when
a row of the described form is erased.
It is easy to see that if a $k$-marker $u$ is the sum of some vectors
in $\BB^k$, then one of them is equal to $u$. Therefore, if the sum
row occurs in some nonsingular submatrix of $M$, we can always replace
it by one of the summand rows.
\qed

Proposition \ref{stack}(i) immediately implies that if $(E,H)$ admits
a reduced boolean 
representation, then
there exists a unique maximal one. The main theorem of this section
provides a more concrete characterization:

\bt
\label{rephc}
Let $(E,H)$ be a simple hereditary collection. Then the following
conditions are equivalent:
\bi
\item[(i)] $(E,H)$ has a boolean representation;
\item[(ii)] $M(${\rm Fl}$\,(E,H),E)$ is a reduced boolean
  representation of $(E,H)$. 
\ei
Moreover, in this case any other reduced boolean
  representation of $(E,H)$ is congruent to a submatrix of $M(${\rm
    Fl}$\,(E,H),E)$. 
\et

\proof
(i) $\Rw$ (ii).
Write $M = M(\flats(E,H),E)$.
Suppose that $(E,H)$ has a boolean representation
$N = (n_{re})$. Then we may assume that the $R \times E$ matrix $N$ is
reduced. By Proposition \ref{zeroes}, we have 
$Z_r \in \flats(E,H)$ for every $r \in R$. For every $e \in E$, we
have $$n_{re} = 0 \iff e \in Z_r \iff \{ e \} \subseteq Z_r \iff
M[Z_r,e] = 0,$$
hence $N$ is (up to permutation of rows) a submatrix of $M$.

We claim that $M$ is also a boolean representation of $(E,H)$. Indeed,
let $X \subseteq E$. If $X \in H$, then $X$ is c-independent with
respect to $N$ since $N$ is a boolean representation of
$(E,H)$, hence $X$ is c-independent with
respect to $M$ since $N$ is a submatrix of $M$. The
converse implication follows from Lemma \ref{indlh}, hence  $M$ is a
boolean representation of $(E,H)$. Naturally, every representation
arising from a $\vee$-generated lattice is reduced.

(ii) $\Rw$ (i). Trivial.
\qed

\section{The lattice of lattice representations}
\label{tlor}

We  define a quasi-order on $\FLg$ by 
$$(L,E) \geq (L',E) \hspace{.7cm}\mbox{if there exists some $\vee$-map
  $\p:L \to L'$ such that $\p|_E = id$}.$$
Note that such $\p$ is necessarily onto: if $\ell' \in L'$, we may
write $\ell' = (e_1 \vee \ldots \vee e_k)$ in $L'$ for some $e_1,
\ldots, e_k \in E$, hence 
$$\ell' = (e_1 \vee \ldots \vee e_k) = (e_1\p \vee \ldots \vee e_k\p)
= (e_1 \vee \ldots \vee e_k)\p \in L\p.$$
Recall that $\flats(L,E) = \{ Z_{\ell} \mid \ell \in L\}$, and $Z_{\ell} =
\ell\dw \cap E$ by Lemma \ref{zl}. 

\bp
\label{repmap}
Let $(L,E),(L',E) \in$ {\rm FLg}. Then
$$(L',E)
\leq (L,E) \iff {\rm Fl}\,(L',E) \subseteq {\rm Fl}\,(L,E).$$
\ep

\proof
Assume first that $(L',E)
\leq (L,E)$. Then there exists some $\vee$-map $\p:L \to L'$ such
that $\p|_E = id$. We show that $\flats(L',E) \subseteq \flats(L,E)$.
Indeed, we claim that
\beq
\label{repmap1}
Z_{\ell'} = Z_{\vee (\ell'\p\inv)}
\eeq
holds for every $\ell' \in L'$. Let $e \in Z_{\ell'}$. Then $e \leq
\ell'$. Write $\ell = \vee (\ell'\p\inv)$. Since $\p$ is onto and a
$\vee$-map, we have $\ell\p = \ell'$. Moreover, $(\ell \vee e)\p =
(\ell\p \vee e\p) = (\ell' \vee e) = \ell'$, hence $(\ell \vee e) \in
\ell'\p\inv$ and so $(\ell \vee e) \leq \max(\ell'\p\inv) =
\ell$. Thus $e \leq \ell$ and so $Z_{\ell'} \subseteq Z_{\vee (\ell'\p\inv)}$.

Conversely, assume that $e \in Z_{\vee (\ell'\p\inv)}$. Then $e \leq
\vee (\ell'\p\inv)$ and so $e = e\p \leq \ell'$. Hence $Z_{\vee
  (\ell'\p\inv)} \subseteq Z_{\ell'}$ and so (\ref{repmap1})
holds. Therefore $\flats(L',E) \subseteq \flats(L,E)$.

Conversely, assume that $\flats(L',E) \subseteq \flats(L,E)$. We build a map
$\p: L \to L'$ as follows. Let $\vee'$ denote join in $L'$. For every
$\ell \in L$, we set $\ell\p = \vee'\{ e \in E \mid e \leq
\ell\}$.

It
is immediate that $\p$ is order-preserving. Hence, given $\ell_1,
\ell_2 \in L$, we have $\ell_i\p \leq (\ell_1 \vee \ell_2)\p$ for $i =
1,2$ and so 
\beq
\label{repmap2}
(\ell_1\p \vee \ell_2\p) \leq (\ell_1 \vee \ell_2)\p.
\eeq
Moreover, for every $e \in E$, we have $e\p = \vee'\{ f \in
E \mid f \leq e \mbox{ in } L\}$. Since $e \leq e$, we get $e \leq
e\p$ in $L'$.

Now take $e \in E$ such that $e \leq (\ell_1 \vee \ell_2)$. Since
$\flats(L',E) \subseteq \flats(L,E)$, we have $Z_{\ell_1\p \vee \ell_2\p} = Z_k$
for some $k \in L$. Let $f \in E$. Since $\p$ is order-preserving, $f
\leq \ell_i$ implies $f \leq f\p \leq \ell_i\p \leq (\ell_1\p \vee
\ell_2\p)$ and so $f \in Z_{\ell_1\p \vee \ell_2\p} = Z_k$. Hence
$\ell_i \leq k$ for $i = 1,2$ and so $(\ell_1 \vee \ell_2) \leq
k$. Thus $e \in Z_k = Z_{\ell_1\p \vee \ell_2\p}$ and so 
$$(\ell_1 \vee \ell_2)\p = \vee'\{ e \in E \mid e \leq
(\ell_1 \vee \ell_2) \} \leq (\ell_1\p \vee \ell_2\p).$$
Together with (\ref{repmap2}), this implies that $\p$ is a
$\vee$-map. 

It remains to be proved that $e\p \leq e$ holds in $L'$ for every $e
\in E$. Since $\flats(L',E) \subseteq \flats(L,E)$, we have
$$\{ f \in E \mid f \leq e \mbox{ in } L'\} =
\{ f \in E \mid f \leq m \mbox{ in } L\}$$
for some $m \in L$. Hence $e \leq m$ holds in $L$. It follows that,
for every $f \in E$, $f \leq e$ in $L$ implies $f \leq m$ in $L$ and
therefore $f \leq e$ in $L'$. Hence
$e\p = \vee'\{ f \in
E \mid f \leq e \mbox{ in } L\} = e$ and so $(L',E)
\leq (L,E)$.
\qed

Recall that, if $E' = \{ Z_e \mid e \in E\}$, then $(\flats(L,E),E')
\cong (L,E)$ holds for every $(L,E) \in \FLg$ by Proposition \ref{ltom}.
We identify $E'$ with $E$ to simplify notation.

Now $\flats(E,H)$ is closed under intersection by Proposition
\ref{trivia}(i).
We say that a $\cap$-semilattice $F$ of $(\flats(E,H),\subseteq)$ is
{\em full} if $\emptyset, E \in F$.
Let $\fisflats(E,H)$ denote the set of all full $\cap$-subsemilattices of
$(\flats(E,H),\subseteq)$. Then $(\fisflats(E,H),\subseteq)$ is a poset
closed under intersection, hence a $\wedge$-semilattice and therefore
a lattice with the determined join
$$(F_1 \vee F_2) = \cap \{ F \in \fisflats(E,H) \mid F_1 \cup F_2
\subseteq F\}.$$

We say that $(L,E) \in \FLg$ is a boolean representation of a simple
hereditary collection $(E,H)$ if $M(L,E)$ is a boolean representation
of $(E,H)$. Let $\BR(E,H)$ denote the class of all $(L,E) \in \FLg$
which are boolean representations of $(E,H)$. 
We restrict to $\BR(E,H)$ the quasi-order previously defined on
$\FLg$. If $(L,E) \in \BR(E,H)$, then by Proposition 
\ref{zeroes} we have $Z_{\ell} \in 
\flats(E,H)$ for every $\ell \in L$. By Theorem \ref{rephc} and
Proposition \ref{repmap}, if
$(E,H)$ is boolean representable, then $(\flats(E,H),E) \geq (L,E)$
for every $(L,E) \in \BR(E,H)$.

It is easy to check that
$$\begin{array}{rcl}
\theta:(\BR(E,H),\leq)&\to&(\fisflats(E,H),\leq)\\
(L,E)&\mapsto&\flats(L,E)
\end{array}$$
is a well-defined map.
Indeed, let $(L,E) \in \BR(E,H)$. Then $\flats(L,E) \subseteq
\flats(E,H)$ by Proposition 
\ref{zeroes}, ant it follows from (\ref{ltom2}) that
$\flats(L,E)$ is a $\cap$-subsemilattice of 
$\flats(E,H)$.
Note that $(L,E) \in \FLg$ implies $E \subseteq L
\setminus \{ B \}$.
Since $\emptyset = Z_{B}$ and $E = Z_T$, we have
$\flats(L,E) \in \fisflats(E,H)$ and so $\theta$ is a well-defined map.

Our next goal is to turn $\theta$ into an isomorphism. A first
obstacle is the fact that $\theta$ is not onto: not every $F \in
\fisflats(E,H)$ is rich enough to represent $(E,H)$. However, we claim that
$\fisflats(E,H) \setminus \im\theta$ is an order ideal of $\fisflats(E,H)$.

Since every $F \in
\fisflats(E,H)$, being a $\cap$-subsemilattice of $\flats(E,H)$,
constitutes a lattice of its own right with the determined join, then, in
view of Lemma \ref{indlh}, the question is whether the matrix arising
from $F$ produces enough witnesses to recognize all the independent
subsets in $H$. Therefore, if $F' \supseteq F$, every witness arising
from $F$ can also be obtained from $F'$ and so $\fisflats(E,H)
\setminus \im\theta$ is an ideal of $\fisflats(E,H)$. Let
$\fisflats_0(E,H)$ denote the Rees quotient (see Subsection \ref{101} of the
Appendix) of $\fisflats(E,H)$ by the
above ideal. By Proposition \ref{rees}, $\fisflats_0(E,H)
= \im\theta \cup \{ B \}$
is a lattice. 

On the other hand, adding a bottom element $B$ to $\BR(E,H)$, we
get a quasi-ordered set $\BR_0(E,H) = \BR(E,H) \cup \{ B \}$ and we
can extend $\theta$ to an onto map $\theta_0: \BR_0(E,H) \to
\fisflats_0(E,H)$ by setting $B\theta_0 = B$. Clearly, Proposition
\ref{repmap} immediately yields:

\bc
\label{psrepmap} 
For all $R,S \in \BR_0(E,H)$, $R \leq S$ if and only if $R\theta_0
\subseteq S\theta_0$.
\ec

Let $\rho$ be the
equivalence in $\BR_0(E,H)$ defined by $\rho = (\leq \cap \geq)$. 
Clearly, two representations $(L,E),(L',E)$ are $\rho$-equivalent if
there exists some lattice isomorphism $\p:L \to L'$ which is the
identity on $E$. Then
the quotient $\BR_0(E,H) / \rho$ becomes a poset and by Corollary
\ref{psrepmap}, the induced mapping
$\oo{\theta_0}: \BR_0(E,H) / \rho \to
\fisflats_0(E,H)$ is a poset isomorphism. Since we have already
remarked that $\fisflats_0(E,H)$ is a lattice (with the determined join),
we have proved the following theorem:

\bt
\label{iss}
Let $(E,H)$ be a simple boolean representable hereditary
collection. Then $\oo{\theta_0}:$ {\rm BR}$_0(E,H) / \rho \to$ {\rm
  FISFl}$_0(E,H)$ is a lattice isomorphism. 
\et

An {\em atom} of a lattice is an element covering the bottom element
$B$.
The atoms of $\BR_0(E,H)$ determine the {\em minimal} boolean
representations of $(E,H)$, and the sji elements of $\BR_0(E,H)$
determine the sji boolean
representations. Clearly, meet is given by intersection in
$\fisflats_0(E,H)$, collapsing into the bottom $B$ if it does not
correspond anymore to a representation of $(E,H)$. But how can join be
characterized in this lattice? 

\bp
\label{joinrep}
Let $(E,H)$ be a simple boolean representable hereditary
collection. Let $F,F' \in$ {\rm
  FISFl}$_0(E,H)$. Then:
\bi
\item[(i)] $(F \vee F') = F \cup F' \cup \{ Z \cap Z' \mid Z \in F,\; Z' \in
F'\}$.
\item[(ii)] If $(L,E)\theta = F$, $(L',E)\theta = F'$ and $(L'',E)\theta =
  (F \vee F')$, then $M(L'',E)$ is the closure of $M(L,E) \oplus_b
  M(L',E)$ under row sum in $\BB^{|E|}$.
\ei
\ep

\proof
(i) Clearly, the right hand side is the (full) $\cap$-subsemilattice
of $\flats(E,H)$ generated by $F \cup F'$.

(ii) Recall the isomorphism from Proposition \ref{ltom}. The rows
$r_Z$ of
$M(L,E)$ (respectively $M(L',E)$, $M(L'',E)$) are determined then by
the flats $Z$ in $F$ (respectively $F'$, $F \vee F'$). It is immediate
that $r_{Z \cap Z'} = r_Z + r_{Z'}$ in $\BB^{|E|}$, hence $M(L'',E)$
must be, up to permutation of rows, the stacking of $M(L,E)$ and
$M(L',E)$, to which we add (if needed) rows which are the sum in
$\BB^{|E|}$ of rows in $M(L',E)$ and $M(L'',E)$.
\qed

Next we introduce the notion of
{\em boolean sum} in $\BR(E,H)$. Given $(L,E), (L',E) \in
\BR(E,H)$, let $(L,E) \oplus_b (L',E)$ denote the $\vee$-subsemilattice
of the direct product $L \times L'$ $\vee$-generated by the diagonal $E_d =
\{ (e,e) \mid e \in E\} \subseteq L \times L'$. Taking the determined
meet, and identifying $E_d$ with $E$ as expectable, it
follows that $(L,E) \oplus_b (L',E) \in \FLg$. In fact, since the
projection $(L,E) \oplus_b (L',E) \to (E,L)$ is a $\vee$-map which is
the identity on $E$, it follows easily that $(L,E) \oplus_b (L',E) \in
\BR(E,H)$. But we can prove more:

\bp
\label{tens}
Let $(E,H)$ be a simple hereditary collection and let $(L,E), (L',E)
\in$ {\rm BR}$(E,H)$. Then 
\beq
\label{tens1}
(L,E)\rho \vee (L',E)\rho = ((L,E) \oplus_b
(L',E))\rho
\eeq 
holds in $\BR_0(E,H) / \rho$. Moreover,
 $M((L,E)\rho \vee (L',E)\rho)$ is the closure of the stacking matrix
 $M(L,E) \oplus_b
  M(L',E)$ under row sum in $\BB^{|E|}$.
\ep

\proof
By the preceding comment, we have $(L,E) \leq (L,E) \oplus_b (L',E)$
and also $(L',E) \leq (L,E) \oplus_b (L',E)$, hence
$$(L,E)\rho \vee (L',E)\rho \leq ((L,E) \oplus_b
(L',E))\rho.$$

Now let $(L'',E) \in
\BR(E,H)$ and suppose that $(L,E),(L',E) \leq (L'',E)$. We must show
that also $(L,E) \oplus_b
(L',E) \leq (L'',E)$. Indeed, there exist $\vee$-maps $\p:L \to L''$
and $\p':L' \to L''$ which fix $E$. Let $\p'': L \times L' \to L''$ be
defined by $(\ell,\ell')\p'' = (\ell\p \vee \ell'\p')$. Since
$((\ell_1,\ell'_1) \vee (\ell_2,\ell'_2)) = (\ell_1 \vee \ell_2,
\ell'_1 \vee \ell'_2)$ in $L \times L'$, it follows easily that $\p''$
is a $\vee$-map. Moreover, since both $\p$ and $\p'$ fix the elements
of $E$, so does $\p''$. Since the restriction of a $\vee$-map to a
$\vee$-subsemilattice is still a $\vee$-map, it follows that $(L,E) \oplus_b
(L',E) \leq (L'',E)$ and (\ref{tens1}) holds.

The last claim follows from Proposition \ref{joinrep}.
\qed

We can now state the following straightforward consequence:

\bc
\label{tours}
Let $(E,H)$ be a simple hereditary collection and consider $(L,E), (L',E)
\in$ {\rm BR}$(E,H)$. Then:
\bi
\item[(i)] $(L,E)$ can be decomposed as
a boolean sum (equivalently, stacking matrices and closing under row
sum) of sji representations;
\item[(ii)] this decomposition is not unique in
general, but becomes so if we take a maximal decomposition by
taking all the sji representations below $(L,E)$.
\ei
\ec

Examples shall be provided in Section \ref{seexe}.

\brem
\label{spcorps}
Let $(E,H)$ be a simple hereditary collection and let $(L,E), (L',E)
\in$ {\rm BR}$(E,H)$. It is reasonable to identify $(L,E)$ and
$(L',E)$ if some bijection of $E$ induces a $\vee$-bijection $L \to
L'$, and list only up to this identification in examples. However, for
purposes of boolean sum decompositions, the bijection on $E$ must be
the identity.
\erem

So we shall devote particular attention to  minimal/sji boolean
representation of $(E,H)$. How do 
these concepts relate to the flats in $\fisflats(E,H)$ and to the
matrices representing them? 

If $L$ is a lattice, we denote by $g_{\vee}(L)$ (respectively
$g_{\wedge}(L)$) the unique minimum set of $\vee$-generators
(respectively $\wedge$-generators) of $L$. Clearly, $g_{\vee}(L)$
equals the set of all sji elements of $L\setminus \{ B \}$. Similarly,
$g_{\wedge}(L)$ equals the set of all smi elements of $L\setminus \{ T \}$.

Given $(L,E) \in \BR(E,H)$, we may view $\flats(L,E)$ as a lattice with the
determined join and define
$\wh{Z}(L,E) = g_{\wedge}(\flats(L,E))$. That is,
$\wh{Z}(L,E)$ consists of all the smi flats in $\flats(L,E) \setminus \{
E \}$, i.e. which cannot be
nontrivially expressed as intersections of other flats in $\flats(L,E)$
(note that $E = \cap \emptyset$, hence $E \notin \wh{Z}(L,E)$).
In view of Proposition \ref{ltom} (and particularly (\ref{ltom2})), we have
\beq
\label{rob1}
\wh{Z}(L,E) = \{ Z_{\ell} \mid \ell \in g_{\wedge}(L) \}.
\eeq
If we transport these concepts into $M(L,E)$, then
$\wh{Z}(L,E)$ corresponds to the submatrix $\wh{M}(L,E)$ determined by
the rows which 
are not sums of other rows in $\BB^{|E|}$, excluding also the row with
just zeroes. By Proposition
\ref{stack}(ii), $\wh{M}(L,E)$ is still a boolean representation of
$(E,H)$. Note that, if we consider $\BB$ ordered by $0 < 1$ and the
direct product partial order in $\BB^{|E|}$, then the rows in
$\wh{M}(L,E)$ are precisely the sji rows of $M(L,E)$ for this partial
order. 

In the following three key propositions, we shall use the concept of
MPS and Proposition \ref{mps}, which the reader can find in 
Subsection \ref{101} of the Appendix. We shall use the notation
$\mathring{L} = L \setminus \{ T,B \}$. 

\bp
\label{flatco}
Let $(E,H)$ be a simple hereditary collection and let $(L,E), (L',E)
\in$ {\rm BR}$(E,H)$. Then the following conditions are equivalent:
\bi
\item[(i)] $(L,E)\rho$ covers $(L',E)\rho$ in {\rm BR}$_0(E,H)$;
\item[(ii)] there exists an MPS $\p:L \to L'$ fixing the elements of $E$;
\item[(iii)] {\rm Fl}$\,(L',E) =$ {\rm Fl}$\,(L,E) \setminus \{ Z_{\ell} \}$ for
  some smi $l 
  \in \mathring{L}$.
\ei
\ep

\proof
(i) $\Rw$ (ii). 
If $(L,E)\rho$ covers $(L',E)\rho$ in $\BR_0(E,H)$, then the
(onto) $\vee$-map $\p:L 
\to L'$ cannot be factorized as the composition of two proper
(onto) $\vee$-maps, and so $\p$ is an MPS. 

(ii) $\Rw$ (iii). By Proposition \ref{mps}, $\ker \p = \rho_{k,\ell}$ for
  some $k,\ell \in L$ such that $k$ covers $\ell$ and $\ell$ is
  smi, hence $\ell \neq T$. Therefore
  we may assume that $L' = L/\rho_{k,\ell}$. Clearly, $Z_{k \vee \ell}
  = Z_k$ and so (\ref{repmap1}) yields $\flats(L',E) = \flats(L,E) \setminus \{
  Z_{\ell} \}$. Note that we are assuming that $(L',E) \in \BR(E,H)
  \subseteq \FLg$, hence $E \subseteq L' \setminus \{ B \}$ and so
  $\ell \neq B$ (since $B$ is covered only by elements of
  $E$). Therefore (iii) holds. as desired.

(iii) $\Rw$ (i). It follows easily from (\ref{repmap1}) that $\ker\p$
has one class with two elements and all the others are singular, hence
$|L'| = |L|-1$ and so $(L,E)\rho$ covers $(L',E)\rho$ in $\BR_0(E,H)$.
\qed

This will help us to characterize minimal and sji boolean
representations of $(E,H)$ in terms of their flats. 

\bp
\label{mimi}
Let $(E,H)$ be a simple hereditary collection and let $(L,E)
\in$ {\rm BR}$(E,H)$. Then the following conditions are equivalent:
\bi
\item[(i)] $(L,E)$ is minimal;
\item[(ii)]
for every MPS $\p:L \to L'$ fixing the elements of $E$, $(L',E)
\notin$ {\rm BR}$(E,H)$;
\item[(iii)] for every smi $l \in \mathring{L}$, the matrix
  obtained by removing 
  the row $\ell$ from $M(L,E)$ is not a matrix boolean
representation of $(E,H)$;
\item[(iv)] for every smi $\ell \in \mathring{L}$, {\rm Fl}$\,(L,E)
  \setminus \{ Z_{\ell} 
  \} \notin$ {\rm Im}$\,\theta$.
\ei
\ep

\proof
(i) $\iff$ (iv). By Proposition \ref{flatco}.
 
(i) $\Rw$ (iii). Let $l \in \mathring{L}$ be an smi, and let $k$ be the unique
element of $L$ covering $\ell$. By Proposition \ref{flatco}, $L' =
L/\rho_{k,\ell}$ is a lattice and $M(L',E)$ is precisely the matrix
obtained by removing the row $\ell$ from $M(L,E)$. If $M(L',E)$ is a
boolean representation of $(E,H)$, then $(L,E)\rho$ covers
  $(L',E)\rho$ in $\BR_0(E,H)$ and so $(L,E)$ is not minimal.

(iii) $\Rw$ (iv). 
Suppose that $\flats(L,E) \setminus \{ Z_{\ell}
  \} \in \im\theta$ for some smi $\ell \in \mathring{L}$. Then $\flats(L,E)
  \setminus \{ Z_{\ell} 
  \} = \flats(L',E)$ for some $(L',E) \in \BR(E,H)$. It is straightforward
  to check that $M(L',E)$ is  the matrix
obtained by removing the row $\ell$ from $M(L,E)$. Thus (iii) fails.
\qed

\bp
\label{sjisji}
Let $(E,H)$ be a simple hereditary collection and let $(L,E)
\in$ {\rm BR}$(E,H)$. Then the following conditions are equivalent:
\bi
\item[(i)] $(L,E)$ is sji;
\item[(ii)]
up to isomorphism, there is at most one MPS $\p:L \to L'$ fixing the
elements of $E$ and such that $(L',E)
\in$ {\rm BR}$(E,H)$;
\item[(iii)] there exists at most one smi $l \in \mathring{L}$ such
  that the matrix 
  obtained by removing 
  the row $\ell$ from $M(L,E)$ is still a matrix boolean
representation of $(E,H)$;
\item[(iv)] there exists at most one smi $\ell \in \mathring{L}$ such
  that {\rm Fl}$\,(L,E)
  \setminus \{ Z_{\ell} 
  \} \in$ {\rm Im}$\,\theta$.
\ei
\ep

\proof
Clearly, $(L,E)$ is sji if and only if $(L,E)\rho$ covers exactly one
element in $\BR_0(E,H)$. Now we apply Proposition \ref{flatco},
proceeding analogously to the proof of Proposition \ref{mimi}.
\qed  

We call a reduced matrix boolean representation $M$ of $(E,H)$ {\em
  rowmin} if any matrix obtained by removing a row of $M$ is no longer
a boolean representation of $(E,H)$.

\bp
\label{doubt}
Let $(E,H)$ be a simple hereditary collection
and let $(L,E)\in$ {\rm BR}$(E,H)$ be minimal. Then $\wh{M}(L,E)$ is rowmin.
\ep

\proof
By Proposition \ref{mimi}, we cannot remove from
$\wh{M}(L,E)$ a row corresponding to an smi element of $\mathring{L}$.
Suppose now that $B$ is smi. Then $B$ is covered in $L$ by a unique
element $e$, necessarily in $E$ since $L$ is $\vee$-generated by $E$ and
so the unique 1 in the column of $M(L,E)$ determined by $e$ occurs at
the row of $B$. Since $\{ e \}$ is independent due to $(E,H)$ being
simple, it follows that the row of $B$ cannot be removed either.
\qed

The next example shows that the converse of Proposition \ref{doubt}
does not hold. 

\be
\label{libourne}
Let $(E,H)$ be the matroid defined by $E = \hat{4}$ and $H = P_3(E)
\setminus \{ 123 \}$. Then 
$$M = \left(
\begin{matrix}
1&&0&&1&&1\\
0&&1&&1&&0\\
0&&0&&0&&1
\end{matrix}
\right)$$
is a rowmin boolean representation of $(E,H)$ but $M \neq \wh{M}(L,E)$
for every minimal $(L,E)\in$ {\rm BR}$(E,H)$.
\ee

Indeed, it follows from the analysis developed later in Example
\ref{bigex} that $M = \wh{M}(L,E)$
for some sji $(L,E)\in \BR(E,H)$ (with $(L,E)\theta = \{ E, 123, 14,
1, 2, \emptyset \}$), being therefore a boolean representation of
$(E,H)$. It must be rowmin since it has only 3 rows and there exist
independent 3-subsets of $E$. However, the description of the minimal
cases in Example
\ref{bigex} shows that $M$ does not arise from any of them.

\section{Closure operators}

In this section we relate the closure operator
induced by a hereditary collection 
with the closure operator induced by a representation. %its lattice of flats:
In the follow-up, we denote by $\cl$  the closure operator on $2^E$
induced by $(E,H)$. Given $(L,E) \in \BR(E,H)$, we denote by $\clos_L$
the closure operator on $2^E$ as defined in Section \ref{cisl}.

\bl
\label{clopes}
Let $(E,H)$ be a simple hereditary collection and let $(L,E) \in$ {\rm
  BR}$(E,H)$. Let $X \subseteq E$. Then:
\bi
\item[(i)] {\rm Cl}$\,X \subseteq$ {\rm Cl}$_LX$;
\item[(ii)] {\rm Cl}$\,X =$ {\rm Cl}$_LX$ if $L =$ {\rm Fl}$\, (E,H)$.
\ei
\el

\proof
(i) 
We have
$\cl X = \cap\{ Z \in \flats(E,H) \mid X \subseteq
Z\}$ and in view of (\ref{clz}) and Proposition \ref{zeroes} also
$$\clos_LX = Z_{\vee X} \in \flats(E,H).$$
Since $X \subseteq \clos_LX$,
we get $\cl X \subseteq \clos_LX$.

(ii) Assume that $L = \flats(E,H)$. Let $Z \in \flats(E,H)$ be such
that $X \subseteq Z$. It suffices to show that $\clos_LX \subseteq
Z$, i.e. $Z_{\vee X} \subseteq Z$. Now in $\flats(E,H)$ we have $\vee
X = \cl X$ and $Z_{\cl X} = \{ e \in E \mid e \in \cl X \} = \cl X$,
hence we must show that 
$\cl X \subseteq Z$. This follows from $X \subseteq Z$ and $Z$ being
closed, therefore we are done.
\qed

\bp
\label{hcl}
Let $(E,H)$ be a simple hereditary collection admitting a boolean
representation $(L,E) \in$ {\rm BR}$(E,H)$ and let $X \subseteq E$. Then the
following conditions 
are equivalent:
\bi
\item[(i)] $X \in H$;
\item[(ii)] $X$ admits an enumeration $x_1,\ldots, x_k$ such that
\beq
\label{hcl2}
{\rm Cl}_L(x_1,\ldots, x_k) \supset {\rm Cl}_L(x_2,\ldots, x_k) \supset \ldots
\supset {\rm Cl}_L(x_k);
\eeq
\item[(iii)] $X$ admits an enumeration $x_1,\ldots, x_k$ such that
\beq
\label{hcl1}
{\rm Cl}(x_1,\ldots, x_k) \supset {\rm Cl}(x_2,\ldots, x_k) \supset \ldots
\supset {\rm Cl}(x_k).
\eeq
\ei
\ep

\proof
(i) $\Rw$ (ii). 
Assume $X \in H$. Since $(L,E) \in \BR(E,H)$, it follows that
$X$ is c-independent.
By Proposition \ref{cindch}, this is equivalent to saying that $X$ admits
an enumeration $x_1,\ldots , x_k$ such that (\ref{hcl2}) holds.

(ii) $\Rw$ (iii). Suppose that $x_i \in \cl (x_{i+1},\ldots,x_k)$ for
some $i$. By Lemma \ref{clopes}(i), we get $x_i \in
\clos_L(x_{i+1},\ldots,x_k)$ and so $\clos_L(x_{i},\ldots,x_k) =
\clos_L(x_{i+1},\ldots,x_k)$, a contradiction. Hence $x_i \notin \cl
(x_{i+1},\ldots,x_k)$ for every $i$ and so (\ref{hcl1}) holds.

(iii) $\Rw$ (i). Consider the representation of $(E,H)$ by
$(L',E) = (\flats(E,H),E)$.  By Lemma \ref{clopes}(ii), we have
$$\clos_{L'}(x_1,\ldots, x_k) \supset \clos_{L'}(x_2,\ldots, x_k) \supset \ldots
\supset \clos_{L'}(x_k).$$
It follows from Proposition
\ref{cindch} that $X$ is c-independent with respect to $M(\flats
(E,H),E)$, and so $X \in H$ by Lemma
\ref{indlh}.
\qed

We can now prove another characterization of boolean
representability:

\bt
\label{indecl}
Let $(E,H)$ be a simple hereditary collection. Then the following
conditions are equivalent:
\bi
\item[(i)] $(E,H)$ admits a boolean representation;
\item[(ii)] every $X \in H$ admits an enumeration $x_1,\ldots , x_k$
  satisfying (\ref{hcl1}).
\ei
\et

\proof
(i) $\Rw$ (ii). 
By Proposition \ref{hcl}.

(ii) $\Rw$ (i).
Let $X \in H$. In view of (\ref{hcl1}), we can use the flats
$\cl(x_i,\ldots,x_k)$ as a witness for $X$, hence $X$ is c-independent
with respect to $M(\flats(E,H),E)$. 
Lemma \ref{indlh} yields the reciprocal implication and so 
$M(\flats(E,H),E)$ is a reduced boolean
  representation of $(E,H)$.
\qed

We can also characterize which lattices provide boolean
representations:

\bp
\label{lbr}
Let $(E,H)$ be a boolean representable simple hereditary collection
and let $\p:(${\rm Fl}$(E,H),E) \to (L,E)$ be a $\vee$-map
fixing the elements of $E$. Then the following
conditions are equivalent:
\bi
\item[(i)] $(L,E) \in$ {\rm BR}$(E,H)$;
\item[(ii)] every $X \in H$ admits an enumeration $x_1,\ldots , x_k$
  satisfying (\ref{hcl2}).
\ei
\ep

\proof
(i) $\Rw$ (ii). 
By Proposition \ref{hcl}.

(ii) $\Rw$ (i). Let $X \subseteq E$. We show that $X \in H$ if and
only if $X$ is c-independent (with respect to $M(L,E)$).

Assume that $X \in H$. Since $\clos_LY = (\vee Y)\dw
\cap E = Z_{\vee Y}$ for every $Y \subseteq E$ by (\ref{clz}), it
follows from (ii) 
that the rows $\clos_L(x_i,\ldots,x_k)$ act as a witness for $X$ in 
$M(L,E)$ and so $X$ is c-independent.

Conversely, assume that $X$ is c-independent. By Proposition
\ref{repmap}, $M(L,E)$ is a submatrix of $M(\flats(E,H),E)$ and so 
$X$ is c-independent with respect to $M(\flats(E,H),E)$. Hence $X \in
H$ by Lemma \ref{indlh}.
\qed

\bc
\label{newlbr}
Let $(E,H)$ be a boolean representable simple hereditary collection
and let $F \in$ {\rm FISFl}$(E,H)$. For every $X \subseteq E$, let
{\rm Cl}$_FX = \cap \{ Z \in F \mid Z \supseteq X \}$.
Then the following
conditions are equivalent:
\bi
\item[(i)] $F =$ {\rm Fl}$\,(L,E)$ for some $(L,E) \in$ {\rm BR}$(E,H)$;
\item[(ii)] every $X \in H$ admits an enumeration $x_1,\ldots , x_k$
  satisfying 
\beq
\label{newlbr1}
{\rm Cl}_F(x_1,\ldots, x_k) \supset {\rm Cl}_F(x_2,\ldots, x_k) \supset \ldots
\supset {\rm Cl}_F(x_k);
\eeq
\ei
\ec

\proof
As noted before, since $F$ is a
$\cap$-subsemilattice of $\flats(E,H)$, it
constitutes a lattice of its own with the determined join 
$$(X \vee Y) = \clos_F(X \cup Y) \quad (X,Y \in F).$$
Identifying $E$ with $\{
\clos_F\{ e \} \mid e \in E\}$, we can view $(F,E)$ as an element of
$\FLg$, isomorphic to $(L,E)$ in view of Proposition \ref{ltom}. Now
we apply 
Proposition \ref{lbr}.
\qed

The important subcase of boolean representations of matroids was
studied in \cite{IR1,IR2} and the following fundamental result was
proved. For the sake of completeness, we include here a short
alternative proof: 

\bt
\label{boom}
{\rm \cite[Theorem 4.1]{IR2}}
Let $(E,H)$ be a simple matroid. Then $M(${\rm Fl}$\,(E,H),E)$ is a boolean
representation of $(E,H)$.
\et

\proof
In view of Theorems \ref{rephc} and \ref{indecl}, it suffices to show
that
every $X \in H$ admits an enumeration $x_1,\ldots , x_k$ satisfying
(\ref{hcl1}). Thus we only have to show that 
\beq
\label{boom1}
x_i \notin
\cl(x_1,\ldots,x_{i-1})
\eeq
for $i = 2,\ldots,k$. Indeed, suppose that $x_i \in
\cl(x_1,\ldots,x_{i-1})$. By Proposition \ref{circuits}, we have $x_i 
\in C \subseteq \{ x_1,\ldots,x_{i} \}$ for some circuit $C$, hence $C
\subseteq X \in H$, a contradiction.
Thus (\ref{boom1}) holds as required. 
\qed

\be
\label{fourpoints}
Let $(E,H)$ be a simple hereditary collection with $|E| = 4$. 
\bi
\item[(i)] If $H$ has 0, 3 or 4 independent 3-subsets, then $(E,H)$ is
  a matroid and therefore boolean representable.
\item[(ii)] If $H$ has 1 independent 3-subset, then $(E,H)$ does not
  satisfy (PR) and so is not boolean representable.
\item[(iii)] If $H$ has 2 independent 3-subsets, then $(E,H)$ is
  not a matroid but it is boolean representable.
\ei
\ee

Indeed, (i) and (ii) are straightforward (in view of \cite[Theorem
5.3]{IR1}). In (iii), we may assume that 
123 and 124 denote the independent 3-subsets. Since 34 is also a
basis, then $(E,H)$ is not a matroid. However, since $\flats(E,H) =
\P_1(E) \cup \{ 12, E\}$,
it follows easily from Theorem \ref{indecl} that $(E,H)$ is
  boolean representable. 

In fact, in this case the lattice of flats can be depicted as
$$
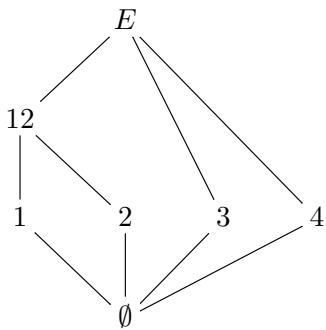
\xymatrix{
&E \ar@{-}[dl] \ar@{-}[ddr] \ar@{-}[ddrr] &&\\
12 \ar@{-}[d] \ar@{-}[dr] &&&\\
1 & 2 & 3 & 4\\ 
&\emptyset \ar@{-}[ul] \ar@{-}[ur] \ar@{-}[u] \ar@{-}[urr] &&
}$$
and so there exist maximal chains $\emptyset \subset 1 \subset 12
\subset E$ and $\emptyset \subset 4 \subset E$ of different
length. Hence $\flats(E,H)$ does not satisfy the {\em Jordan-Dedekind}
condition and so is not semimodular by \cite[Theorem 374]{Gra}. We
recall that a lattice $L$ is said to be 
{\em semimodular} if there is no sublattice of the form
$$
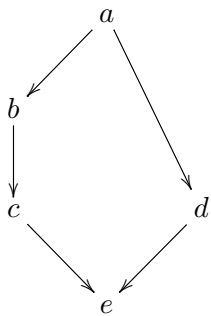
\xymatrix{
& a \ar[dl] \ar[ddr] & \\
b \ar[d] && \\
c \ar[dr] && d \ar[dl] \\
& e &
}$$
with $d$ covering $e$. 

A lattice is called {\em atomic} if every element
is a join of atoms ($B$ being the join of the empty set).
It is said to be {\em geometric} if it is both semimodular and
atomic. 
It is well known that a lattice is geometric if and only if it is
isometric to the lattice of flats of some matroid \cite[Theorem 1.7.5]{Oxl}.

Hence the above example shows that properties such as semimodularity
or the Jordan-Dedekind condition, which hold in the lattice of flats
of a matroid, may fail in the lattice of flats
of a boolean representable hereditary collection, even though it is
simple and paving (see Subsection \ref{spav}).

\section{Examples}
\label{seexe}

We present now some examples where we succeed on identifying all the
minimal and sji boolean representations.

\be
\label{bigex}
Let $(E,H)$ be the matroid defined by $E = \hat{4}$ and $H = P_3(E)
\setminus \{ 123 \}$. We compute all the minimal and sji
representations of $(E,H)$.
\ee

It is routine to compute $\flats(E,H) = P_1(E) \cup \{ 14, 24, 34,
123, 1234 \}$. Which $F \in \fisflats(E,H)$ correspond to lattice
representations (i.e. $F \in \im\theta$)? We claim that $F \in
\im\theta$ if and only if one of 
the following conditions is satisfied:
\beq
\label{bigex1}
123 \in F \mbox{ and } |\{ 1, 2, 3 \} \cap F| \geq 2,
\eeq
\beq
\label{bigex2}
|\{ 14, 24, 34 \} \cap F| \geq 2.
\eeq
In view of Corollary \ref{newlbr}, it is easy to see that any of
the conditions implies $F \in \im\theta$ (note that $4 = 14 \cap 24
\in F$ in the case (\ref{bigex2})).

Conversely, assume that $F \in \im\theta$ and suppose that $|\{ 14,
24, 34 \} \cap F| \leq 1$. Without loss of generality, we may assume
that $24, 34 \notin F$. Since $234 \in H$, it follows from Corollary
\ref{newlbr} that there exists an enumeration $x,y,z$ of 234 such that
$$\clos_F(x) \subset \clos_F(xy) \subset \clos_F(xyz).$$
The only possibility for $\clos_F(xy)$ in $F$ is now $123$. Hence
$\clos_F(x) \in \{ 2,3\}$. Out of symmetry, we may assume that $2 \in
F$. On the other hand, since $13 \in H$, there exists an enumeration
$a,b$ of 13 such that 
2$$\clos_F(a) \subset \clos_F(ab) = 123.$$
The only possibilities for $\clos_F(a)$ in $F$ are now $1$ or 3, hence
(\ref{bigex1}) holds.

We consider now the minimal case.
Assume that $F \in \im\theta$. Since we can view $(F,E)$ as a lattice
projecting onto $F$ through $\theta$, and by Proposition \ref{mimi},
the key lies with the smi elements of $F$ (with respect to
intersection). More precisely, $(F,E) \in \BR(E,H)$ is
minimal if and only if removal of an smi element of $F \setminus \{
E, \emptyset\}$ takes us outside $\im\theta$. It follows easily from our
characterization of $\im\theta$ that this corresponds to having
$$F = \{ E, 123, i, j, \emptyset\} \quad \mbox{or} \quad F = \{
E, i4, j4, 4, \emptyset \}$$
for some distinct $i,j \in \hat{3}$, leading to the lattices
$$
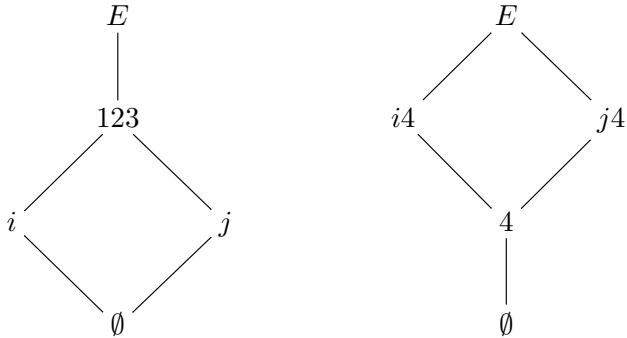
\xymatrix{
&E \ar@{-}[d] &&&&E \ar@{-}[dl] \ar@{-}[dr] & \\
&123 \ar@{-}[dl] \ar@{-}[dr] &&&i4 \ar@{-}[dr] &&j4 \ar@{-}[dl] \\
i \ar@{-}[dr] && j \ar@{-}[dl] &&&4 \ar@{-}[d] &\\
&\emptyset &&&&\emptyset &
}$$
Note that, if we wish to identify the generators $E$ in these
lattices, we only have to look for $\cl e$ for each $e \in E$. For
instance, in the first lattice, the top element corresponds to the
generator 4. 

Following Remark \ref{spcorps}, we can count the number of minimal
lattice representations
\bi
\item
up to identity in the $\vee$-generating set $E$ $(3+3 = 6)$;
\item
up to some bijection of $E$ inducing a $\vee$-bijection on the
lattices $(1+1 = 2)$.
\ei

With respect to the sji representations, it follows from 
Proposition \ref{sjisji} that $(F,E) \in \BR(E,H)$ is
sji if and only if there is at most one smi element of $F \setminus \{
E, \emptyset\}$ whose removal keeps us inside $\im\theta$. We claim
that this corresponds to either the minimal case or having  
\beq
\label{bigex3}
F = \{ E, 123, i, j, 4, \emptyset\} \quad \mbox{or} \quad
F = \{ E, 123, i4, i, j, \emptyset\} \quad \mbox{or} \quad F = \{
E, i4, j4, k, 4, \emptyset \}
\eeq
for some $i,j,k \in \hat{3}$ with $i \neq j$. It is immediate that the
cases in (\ref{bigex3}) lead to $(F,E)$ sji, the only possible
removals being respectively $4$, $i4$ and $k$. Conversely, assume that $F$
corresponds to an sji non minimal case. Suppose first that $F$
satisfies (\ref{bigex1}). If none of the pairs $k4$ is in $F$, then
$F$ must contain precisely three singletons to avoid the minimal case,
and one of them must be 4 to avoid having a mutiple choice in the
occasion of removing one of them. Hence we may
assume that $i4 \in F$ and so also $i = 123 \cap i4$. If $j4 \in F$
for another $j \in \hat{3}$, then also $j \in F$ and we would have the
option of removing either $i4$ or $j4$. Hence $F$ contains
$E, 123, i4, i, \emptyset,$ and possibly any other singletons. In
fact, it must contain at least one in view of (\ref{bigex1}) but
obviously not both. Thus $F$ is of the first form in (\ref{bigex3}) in
this case.

Now assume that $F$
satisfies (\ref{bigex2}) but not (\ref{bigex1}). Assume that $i4,j4
\in F$ for some distinct $i,j,k \in \hat{3}$.
Then $123 \notin F$, otherwise $i,j \in F$ and we can remove either
$i4$ or $j4$. Clearly, a third pair $k4$ is forbidden, otherwise we
could remove any one of the three pairs. Thus $F$ contains
$E, i4, j4, 4, \emptyset,$ and possibly any other singletons. In
fact, it must contain at least one to avoid the minimal case but
obviously not two, since any of them could be removed. Thus $F$ is of
the second form in (\ref{bigex3}) and we have identified all the sji cases.
In the second case of (\ref{bigex3}), we must separate the subcases $k
= i$ and $k \notin \{ i,j \}$. Thus the sji non minimal cases lead to
the lattices 
$$
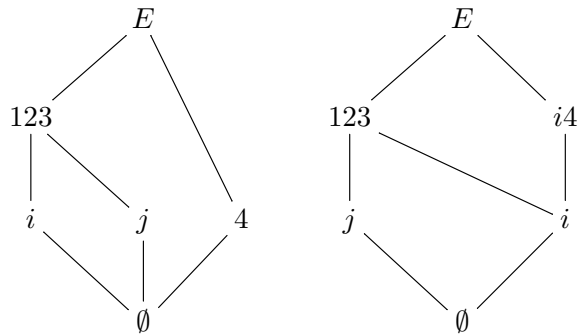
\xymatrix{
&E \ar@{-}[dl] \ar@{-}[ddr] &&&E \ar@{-}[dl] \ar@{-}[dr] &&\\
123 \ar@{-}[d] \ar@{-}[dr] &&& 123 \ar@{-}[d] \ar@{-}[drr] && i4 \ar@{-}[d] & \\
i & j & 4 & j && i &\\ 
&\emptyset \ar@{-}[ul] \ar@{-}[ur] \ar@{-}[u]&& &\emptyset \ar@{-}[ul]
\ar@{-}[ur] && 
}
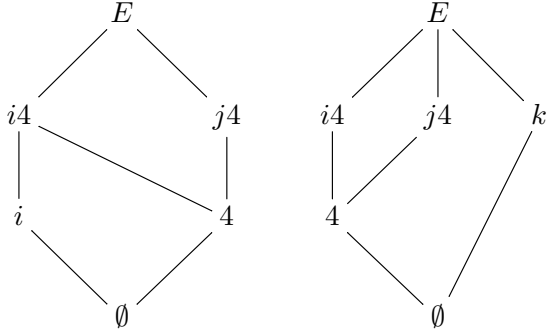
$$
$$\xymatrix{
&E \ar@{-}[dl] \ar@{-}[dr] &&&E
\ar@{-}[dl] \ar@{-}[dr] \ar@{-}[d] &\\
i4 \ar@{-}[d] \ar@{-}[drr]
&&j4 \ar@{-}[d] & i4
\ar@{-}[d] & j4 \ar@{-}[dl] & k \ar@{-}[ddl] \\
i && 4 & 4 &&\\ 
& \emptyset \ar@{-}[ul]
\ar@{-}[ur] &&& \emptyset \ar@{-}[ul] &
}$$

Following Remark \ref{spcorps} as in the minimal case, the number of sji
lattice representations in both counts (which includes the minimal
ones) is respectively $6+3 + 6 + 6 + 3
= 24$ and $1+1+1+1+1 = 5$.

It is easy to see that 
$$\flats(E,H) = \{ E, 123, 34, 2, 3, \emptyset \} \cup \{ E, 14, 24,
1, 4, \emptyset \}$$ 
provides a decomposition of the top boolean representation
$\flats(E,H)$ as the join of two sji's. In matrix form, and in view of
Proposition \ref{tens}, this corresponds to express the matrix
$$M(\flats(E,H),E) = \left(
\begin{matrix}
0&&0&&0&&0\\
0&&0&&0&&1\\
0&&1&&1&&0\\
1&&0&&1&&0\\
1&&1&&0&&0\\
0&&1&&1&&1\\
1&&0&&1&&1\\
1&&1&&0&&1\\
1&&1&&1&&0\\
1&&1&&1&&1
\end{matrix}
\right)$$
as the stacking of the matrices
$$\left(
\begin{matrix}
0&&0&&0&&0\\
0&&0&&0&&1\\
1&&1&&0&&0\\
1&&0&&1&&1\\
1&&1&&0&&1\\
1&&1&&1&&1
\end{matrix}
\right)
\hspace{2cm}
\left(
\begin{matrix}
0&&0&&0&&0\\
0&&1&&1&&0\\
1&&0&&1&&0\\
0&&1&&1&&1\\
1&&1&&1&&0\\
1&&1&&1&&1
\end{matrix}
\right)$$
Note that the maximal decomposition of $\flats(E,H)$ as
join of sji's would include 24 factors.

It is also easy to see, in view of 
% the characterization of the minimal case and 
Proposition \ref{stack}(ii) (which allows us to discard those
rows corresponding to non smi elements) that $\mindeg(E,H) = 3$: we
take the minimal representation defined by $F = \{
E, 14, 24, 4, \emptyset \}$ and discard the row corrresponding to $4 = 14 \cap 
24$. We can also discard the useless row of zeroes corresponding to
$\emptyset$ to get the matrix
$$\left(
\begin{matrix}
0&&1&&1&&0\\
1&&0&&1&&0\\
1&&1&&1&&1
\end{matrix}
\right)$$

We cannot do better than this since there are independent 3-sets in
$H$. Therefore $\mindeg(E,H) = 3$.

\be
\label{fano}
Let $(E,H)$ be the Fano matroid defined by $E = \hat{7}$ and $H = P_3(E)
\setminus \{ 125, 137$, $146, 236, 247, 345, 567 \}$. We compute all the
minimal and sji 
representations of $(E,H)$.
\ee

Note that $\L = P_3(E) \setminus H$ is precisely the set of {\em lines} in
the Fano plane \cite{W7} (the projective plane of order 2 over the two
element field):
$$
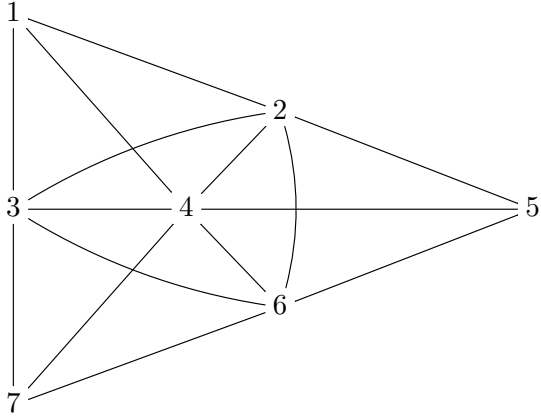
\xymatrix{
1 \ar@{-}[ddrr] \ar@{-}[dd] \ar@{-}[drrr] &&&&&& \\
&&& 2 \ar@{-}[dl] \ar@{-}@/_/[dlll] &&& \\
3 \ar@{-}[dd] \ar@{-}[rr] \ar@{-}@/_/[drrr] && 4 \ar@{-}[rrrr]
\ar@{-}[dr] &&&& 5 
\ar@{-}[dlll] \ar@{-}[ulll] \\
&&& 6 \ar@{-}[dlll] \ar@{-}@/_/[uu] &&& \\
7 \ar@{-}[uurr] &&&&&&
}$$

We can view the Fano plane as a PEG (see Subsection \ref{sgeom} of the
Appendix). We list a few of its properties: 
\bi
\item[(F1)]
Any two distinct lines intersect at a single point.
\item[(F2)]
Every point belongs to exactly three lines.
\item[(F3)]
Any two points belong to some line.
\item[(F4)]
If $K$ consists of 5 lines, then there exist two points $a,b$ such
that $K$ consists of all lines containing either $a$ or $b$.
% \item[(F4)]
% Given 4 points, we can always exclude one and not get a line.
\ei
We note that (F4) follows from (F1) since the two lines not in $K$
must cover exactly 5 points, and we may take $a,b$ as the two
remaining points.

It is easy to check that $\flats(E,H) = P_1(E) \cup \L$: the lines are
obviously closed, the 2-subsets are not, and every 4-subset of $E$
contains an independent set and has therefore closure $E$ by
Proposition \ref{newtrivia}.

Let $F \in
\fisflats(E,H)$. We claim that $F \in 
\im\theta$ if and only if $|F \cap \L| \geq 5$ or
\beq
\label{fano1}
|F \cap \L| = 4 \mbox{ and no 3 lines of $F \cap \L$ intersect at a point.}
\eeq

Assume that $F \in \im\theta$. Suppose
that $|F \cap \L| < 5$ and (\ref{fano1}) fails. If $|F \cap \L| \leq
3$, we can extend $F$ 
to some set $F'$ satisfying 
$|F' \cap \L| = 4$ and having 3 lines intersecting at a
point. Otherwise, let $F' = F$.

Suppose that $F' \cap \L = \{ X_1,X_2,X_3,X_4\}$ and
$X_1,X_2,X_3$ intersect at a certain point $p$. For $i = 1,2,3$, take some
point $x_i \in X_i\setminus (X_4 \cup \{ p\})$. The points are all
distinct by (F1). Let $Y = x_1x_2x_3$. Since $X_4 \cap Y = \emptyset$,
it follows from (F1) that $Y \in H$.  
It follows from Corollary
\ref{newlbr} that there exists an enumeration $y_1,y_2,y_3$ of Y such that
$$\clos_{F'}(y_1) \subset \clos_{F'}(y_1y_2) \subset
\clos_{F'}(y_1y_2y_3).$$
Hence $\clos_{F'}(y_1y_2) \in F' \cap \L$. Since $Y \cap X_4 =
\emptyset$, it follows that $|Y \cap X_i| \geq 2$ for some $i \in
\hat{3}$. Hence $x_j \in X_i$ for some $j \neq i$, yielding $|X_i \cap
X_j| \geq 2$ and contradicting (F1).

Assume now that (\ref{fano1}) holds. 
Let $X = xyz \in H$. By (F3), we
may write $x'yz, xy'z, xyz' \in \L$ for some $x',y',z' \in E$. 
Since $x' \neq x$ due to $xyz \in H$, $x'yz \in \L$ (and
similarly, $y' \neq y$), it follows that the three lines $x'yz, xy'z,
xyz'$ are distinct.

Suppose
that $x'yz, xy'z, xyz' \notin F$. Since $|x,y,z,x',y',z'| \leq 6$,
there exists some point $p$ which occurs in no line among $x'yz, xy'z,
xyz'$. Since $|\L| = 7$, it follows that the three lines
of $\L$ containing $p$ must be all in $F \cap L = \L \setminus \{
x'yz, xy'z, xyz' \}$, contradicting (\ref{fano1}). Thus we may assume
without loss of generality that $x'yz \in F$. 

On the other hand, in view of (F1) and (\ref{fano1}), the number of
intersections of lines in $F$ is precisely $\binom{4}{2} = 6$, 
which implies that $F$ contains at least 6 points among the 7 of
$E$. Hence $y \in F$ or $z \in F$ and we may assume without loss of
generality that $y \in F$. Thus
$$\clos_F(y) = y \subset x'yz = \clos_F(yz) \subset E = \clos_F(xyz)$$
and so $F \in 
\im\theta$ by Corollary \ref{newlbr}.

Finally, if $|F \cap \L| \geq 5$, it suffices to show that $F$ contains some $F'
\in \fisflats(E,H)$ satisfying (\ref{fano1}). Let $K$ be any 5-subset
of lines in $F$, and let $a,b \in E$ be given by (F4). Let $F'$ be obtained
from $F$ by 
excluding the lines not in $K$ and the line containing $a,b$. It is
easy to check that $F'$ satisfies (\ref{fano1}).
This completes the proof for the characterization of $\im\theta$.

Now, similarly to Example \ref{bigex}, the lattice $(F,E)$ is a
minimal boolean representation if and only if removal of an smi
element of $F \setminus \{ E, \emptyset\}$ takes us outside
$\im\theta$. It follows easily from our 
characterization of $\im\theta$ (and the fact that $|F \cap \L| \geq
5$ implies that $F$ contains some $F'
\in \fisflats(E,H)$ satisfying (\ref{fano1})) that this corresponds to $F$
satisfying (\ref{fano1}) and having only the 6 points that are
necessarily present as the outcome of intersections of lines in $F$.
Writing $F \cap \L = \{ p,q,r,s\}$ and denoting by $xy$ the
intersection of $x,y \in F \cap \L$, we see that the minimal boolean
representations are, up to isomorphism, given by the lattice
$$
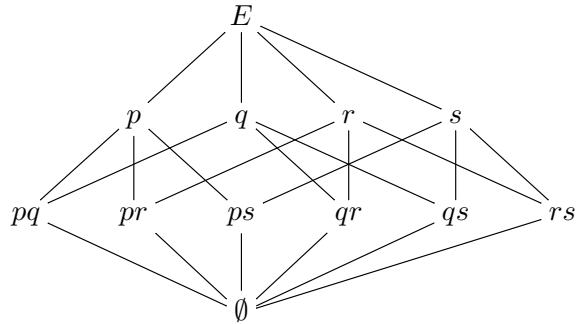
\xymatrix{
&&E \ar@{-}[d] \ar@{-}[dl] \ar@{-}[dr] \ar@{-}[drr] &&&\\
&p &q &r &s &\\
pq \ar@{-}[ur] \ar@{-}[urr] & pr \ar@{-}[u] \ar@{-}[urr] & ps \ar@{-}[ul]
\ar@{-}[urr] & qr \ar@{-}[ul] \ar@{-}[u] & qs \ar@{-}[ull] \ar@{-}[u] & rs
\ar@{-}[ull] \ar@{-}[ul] \\
&&\emptyset \ar@{-}[ull] \ar@{-}[ul] \ar@{-}[u] \ar@{-}[ur] \ar@{-}[urr]
\ar@{-}[urrr] &&&
}$$
Up to permuting rows/columns, the matrix representation $\wh{M}(F,E)$
(where we keep only the rows of $M(F,E)$ corresponding to the smi
elements of $F \setminus \{ E \}$) is then of the form
\beq
\label{fano9}
\left(
\begin{matrix}
0&0&1&1&0&1&1\\
0&1&1&0&1&0&1\\
1&0&0&1&1&0&1\\
1&1&0&0&0&1&1
\end{matrix}
\right)
\eeq

Following Remark \ref{spcorps} as in Example \ref{bigex}, the number of minimal
lattice representations in both counts is respectively $7$ and 1.

It follows from 
Proposition \ref{sjisji} that $(F,E) \in \BR(E,H)$ is
sji if and only if there is at most one smi element of $F \setminus \{
E, \emptyset\}$ whose removal keeps us inside $\im\theta$. We claim
that this corresponds to either the minimal case or one of the following:  
\bi
\item[(A)] $F$ satisfies (\ref{fano1}) and contains all the points;
\item[(B)] $|F \cap L| = 5$ and $F$ contains only 6 points.
\ei
Assume that (A) holds. We have already remarked that the intersections
of the 4 lines in $F$ yield 6 distinct points, hence the unique
possible removal is the 7th point.

Assume now that (B) holds. Clearly, we cannot remove a point. Let
$a,b \in E$ be given by (F4). Then the
unique line we can remove is the line containing $a$ and $b$, in
order to satisfy (\ref{fano1}). Thus both (A) and (B) correspond to
sji (non minimal) cases.

Assume now that $F$ corresponds to an sji non minimal case. It is
imediate that $F$ must contain 4 or 5 lines. Suppose that $|F \cap L|
= 4$. Then $F$ satisfies (\ref{fano1}) and must contain all the points
to avoid the minimal case. Hence (A) holds. Finally, we suppose that $|F \cap L|
= 5$. 
Let $a,b \in E$ be given by (F4).
Since we have $\binom{5}{2} = 10$ pairs of lines and 6 of
these pairs intersect in either $a$ or $b$, $F$ must contain only 6
points. Therefore (B) holds.

Which lattice corresponds to (A)? Let us use the same notation as in
the minimal case and denote the seventh point by $z$. Since every line
in $F$ must intersect each of the other three, and always at different
points, it follows that $z$ does not belong to any line in
$F$. Therefore we obtain the lattice
$$
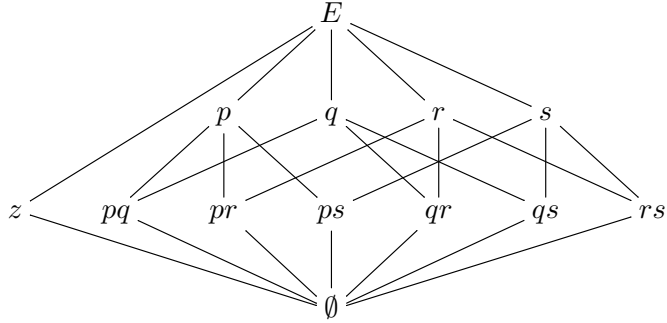
\xymatrix{
&&&E \ar@{-}[d] \ar@{-}[dl] \ar@{-}[ddlll] \ar@{-}[dr] \ar@{-}[drr] &&&\\
&&p &q &r &s &\\
z & pq \ar@{-}[ur] \ar@{-}[urr] & pr \ar@{-}[u] \ar@{-}[urr] & ps \ar@{-}[ul]
\ar@{-}[urr] & qr \ar@{-}[ul] \ar@{-}[u] & qs \ar@{-}[ull] \ar@{-}[u] & rs
\ar@{-}[ull] \ar@{-}[ul] \\
&&&\emptyset \ar@{-}[ulll] \ar@{-}[ull] \ar@{-}[ul] \ar@{-}[u]
\ar@{-}[ur] \ar@{-}[urr] 
\ar@{-}[urrr] &&&
}$$
Up to permuting rows/columns, the matrix representation $\wh{M}(F,E)$
is then of the form
$$\left(
\begin{matrix}
0&0&1&1&0&1&1\\
0&1&1&0&1&0&1\\
1&0&0&1&1&0&1\\
1&1&0&0&0&1&1\\
1&1&1&1&1&1&0
\end{matrix}
\right)$$

In case (B), let $r$ denote the line containing $a$ and $b$, and let
$p,q$ (respectively $s,t$) be the two other lines containing $a$
(respectively $b$). Using the same notation for intersection of
lines, we get the lattice
$$
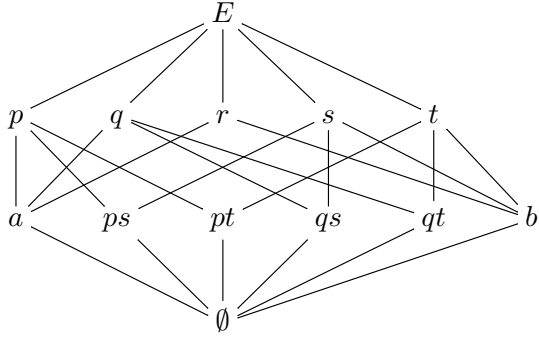
\xymatrix{
&&E \ar@{-}[d] \ar@{-}[dl] \ar@{-}[dll] \ar@{-}[dr] \ar@{-}[drr] &&&\\
p &q &r &s &t &\\
a \ar@{-}[u] \ar@{-}[ur] \ar@{-}[urr] & ps \ar@{-}[ul] \ar@{-}[urr] &
pt \ar@{-}[ull] 
\ar@{-}[urr] & qs \ar@{-}[ull] \ar@{-}[u] & qt \ar@{-}[ulll] \ar@{-}[u] & b
\ar@{-}[ulll] \ar@{-}[ull] \ar@{-}[ul] \\
&&\emptyset \ar@{-}[ull] \ar@{-}[ul] \ar@{-}[u] \ar@{-}[ur] \ar@{-}[urr]
\ar@{-}[urrr] &&&
}$$
Up to permuting rows/columns, the matrix representation $\wh{M}(F,E)$
is then of the form
$$\left(
\begin{matrix}
0&0&1&1&0&1&1\\
0&1&1&0&1&0&1\\
1&0&0&1&1&0&1\\
1&1&0&0&0&1&1\\
1&1&1&1&0&0&0
\end{matrix}
\right)$$

Following Remark \ref{spcorps} as in the minimal case, the number of sji
lattice representations in both counts is respectively $7+7+21 = 35$ and
$1+1+1 = 3$.

It is easy to see that 
$$\begin{array}{lll}
\flats(E,H)&=&\{ E, 125, 146, 236, 345, 567, 1, 2, 3, 4, 5, 6,
\emptyset \} \\
&\cup &\{ E, 137, 146, 236, 247, 567, 1, 2, 3, 4, 6, 7
\emptyset \}
\end{array}$$ 
provides a decomposition of the top boolean representation
$\flats(E,H)$ as the join of two sji's. In matrix form, and in view of
Proposition \ref{tens}, this corresponds to the stacking of matrices
$$\left(
\begin{matrix}
0&0&1&1&0&1&1\\
0&1&0&1&1&1&0\\
0&1&1&0&1&0&1\\
1&0&0&1&1&0&1\\
1&0&1&0&1&1&0\\
1&1&0&0&0&1&1\\
1&1&1&1&0&0&0
\end{matrix}
\right)
= 
\left(
\begin{matrix}
0&0&1&1&0&1&1\\
0&1&1&0&1&0&1\\
1&0&0&1&1&0&1\\
1&1&0&0&0&1&1\\
1&1&1&1&0&0&0
\end{matrix}
\right)
\oplus_b
\left(
\begin{matrix}
0&1&0&1&1&1&0\\
0&1&1&0&1&0&1\\
1&0&0&1&1&0&1\\
1&0&1&0&1&1&0\\
1&1&1&1&0&0&0
\end{matrix}
\right)
$$
where we depict only the rows corresponding to the smi elements of the
lattices (minus the top).

Note that the maximal decomposition of $\flats(E,H)$ as
join of sji's would include 35 factors.

We claim that $\mindeg(E,H) = 4$ taking the
matrix representation (\ref{fano9}) given for the minimal case.
Indeed, suppose that $M$ is a 3-row matrix representation of
$(E,H)$. Since all the 2-subsets are c-independent, all columns must
be different and nonzero in view of Proposition \ref{nons}. Since each
of our 7 columns has 3 entries, there exist 4 columns 
(corresponding to some distinct $a,b,c,d \in E$) 
having at most one zero. By Proposition \ref{nons}, all the 3-subsets
of $\{ a,b,c,d\}$ must be dependent, thus lines. In particular, two
lines may have two points in common, contradicting (F1). Therefore 
$\mindeg(E,H) = 4$.

Further information on the Fano plane can be found in \cite{RSil}. 

\smallskip

Before presenting the next example, we need to recall two standard concepts
from graph theory.

Given a
finite graph $G$, the {\em girth} of $G$, denoted by $\gth G$, is the
length of the shortest cycle in $G$ (assumed to be $\infty$ is $G$ is
acyclic). Note that $\gth G \geq 3$ for any finite graph. We denote
the maximum degree of a vertex in $G$ by $\maxdeg G$.

If $G = (V,E)$ is a connected graph, we can define a metric
$d$ on $V$ by
$$d(v,w) = \mbox{ length of the shortest path connecting $v$ and $w$
  (counting edges)}.$$
The {\em diameter} of $G$, denoted by $\diam G$, is the maximum value
in the image of $d$. If $G$ is not connected, we define $\diam G =
\infty$. 

Finally, let $K_{m,n}$ denote the complete bipartite graph on $m+n$
vertices.
% and let $K'_{m,n}$ denote the graph obtained by removing an edge
% from $K_{m,n}$. 

\be
\label{uthreeb}
Let $b \geq 5$ and $(E,H) = U_{3,b}$. We compute all the minimal and sji
representations of $(E,H)$.
\ee

It is immediate that $\flats(E,H) = P_2(E) \cup \{ E \}$. Given $F \in
\fisflats(E,H)$, we define a finite undirected graph $F\gamma$ with
vertex set $E = \hat{b}$ and edges $p \edge q$ whenever $p,q$ are
distinct and $pq \notin F$. We claim that 
\beq
\label{uthreeb1}
F \in \im\theta \; \iff \; (\gth F\gamma > 3 \mbox{ and } |E \setminus
F| \leq 1). 
\eeq

Indeed, assume that $F \in \im\theta$. Suppose that $\gth F\gamma =
3$. Then there exist distinct $p,q,r \in E$ 
such that $pq, pr, qr \notin F$. Hence $\clos_F(xy) = E$ for all distinct
$x,y \in \{ p,q,r\}$. Since $pqr \in H$, this contradicts $F \notin
\im\theta$ in view of Corollary
\ref{newlbr}. Thus $\gth F\gamma > 3$.

Suppose next that $x,y \in E \setminus
F$ are distinct. Let $z,t,w \in E \setminus \{ x,y \}$ be
distinct. By Corollary
\ref{newlbr}, $xyz$ admits an enumeration $x_1, x_2, x_3$ satisfying
\beq
\label{uthreeb2}
\clos_F(x_1) \subset \clos_F(x_1x_2) \subset
\clos_F(x_1x_2x_3).
\eeq
Since $F \subseteq P_2(E) \cup \{ E \}$, we get $x_1, x_1x_2 \in F$
and so $x_1 = z$ and $i_zz \in F$ for some $i_z \in \{ x, y
\}$. Similarly, $i_tt, i_ww \in F$ for some $i_t, i_w \in \{ x, y
\}$. We may thus assume that $i_z = i_t = x$, hence $x = i_zz \cap
i_tt \in F$, a contradiction. Therefore $|E \setminus
F| \leq 1$.

Conversely, assume that $\gth F\gamma > 3$ and $|E \setminus
F| \leq 1$. Let $x,y,z \in E$ be distinct. By Corollary
\ref{newlbr}, it suffices to show that
$xyz$ admits an enumeration $x_1,x_2,x_3$ satisfying
(\ref{uthreeb2}). Since $\gth F\gamma > 3$, we have $\{ xy, xz, yz \}
\cap F \neq \emptyset$. We may assume that $xy \in F$. Since $|E \setminus
F| \leq 1$, we have either $x \in F$ or $y \in F$. In any case,
(\ref{uthreeb2}) is satisfied by some enumeration of $x,y,z$ and so
(\ref{uthreeb1}) holds.

The minimal cases are once more characterized by the following
property: removal of an smi element of $F \setminus \{
E, \emptyset\}$ must make (\ref{uthreeb1}) fail. It is easy to see that
the smi elements of $F \setminus \{
E, \emptyset\}$ are precisely the 2-sets and the points which are not
intersections of 2-sets in $F$, i.e. vertices of degree $\geq |E|-2$
in $F\gamma$.
We claim that $(F,E)$ is
minimal if and only
\bi
\item[(M)] $\gth F\gamma > 3$, $\diam F\gamma = 2$ and 
\beq
\label{uthreeb3}
\maxdeg F\gamma \geq |E|-2 \Rw |E \setminus F| = 1.
\eeq
\ei

Assume that $(F,E)$ is
minimal. Then $\gth F\gamma > 3$ and $|E \setminus F| = 1$ by
(\ref{uthreeb1}). Suppose that $\maxdeg F\gamma \geq |E|-2$. Then some
points of $E$ are smi elements of $F\gamma$. If $E \subseteq F$, we
could remove one of these smi points and still satisfy
(\ref{uthreeb1}). Hence $|E \setminus F| = 1$. 

Finally, since $|E| > 2$ and $\gth F\gamma > 3$, we have $\diam
F\gamma \geq 2$. Suppose that $x,y \in E$ lie at distance $> 2$ in
$F\gamma$. Then we could add an edge $x \edge y$ and still satisfy
(\ref{uthreeb1}). Since adding an edge corresponds to removal of the
smi $xy$ from $F$, this contradicts $(L,E)$ being minimal.

Conversely, assume that $\gth F\gamma > 3$, $\diam F\gamma = 2$ and
(\ref{uthreeb3}) holds. Since $\diam F\gamma = 2$, it is clear that we
cannot add any extra edge and keep $\gth F\gamma > 3$, hence removal
of 2-sets from $F$ is forbidden. On the other hand, in view of
(\ref{uthreeb1}), we can only remove a point from $F$ if $E \subseteq
F$, and by (\ref{uthreeb3}) this can only happen if $\maxdeg F\gamma <
|E|-2$. However, as remarked before, this implies that no point is an
smi element of $F$. Therefore $(L,E)$ is
minimal as claimed.

Next we show that $(F,E)$ is sji if and only if $\gth F\gamma > 3$ and
one of the following cases holds:
\bi
\item[(A)] 
There exists a unique 2-subset $\{u,v \}$ of $E$ such that $d(u,v) >
2$ in $F\gamma$, and (\ref{uthreeb3}) holds.  
\item[(B)] $\diam F\gamma = 2$ and
\beq
\label{uthreeb5}
F\gamma \cong K_{2,n} \Rw |E \setminus F| = 1.
\eeq
\ei

Indeed, assume that (A) holds. Clearly, the unique edge that can
be added to the graph and keep its girth above 3 is $u \edge v$.
On the other hand, since (\ref{uthreeb3}) holds, the possibility of
removal of an smi point is excluded. Thus
$(F,E)$ is sji in this case.

Assume now that (B) holds. We cannot remove a 2-set from $F$, since
adding an edge to a graph of diameter 2 brings along girth 3. On the
other hand, having an option on removing an smi point would imply the
existence of two points of degree $\geq |E|-2$, which implies $F\gamma
\cong K_{2,n}$. But in view of (\ref{uthreeb5}), only one of these
points can be present on $F$. Thus $(L,E)$ is sji also in this case.

Conversely, assume that $(F,E)$ is sji. Suppose first that $\diam F\gamma =
2$. As remarked before, we cannot remove a 2-set from $F$, and smi
points correspond to degree $\geq |E|-2$. Therefore there is at most
one such vertex in $F$. Since $K_{2,n}$ has two, then (\ref{uthreeb5})
holds and we fall into case (B).

Finally, assume that $\diam F\gamma > 2$. Then there exist $u,v \in E$
at distance 3 in $F\gamma$, and adding an edge $u \edge v$ does not
spoil (\ref{uthreeb1}). Since $(L,E)$ is sji, then the pair $u,v$ is
unique. Similarly to the characterization of the minimal case,
(\ref{uthreeb3}) must hold to prevent removal of an smi
point. Therefore (A) holds.

We prove next that 
\beq
\label{tojuly}
\mindeg U_{3,2b} = b(b-1)
\eeq
holds for every $b \geq 3$. Indeed, assume that $M$ is an $R \times E$
boolean representation of $U_{3,2b}$ with minimum
degree. By Proposition \ref{stack}(ii), 
we can add all the boolean sums of rows in $M$ and have still a
boolean representation of $U_{3,2b}$, and we can even add a row of
zeroes (we are in fact building the matrix $M^{\nu} \in \M$ from Section
\ref{mvsl}). Now by Proposition \ref{mtol} we have $M^{\nu} = M(L,E)$
for some $(L,E) \in \BR \, U_{3,2b}$, and so $F = (L,E)\theta$ satisfies
$\gth F\gamma > 3$ by (\ref{uthreeb1}). By Tur\'an's Theorem
\cite[Theorem 7.1.1]{Die}, the maximum number of edges in a triangle-free graph
with $2b$ vertices is reached by the complete bipartite graph
$K_{b,b}$ which has $b^2$ edges. Therefore $F\gamma$ has at most $b^2$
edges. Since $2^{E}$ has $\binom{2b}{2} = b(2b-1)$ 2-sets, it follows
that $F$ has at least $b(2b-1) - b^2 = b(b-1)$ 2-sets. Since the
2-sets repesent necessarily smi elements of $M^{\nu}$, it follows that
$M = \wh{M}(L,E)$ has at least $b(b-1)$ elements and so $\mindeg
U_{3,2b} \geq b(b-1)$. Equality is now realized through $F\gamma =
K_{b,b}$. Note that in this case no vertex has degree $\geq |E|-2$,
hence all the points are meets of closed 2-sets and the smi rows of
the matrix are precisely the $b(b-1)$ rows defined by the complement
graph of $K_{b,b}$. Therefore $\mindeg U_{3,2b} = b(b-1)$.

With respect to the odd case, we show that
\beq
\label{tojuly1}
\mindeg U_{3,2b+1} = b^2
\eeq
holds for every $b \geq 3$.
The argument is similar to the the proof of (\ref{tojuly}).
By Tur\'an's Theorem
\cite[Theorem 7.1.1]{Die}, the maximum number of edges in a triangle-free graph
with $2b+1$ vertices is reached by the complete bipartite graph
$K_{b,b+1}$ which has $(b+1)b$ edges. Therefore $F\gamma$ has at most $(b+1)b$
edges. Since $2^{E}$ has $\binom{2b+1}{2} = (2b+1)b$ 2-sets, it follows
that $F$ has at least $(2b+1)b - (b+1)b = b^2$ 2-sets. Therefore
(\ref{tojuly1}) holds.

It is now a simple exercise, for instance, to check that the minimal
representations of $U_{3,6}$ correspond (up to permutation of
vertices) to the graphs
$$
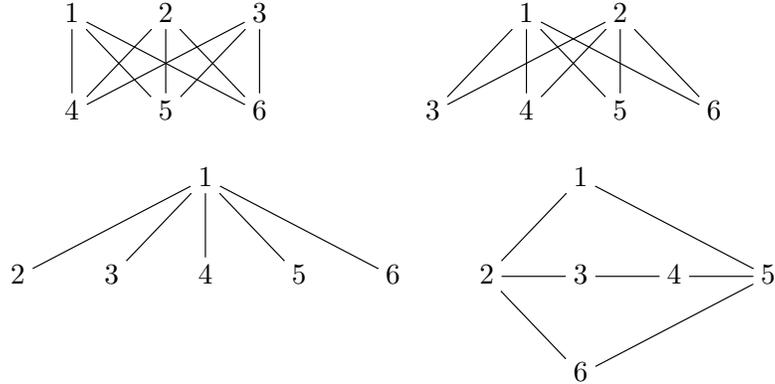
\xymatrix{
1 \ar@{-}[d] \ar@{-}[dr] \ar@{-}[drr] & 2 \ar@{-}[dl] \ar@{-}[d] \ar@{-}[dr] & 3
\ar@{-}[dll] \ar@{-}[dl] \ar@{-}[d] &&& 1 \ar@{-}[dl] \ar@{-}[d] \ar@{-}[dr]
\ar@{-}[drr]  & 2 \ar@{-}[dll] \ar@{-}[dl] \ar@{-}[d] \ar@{-}[dr] &\\
4 & 5 & 6 && 3 & 4 & 5 & 6
}$$
$$\xymatrix{
&&1 \ar@{-}[d] \ar@{-}[dr] \ar@{-}[drr] \ar@{-}[dll] \ar@{-}[dl] &&&&
1 \ar@{-}[dl] \ar@{-}[drr] &&\\ 
2 & 3 & 4 & 5 & 6 & 2 \ar@{-}[r] & 3 \ar@{-}[r] & 4 \ar@{-}[r] & 5 \\
&&&&&& 6 \ar@{-}[ul] \ar@{-}[urr] &&
}$$
and to $F_1,F_2,F_3 \in \fisflats(E,H)$ given respectively by
\bi
\item[]
$F_1 = \{ E, 12, 13,23,45,46,56,1,2,3,4,5,6, \emptyset\}$;
\item[]
$F_2 = \{ E, 12, 34, 35,36,45,46,56,1,3,4,5,6, \emptyset\}$;
\item[]
$F_3 = \{ E, 23,24,25,26, 34, 35, 36,45,46,56,2,3,4,5,6, \emptyset\}$;
\item[]
$F_4 = \{ E, 13,14,16,24,25, 35, 36,46,1,2,3,4,5,6, \emptyset\}$.
\ei
The corresponding lattices are now
$$
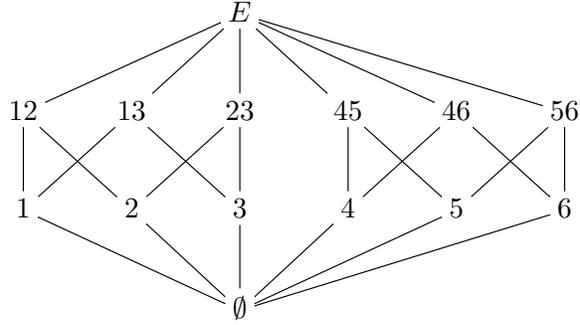
\xymatrix{
&&E \ar@{-}[dll] \ar@{-}[dl] \ar@{-}[d] \ar@{-}[dr] \ar@{-}[drr]
\ar@{-}[drrr] &&&\\
12 \ar@{-}[d] \ar@{-}[dr] & 13 \ar@{-}[dl] \ar@{-}[dr] & 23 \ar@{-}[dl]
\ar@{-}[d] & 45 \ar@{-}[d] \ar@{-}[dr] & 46 \ar@{-}[dl] \ar@{-}[dr] & 56
\ar@{-}[dl] \ar@{-}[d] \\
1 & 2 & 3 & 4 & 5 & 6 \\
&&\emptyset \ar@{-}[ull] \ar@{-}[ul] \ar@{-}[u] \ar@{-}[ur] \ar@{-}[urr]
\ar@{-}[urrr] &&&
}
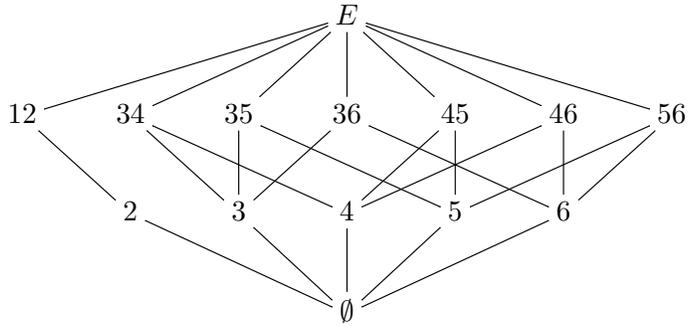
$$
$$\xymatrix{
&&&E \ar@{-}[dlll] \ar@{-}[dll] \ar@{-}[dl] \ar@{-}[d] \ar@{-}[dr] \ar@{-}[drr]
\ar@{-}[drrr] &&&\\
12 \ar@{-}[dr] & 34 \ar@{-}[drr] \ar@{-}[dr] & 35 \ar@{-}[d]
\ar@{-}[drr] & 36 \ar@{-}[dl] \ar@{-}[drr] & 45 \ar@{-}[d] \ar@{-}[dl] &
46 \ar@{-}[d] \ar@{-}[dll] & 56 \ar@{-}[dll] \ar@{-}[dl] \\
& 2 & 3 & 4 & 5 & 6 &\\
&&&\emptyset \ar@{-}[ull] \ar@{-}[ul] \ar@{-}[u]
\ar@{-}[ur] \ar@{-}[urr] &&&
}
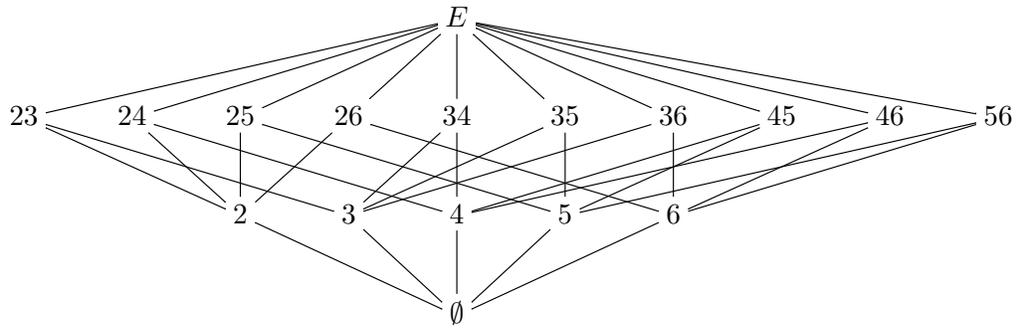
$$
$$\xymatrix{
&&&&E \ar@{-}[dllll] \ar@{-}[dlll] \ar@{-}[dll] \ar@{-}[dl] \ar@{-}[d]
\ar@{-}[dr] \ar@{-}[drr] 
\ar@{-}[drrr] \ar@{-}[drrrr] \ar@{-}[drrrrr] &&&&&\\
23 \ar@{-}[drr] \ar@{-}[drrr] & 24 \ar@{-}[drrr] \ar@{-}[dr] & 25
\ar@{-}[d] \ar@{-}[drrr] & 26 \ar@{-}[dl] 
\ar@{-}[drrr] & 34 \ar@{-}[dl] \ar@{-}[d] & 35 \ar@{-}[dll] \ar@{-}[d] &
36 \ar@{-}[d] \ar@{-}[dlll] & 45 \ar@{-}[dll] \ar@{-}[dlll] & 46 \ar@{-}[dllll]
\ar@{-}[dll] & 56 \ar@{-}[dlll] \ar@{-}[dllll] \\
&& 2 & 3 & 4 & 5 & 6 &&&\\
&&&&\emptyset \ar@{-}[ull] \ar@{-}[ul] \ar@{-}[u]
\ar@{-}[ur] \ar@{-}[urr] &&&&&
}
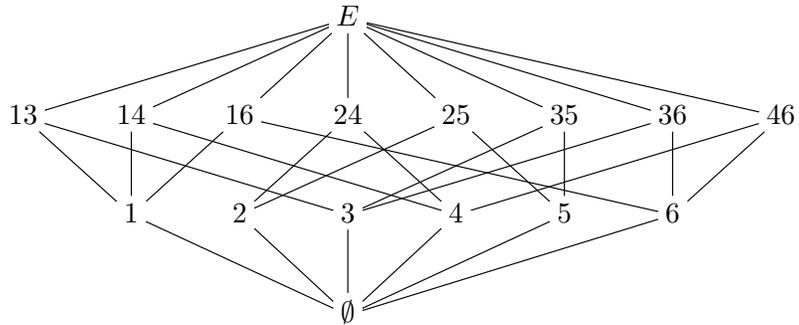
$$
$$\xymatrix{
&&&E \ar@{-}[dlll] \ar@{-}[dll] \ar@{-}[dl] \ar@{-}[d] \ar@{-}[dr] \ar@{-}[drr]
\ar@{-}[drrr] \ar@{-}[drrrr] &&&&\\
13 \ar@{-}[drrr] \ar@{-}[dr] & 14 \ar@{-}[d] \ar@{-}[drrr] & 16 \ar@{-}[dl]
\ar@{-}[drrrr] & 24 \ar@{-}[dl] \ar@{-}[dr] & 25 \ar@{-}[dll] \ar@{-}[dr] &
35 \ar@{-}[d] \ar@{-}[dll] & 36 \ar@{-}[d] \ar@{-}[dlll] & 46 \ar@{-}[dlll]
\ar@{-}[dl] \\
& 1 & 2 & 3 & 4 & 5 & 6 &\\
&&&\emptyset \ar@{-}[ull] \ar@{-}[ul] \ar@{-}[u]
\ar@{-}[ur] \ar@{-}[urr] 
\ar@{-}[urrr] &&&&
}$$

The non minimal sji representations of $U_{3,6}$ can be easily computed. 
In fact, it is easy to see that if (A) holds, then by adding an edge $u
\edge v$ to the graph $F\gamma$ we get a graph of diameter 2 and still
girth $> 3$. The converse is not true, but a brief analysis of all the
possible removals of one edge from a minimal case graph to reach (A)
gives us all such sji representations.

Those of type (B) are obtained by adding the seventh point to the
minimal representation given by $K_{1,5}$ (the other types already
have the seven points or are excluded by (\ref{uthreeb5}).

Therefore the graphs corresponding to the sji representations of type
(A) are
$$
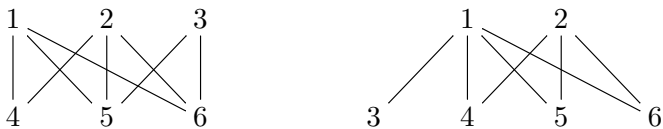
\xymatrix{
1 \ar@{-}[d] \ar@{-}[dr] \ar@{-}[drr] & 2 \ar@{-}[dl] \ar@{-}[d] \ar@{-}[dr] & 3
\ar@{-}[dl] \ar@{-}[d] &&& 1 \ar@{-}[dl] \ar@{-}[d] \ar@{-}[dr]
\ar@{-}[drr]  & 2 \ar@{-}[dl] \ar@{-}[d] \ar@{-}[dr] &\\
4 & 5 & 6 && 3 & 4 & 5 & 6
}$$
obtained by removing an edge from $K_{3,3}$ and $K_{2,4}$,
respectively. Adding the case (B) representation, we obtain types
\bi
\item[]
$F_5 = \{ E, 12, 13,23,34, 45,46,56,1,2,3,4,5,6, \emptyset\}$;
\item[]
$F_6 = \{ E, 12, 23, 34, 35,36,45,46,56,1,3,4,5,6, \emptyset\}$;
\item[]
$F_7 = \{ E, 23,24,25,26, 34, 35, 36,45,46,56,1,2,3,4,5,6, \emptyset\}$.
\ei
The corresponding lattices are 
$$
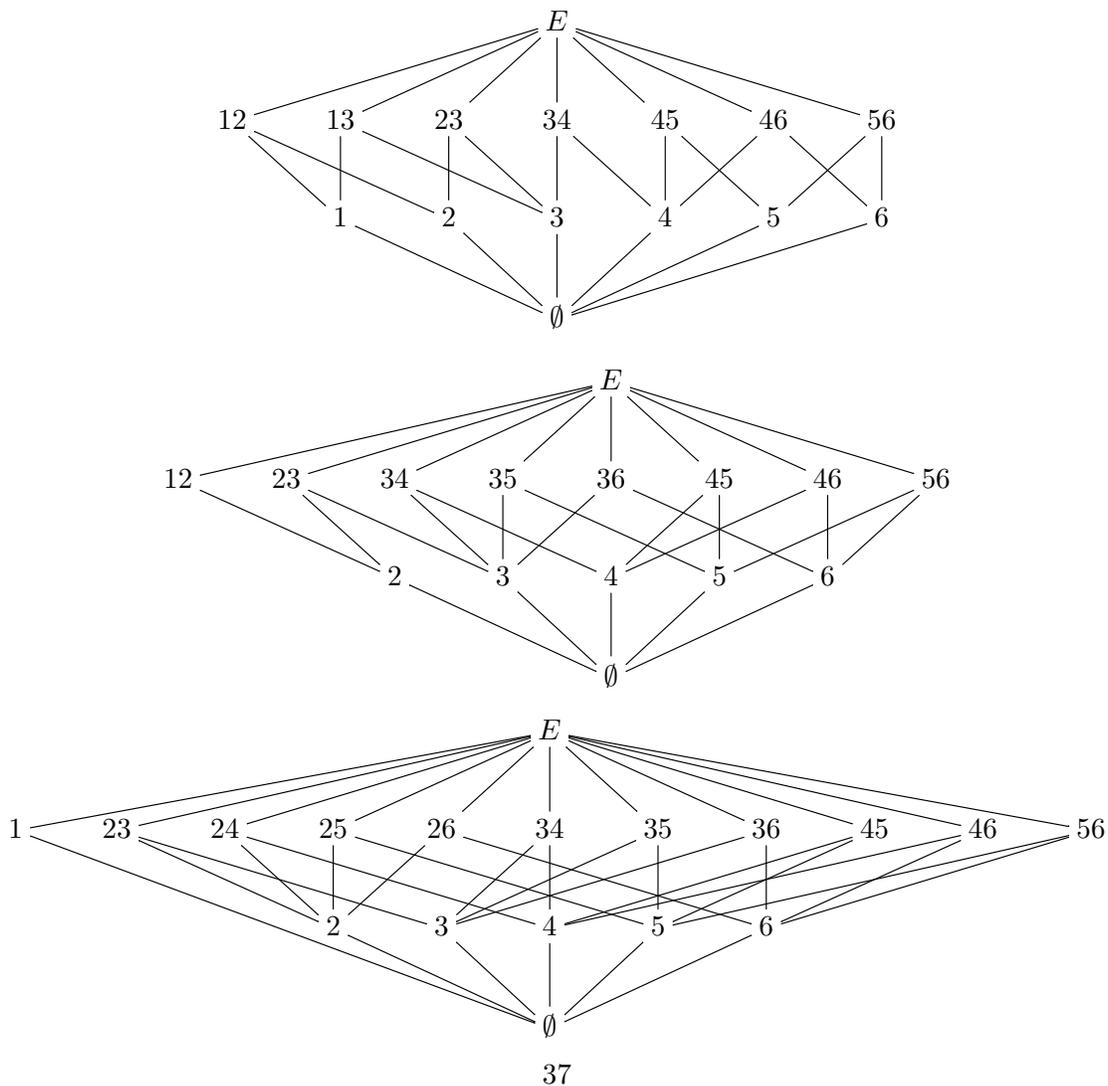
\xymatrix{
&&&E \ar@{-}[dlll] \ar@{-}[dll] \ar@{-}[dl] \ar@{-}[d] \ar@{-}[dr] \ar@{-}[drr]
\ar@{-}[drrr] &&&\\
12 \ar@{-}[dr] \ar@{-}[drr] & 13 \ar@{-}[d] \ar@{-}[drr] & 23 \ar@{-}[d]
\ar@{-}[dr] & 34 \ar@{-}[d] \ar@{-}[dr] & 45 \ar@{-}[d] \ar@{-}[dr] &
46 \ar@{-}[dl] \ar@{-}[dr] & 56 
\ar@{-}[dl] \ar@{-}[d] \\
& 1 & 2 & 3 & 4 & 5 & 6 \\
&&&\emptyset \ar@{-}[ull] \ar@{-}[ul] \ar@{-}[u] \ar@{-}[ur] \ar@{-}[urr]
\ar@{-}[urrr] &&&
}$$
$$\xymatrix{
&&&&E \ar@{-}[dllll] \ar@{-}[dlll] \ar@{-}[dll] \ar@{-}[dl] \ar@{-}[d]
\ar@{-}[dr] \ar@{-}[drr] 
\ar@{-}[drrr] &&&\\
12 \ar@{-}[drr] & 23 \ar@{-}[drr] \ar@{-}[dr] &  34 \ar@{-}[drr]
\ar@{-}[dr] & 35 \ar@{-}[d] 
\ar@{-}[drr] & 36 \ar@{-}[dl] \ar@{-}[drr] & 45 \ar@{-}[d] \ar@{-}[dl] &
46 \ar@{-}[d] \ar@{-}[dll] & 56 \ar@{-}[dll] \ar@{-}[dl] \\
&& 2 & 3 & 4 & 5 & 6 &\\
&&&&\emptyset \ar@{-}[ull] \ar@{-}[ul] \ar@{-}[u]
\ar@{-}[ur] \ar@{-}[urr] &&&
}$$
$$\xymatrix{
&&&&&E \ar@{-}[dlllll] \ar@{-}[dllll] \ar@{-}[dlll] \ar@{-}[dll]
\ar@{-}[dl] \ar@{-}[d] 
\ar@{-}[dr] \ar@{-}[drr] 
\ar@{-}[drrr] \ar@{-}[drrrr] \ar@{-}[drrrrr] &&&&&\\
1 & 23 \ar@{-}[drr] \ar@{-}[drrr] & 24 \ar@{-}[drrr] \ar@{-}[dr] & 25
\ar@{-}[d] \ar@{-}[drrr] & 26 \ar@{-}[dl] 
\ar@{-}[drrr] & 34 \ar@{-}[dl] \ar@{-}[d] & 35 \ar@{-}[dll] \ar@{-}[d] &
36 \ar@{-}[d] \ar@{-}[dlll] & 45 \ar@{-}[dll] \ar@{-}[dlll] & 46 \ar@{-}[dllll]
\ar@{-}[dll] & 56 \ar@{-}[dlll] \ar@{-}[dllll] \\
&&& 2 & 3 & 4 & 5 & 6 &&&\\
&&&&&\emptyset \ar@{-}[uulllll] \ar@{-}[ull] \ar@{-}[ul] \ar@{-}[u]
\ar@{-}[ur] \ar@{-}[urr] &&&&&
}$$

It is
easy to count $20 +15+6+180 = 221$ minimal lattice representations for
$U_{3,6}$ only (but they reduce to $1+1+1+1 = 4$ in the alternative
counting of Remark \ref{spcorps})! The sji's (including the minimal
cases) amount to $221+180+120+6 = 527$ and $4+1+1+1 = 7$ in both
countings. 
Note also that $\mindeg U_{3,6} = 6$ by (\ref{tojuly}).

Note that the lattices in the examples in which $E \subseteq F$, after
removal of the top and bottom elements, are essentially the {\em Levi
  graphs} of the graphs $F\gamma$. The Levi graph of $F\gamma$ can be
obtained by introducing a new 
vertex at the midpoint of every edge (breaking thus the original edge
into two), and the connection to the lattice is established by
considering that each of the new vertices lies above its two adjacent
neighbours. 

Note also that famous graphs of
girth $> 3$ and diameter 2 such as the {\em Petersen graph} \cite{W9}
turn out to encode minimal respresentations via the function $\gamma$
(in $U_{3,10}$, since the Petersen graph has 10 vertices).

\section{Additions}

This section contains results which are relevant to the theory but
were not needed for the main sections.

\subsection{Rank functions}

Let $(E,H)$ be a hereditary collection. The {\em rank function} $r_H:
2^E \to \N$ is defined by
$$Xr_H = \max\{ |I|: I \in 2^X \cap H\}.$$
Given a function $f:2^E \to \N$, consider the following axioms for all
$X,Y \subseteq E$:
\bi
\item[(A1)] $X \subseteq Y \hspace{.3cm}\Rw \hspace{.3cm}
  Xf \leq Yf$;
\item[(A2)] $\exists I \subseteq X : |I| = If = Xf$;
\item[(A3)] $(Xf = |X| \; \wedge\; Y \subseteq X) \hspace{.3cm}\Rw
  \hspace{.3cm} Yf = |Y|$. 
\ei
It is easy to see that the three axioms are independent.

\bp
\label{arf}
Given a function $f:2^E \to \N$, the following conditions are
equivalent:
\bi
\item[(i)] $f = r_H$ for some hereditary collection $(E,H)$;
\item[(ii)] $f$ satisfies axioms (A1)--(A3).
\ei
\ep

\proof
(i) $\Rw$ (ii). It follows immediately from the equivalence
$$Xr_H = |X| \iff X \in H.$$

(ii) $\Rw$ (i). Let $H = \{ I \subseteq E: If = |I|\}$. By (A3), $H$
is closed under taking subsets. Taking $X = \emptyset$ in (A2), we get
$\emptyset f = 0$, hence $\emptyset \in H$ and so $(E,H)$ is a
hereditary collection. Now, for every $X \in E$, we have 
$$Xr_H = \max\{ |I|: I \in 2^X \cap H\} 
= \max\{ |I|: I \subseteq X ,\; If = |I|\}.$$
By (A2), we get $Xr_H \geq Xf$, and $Xr_H \leq Xf$ follows from
(A1). Hence $f = r_H$ as required.
\qed

We collect next some elementary properties of rank functions:

\bp
\label{proprf}
Let $(E,H)$ be a hereditary collection and let $X,Y \subseteq E$. Then:
\bi
\item[(i)] $Xr_H \leq |X|$;
\item[(ii)] $Xr_H + Yr_H \geq (X\cup Y)r_H$;
\item[(iii)] $Xr_H + Yr_H \geq (X\cup Y)r_H + (X\cap Y)r_H$ if some
  independent subset of maximum size of $X\cap Y$ can be extended to
  some  independent subset of maximum size of $X\cup Y$;
\item[(iv)] $Xr_H + Yr_H \geq (X\cup Y)r_H + (X\cap Y)r_H$ if $(E,H)$
  is a matroid.
\ei
\ep

\proof
(i) By (A2).

(ii) Assume that $(X\cup Y)r_H = |I|$ with $I \in 2^{X\cup Y} \cap
H$. Then $I \cap X, I \cap Y \in H$ and so
$(X\cup Y)r_H = |I| \leq |I\cap X| + |I\cap Y| \leq Xr_H + Yr_H$.

(iii) We may assume that $(X\cup Y)r_H = |I|$ and $(X\cap Y)r_H = |I
\cap X \cap Y|$. It follows that 
$(X\cup Y)r_H + (X\cap Y)r_H = |I| + |I
\cap X \cap Y| = |I\cap X| + |I\cap Y| \leq Xr_H + Yr_H$.

(iv) This is well known, but we can include a short deduction from
(iii).

Let $K \subseteq L \subseteq E$, and assume that
  $J$ is an independent subset of maximum size of $K$. Let $J'$ be a
  maximal independent subset of $L$ containing $J$. If $(E,H)$
  is a matroid, it follows from the exchange property that $|J'| =
  Lr_H$. Now we apply part (iii) to $K = X\cap Y$ and $L = X\cup Y$.
\qed

\bp
\label{epflats}
Let $(E,H)$ be a hereditary collection of rank $r$.
\bi
\item[(i)] If $X,Y \in$ {\rm Fl}$(E,H)$ and $Xr_H = Yr_H$, then
$$X \subseteq Y \hspace{.3cm} \iff \hspace{.3cm} X = Y.$$
\item[(ii)] $E$ is the unique flat of rank $r$.
\ei
\ep

\proof
(i) Let $X,Y \in \flats(E,H)$. Assume that $X \subseteq Y$ and let $I
\subseteq X$ satisfy $I \in H$ and $|I| = Xr_H = Yr_H$. If $p \in Y
\setminus X$, then $X$ closed yields $I \cup \{ p \} \in H$ and $Yr_H
> |I| = Xr_H$, a contradiction. Therefore $X = Y$ and (i) holds.

(ii) By part (i).
\qed

It follows that the flats of rank $r-1$ are maximal in $\flats(E,H)
\setminus \{ E \}$. Such flats are called {\em hyperplanes}.

The following result relates the rank function with the closure
operator $\cl$ induced by a (simple) hereditary collection.

\bp
\label{ranklat}
Let $(E,H)$ be a simple hereditary collection admitting a boolean
representation and let $X \subseteq E$. Write $L =$ {\rm Fl}$(E,H)$.
Then $Xr_H$ is the maximum
$k$ such that (\ref{hcl1}) holds for some $x_1,\ldots,x_k \in X$, and 
the maximum
$k$ such that (\ref{hcl2}) holds for some $x_1,\ldots,x_k \in X$.
\ep

\proof
The first equality follows from Proposition \ref{hcl} and the
definition of $r_H$. 
The second follows from Theorem \ref{rephc} and Lemma \ref{clopes}(ii).
\qed

\subsection{Paving hereditary collections}
\label{spav}

A hereditary collection of rank $r > 2$ is said to be {\em paving} if
it has no circuits of size $< r$ (or equivalently, of rank less than
$r - 1$). 
%Note that if $C$ is a circuit, then $Cr_H = |C|-1$.

\bl
\label{paving}
Let $(E,H)$ be a hereditary collection of rank $r > 2$. Then the
following conditions are equivalent:
\bi
\item[(i)] $(E,H)$ is paving;
\item[(ii)] $P_{r-1}(E) \subseteq H$;
\item[(iii)] $P_{r-2}(E) \subseteq$ {\rm Fl}$(E,H)$.
\ei
\el

\proof
(i) $\Rw$ (ii). Since every dependent subset of $(E,H)$ must contain a
circuit.

(ii) $\Rw$ (i). Trivial.

(ii) $\Rw$ (iii). By Proposition \ref{trivia}(ii).

(iii) $\Rw$ (i). Suppose that $C$ is a circuit of rank $< r$ and let
$x \in C$. Then $|C \setminus \{ x \}| \leq r-2$, hence $C \setminus
\{ x \}$ is closed and $C \setminus \{ x \} \in H$ yields $C \in H$, a
contradiction. Thus $(E,H)$ is paving.
\qed

Next we provide a simple characterization of boolean representable
paving hereditary collections:

\bp
\label{brp}
Let $(E,H)$ be a paving hereditary collection of rank $r$. Then the
following conditions are equivalent:
\bi
\item[(i)] $(E,H)$ is boolean representable;
\item[(ii)] $\forall X \in H\hspace{.3cm} |X| = r \Rw \exists x \in X:
  x \notin$ {\rm Cl}$(X \setminus \{ x \})$;
\item[(iii)] $\forall X \in H\hspace{.3cm} |X| = r \Rw \exists x \in
  X:$ {\rm Cl}$(X \setminus \{ x \}) \neq E$.
\ei
\ep

\proof
(i) $\Rw$ (ii). By Proposition \ref{hcl}.

(ii) $\Rw$ (iii). Immediate.

(iii) $\Rw$ (i). By Theorem \ref{indecl}, it suffices to show that
every $X \in H$ admits an enumeration $x_1,\ldots , x_k$ such that
$$\clos(x_1,\ldots, x_k) \supset \clos(x_2,\ldots, x_k) \supset \ldots
\supset \clos(x_k).$$
By condition (iii) in Lemma \ref{paving}, this condition is satisfied
if $|X| < r$. Hence we may assume that $|X| = r$ and we only need to
show that there exists some enumeration $x_1,\ldots , x_r$ of $X$ such
that 
$$\cl X \supset \clos(X \setminus \{ x_1 \}) \supset \{ x_3,
\ldots,x_r \} \supset \ldots \supset \{ x_r \}.$$
Since $\cl X = E$ by Proposition \ref{newtrivia}, condition (iii)
yields the required inclusion.
\qed

In connection with Proposition \ref{brp}, we can mention several
equivalent characterizations of matroids 
among paving hereditary collections: 

\bp
Let $(E,H)$ be a paving hereditary collection of rank $r$. Then the
following conditions are equivalent:
\bi
\item[(i)]
$(E,H)$ is a matroid;
\item[(ii)]
every $(r-1)$-subset of $E$ is contained in a unique hyperplane;
\item[(iii)]
if $X$ is an $(r-1)$-subset of $E$, then {\rm Cl}$\, X \neq E$.
\ei
\ep

\proof
(i) $\Rw$ (ii). By \cite[Proposition 2.1.21]{Oxl}.

(ii) $\Rw$ (iii). Immediate.

(iii) $\Rw$ (i). First, note that $(E,H)$ is boolean representable by
Proposition \ref{brp}. Let $I,J \in H$ with $|I| = |J| + 1$. We must show
that $J \cup \{ i \} \in H$ for some $i \in I \setminus J$. Since 
$P_{r-1}(E) \subseteq H$ by Lemma \ref{paving}, we may assume that
$|J| = r-1$. Since $\cl J \neq E$ and $\cl I = E$ by Proposition
\ref{newtrivia}, we get $I \not\subseteq \cl J$. Take $i \in I
\setminus \cl J$. 
Since $P_{r-2}(E) \subseteq \flats (E,H)$ by  Lemma \ref{paving} and
$\cl J \subset \clos(J \cup \{ i \})$, it follows from Proposition
\ref{hcl} that $J \cup \{ i \}
\in H$. Therefore $(E,H)$ is a matroid.
\qed

Note that $(E,H)$ being boolean representable and all its bases having
rank $r$ does {\em not} imply that $(E,H)$ is a matroid, a
counterexample being provided by $E = \hat{6}$ and $H = P_4(E)
\setminus \{ 2456, 3456 \}$.

\subsection{Boolean operations}

Boolean representability behaves badly with respect to intersection
and union, as we show next.

First, we recall a well-known fact: every hereditary collection
$(E,H)$ is the intersection of matroids on $E$, namely the intersection of the
matroids $M_X$ over all circuits $X$ of $(E,H)$, where $M_X$ is the matroid
consisting of all subsets of $E$ not containing $X$ \cite{Oxl}. Since
all simple matroids are boolean representable by Theorem \ref{boom},
it follows that all simple hereditary collections are the intersection
of boolean representable hereditary collections. Therefore  boolean
representable hereditary collections are not closed under intersection.

\be
\label{unio}
Let $E = \hat{6}$ and $J_1 = P_3(E) \setminus \{ 123, 125, 135, 235,
146, 246, 346, 456 \}$, $J_2 = P_2(E) \cup \{ 123, 124, 125, 126 \}$. Then
$(E,J_1)$ and $(E,J_2)$ are both boolean 
representable hereditary collections, but $(E,J_1 \cup J_2)$ is not.
\ee

It is easy to check that $1235 \in
\flats(E,J_1)$. Since $|xyz \cap 1235| = 2$ for every $xyz \in J_1$,
it follows from Proposition \ref{brp} that $(E,J_1)$ is boolean 
representable. Similarly, since $12 \in \flats(E,J_2)$, we show that
that $(E,J_2)$ is boolean representable.

Now $J_1 \cup J_2 = P_3(E) \setminus \{ 135, 235,
146, 246, 346, 456 \}$ and it is straightforward to check that $\cl 13
= \cl 15 = \cl 35 = 1235$. By Proposition \ref{brp}, $(E,J_1\cup J_2)$
is not boolean representable.

\smallskip

However, closure under union can be satisfied in particular circumstances:

\bp
\label{unith}
Let $(E,H_1),(E,H_2)$ be simple boolean representable hereditary
collections of rank 3. If
\beq
\label{unith1}
X \in {\rm Fl}(E,H_i)\setminus \{ E \} \Rw |X| \leq 3
\eeq
holds for $i=1,2$, then $(E,H_1\cup H_2)$ is boolean representable.
\ep

\proof
For $i = 1,2$, write $F_i = \flats(E,H_i)$ and let $L_i = \{ X \in
F_i: |X| = 3 \}$. We define also the set of {\em potential lines} 
$$P_i = \{ X \subseteq E: |X| = 3 \mbox{ and $|X \cap Y| \leq 1$ for
  every }Y \in F_i\setminus \{ E \} \}.$$
It is easy to see that $H_i \cap P_i = \emptyset$. Indeed, suppose
that $xyz \in H_i \cap P_i$. Since $(E,H_i)$, we may apply Proposition \ref{brp}
and assume, without loss of generality, that there exists some $Y \in
F_i\setminus \{ E \}$ containing $xy$. This contradicts $xyz \in P_i$,
hence $H_i \cap P_i = \emptyset$.

Let $F = \flats(E,H_1 \cup H_2)$ and define
\bi
\item[] $W = (L_1 \cap L_2) \cup (L_1 \cap P_2) \cup (P_1 \cap L_2)$; 
% \item[] $W' = W \cup (P_1 \cap P_2)$;
\item[] $W' = \{ X \subseteq E: |X| = 2$ and $X \cup \{ q \} \notin
  (L_1 \cup P_1)\cap(L_2 \cup P_2)$ for every $q \in E \}$.
\ei
We claim that 
\beq
\label{unith2}
W \cup W' \subseteq F.
\eeq
Let $xyz \in W$. We may assume that $xyz \in L_1$. Suppose first that
$xyz \in L_2$. If $I \subseteq xyz$, $I \in H_1 \cup H_2$ and $p
\notin xyz$, then $I \in H_i$ for some $i$ and so $xyz \in
L_i$ yields $I \cup \{ p \} \in H_i \subseteq H_1 \cup H_2$ as
required.

Hence we may assume that $xyz \in P_2$. Let $I \subseteq xyz$ and $p
\notin xyz$. If $I \in H_1$, all is similar to the preceding case,
hence we assume that $I \in H_2 \setminus H_1$. Since $H_2 \cap P_2 =
\emptyset$, we have $|I| \leq 2$ and so $I \in H_1$. If  $p
\notin xyz$, It follows from $xyz \in F_1$ that $I \cup \{ p \} \in
H_1 \subseteq H_1 \cup H_2$. Thus $W \subseteq F$.

Assume now that $xy \in W'$. Let $z \notin xy$. It suffices to show
that $xyz \in H_1 \cup H_2$. 
Since $xy \in W'$, we have $xyz \notin L_i \cup
P_i$ for some $i$. Thus $|xyz \cap Y| = 2$ for some $Y \in
F_i\setminus \{ E \}$, and $Y$ closed yields $xyz \in H_i \subseteq
H_1 \cup H_2$ as required. Therefore $xy \in F$ and (\ref{unith2})
holds.

Now let $\cl X$ (respectively $\clos_iX$) denote the closure of $X
\subseteq E$ in $(E,H_1 
\cup H_2)$ (respectively $(E,H_i)$). 
Let $xyz \in H_1 \cup H_2$. By Proposition \ref{brp}, we must show
that $\cl (xy) \neq E$ or $\cl (xz) \neq E$ or $\cl (yz) \neq E$. We may
assume that $xyz \in 
H_1$. By Proposition \ref{brp}, we may assume also that $\clos_1(xy) \neq
E$. 

Since $xy \in W'$ implies $xy \in F$ by (\ref{unith2}), we have
that $xyq \in (L_1 \cup P_1)\cap(L_2 \cup P_2)$ for some $q \in E$. If
$xyq \in W \subseteq F$, we immediately get $\cl (xy) \neq E$. It
remains to consider the case 
$xyq \in P_1 \cap P_2$. Since $|xyq \cap \clos_1(xy)| \geq 2$ and
$\clos_1(xy) \neq E$, we reach a contradiction. Thus $\cl (xy) \neq E$
and $(E,H_1 \cup H_2)$ is boolean representable.
\qed

\subsection{Truncation}

Given a hereditary collection $(E,H)$ and $k \geq 0$, the $k$-{\em
  truncation} of $(E,H)$ is the hereditary collection $(E,H_k)$
defined by $H_k = \{ X \subseteq E : |X| \leq k \}$. 

\bp
\label{tru}
Let  $(E,H)$ be a hereditary collection and let $k \geq 0$. Then:
\bi
\item[(i)] {\rm Fl}$\,(E,H_k) \subseteq$ {\rm Fl}$\,(E,H)$;
\item[(ii)] If $X \subset E$, then $X \in$ {\rm Fl}$\,(E,H_k)$ if and
  only if $X \in$ {\rm Fl}$\,(E,H)$ and $X$ does not contain a basis of
  $H_k$.
\ei
\ep

\proof
It suffices to prove (ii). Let $X \subset E$. Assume first that $X \in
\flats(E,H_k)$. By Proposition \ref{newtrivia}, $X$ does not
contain a
basis of $(E,H_k)$. Let $p \in E \setminus X$ and let $I \subseteq X$ be
such that $I \in H$. Since $I$ is not a
basis of $(E,H_k)$, we have $|I| < k$ and so $I \in H_k$. Now $X \in
\flats(E,H_k)$ yields $I \cup \{ p \} \in H_k \subseteq H$. Therefore $X \in
\flats(E,H)$.

Conversely, assume that $X \in
\flats(E,H)$ and $X$ does not contain a basis of $(E,H_k)$. Let $p \in
E \setminus X$ and let $I \subseteq X$ be such that $I \in H_k$. Since
$H_k \subseteq H$ and $X \in
\flats(E,H)$, we get $I \cup \{ p \} \in H$. But $I$ is not a basis of
$(E,H_k)$, hence $|I| < k$ and so $I \cup \{ p \} \in H_k$. Thus $X \in
\flats(E,H_k)$ as required.
\qed

The next example shows that boolean representability is not preserved
under truncation, even in the simple case.

\be
\label{truno}
Let $E = \hat{6}$ and let $H$ be the hereditary collection defined by
$H = (P_3(E) \setminus \{ 135, 235,
146, 246, 346, 456 \}) \cup \{ 1234, 1236, 1245, 1256\}$. Then $(E,H)$
is boolean representable, but $(E,H_3)$ is not.
\ee

It is easy to check that $P_1(E) \cup \{ 12, 1235\}
\subseteq \flats(E,H)$. By
Theorem \ref{indecl}, to show that $(E,H)$ is boolean 
representable it suffices to show that every $X \in H$ admits
an enumeration $x_1,\ldots , x_k$ satisfying (\ref{hcl1}). We may of
course assume that $|X| > 2$. Hence $X$ cannot contain both 4 and
6. Since $1235$ is 
closed, we may assume that $X \subseteq 1235$. Since we may assume
that $X$ is a 3-set, we are reduced to the cases $X \in \{ 123, 125
\}$. Now $1 \subset 12 \subset 1235$ yields the desired chain of
flats, and so $(E,H)$ is boolean 
representable.

On the other hand, $H_3$ is the collection $J_1 \cup J_2$ of Example
\ref{unio}, already proved not to be boolean 
representable.

% \addcontentsline{toc}{section}{Appendix}
\section{Appendix}

We gather here several results which, although not essential for
obtaining our main results, can help the interested reader to gain
further insight into our approach and methods.

% \addcontentsline{toc}{subsection}{A. Categoric alternatives}
\subsection{Categoric alternatives}
\label{101}

We note that the category FL is isomorphic to some other categories that bring
different viewpoints into our discussions. 

Recall that a structure $(S,+,\cdot,0)$
is a {\em semiring} if: 
\bi
\item
$(S,+,0)$ is a commutative monoid;
\item
$(S,\cdot,0)$ is a semigroup with zero
\item
$x(y+z) = xy + xz$ and $(y+z)x = yx + zx$ for all $x,y,z \in S$.
\ei
The semiring $S$ is {\em idempotent} if $x+x =
x$ for every $x \in S$. It is {\em null} if $xy = 0$ for all $x,y \in
S$. Morphisms and modules over semirings are defined the obvious way
(see \cite[Chapter 9]{RS}).

We introduce the following notation:
\bi
\item[] FICM: the category of {\em finite idempotent commutative monoids}
  together with monoid morphisms;
\item[] FINS: the category of {\em finite idempotent null semirings}
  together with semiring morphisms; 
\item[] FBM: the category of {\em finite unitary right $\BB$-modules}
  together with
  $\BB$-module morphisms.
\ei 

\bp
\label{isocat}
The categories {\rm FL}, {\rm FICM}, {\rm FINS} and {\rm FBM} are isomorphic.
\ep

\proof
It is well known that the functor FL $\to$ FICM defined by $(L,\leq) \mapsto
(L,\vee)$ and identity on arrows defines an isomorphism of
categories. 

Clearly, the forgetful functor FINS $\to$ FICM is also an isomorphism of
categories. The same happens for the forgetful functor FBM $\to$
FICM. Indeed, it is easy to see that each $\BB$-module is necessarily
idempotent since $x = 1x = (1+1)x = 1x+1x = x+x$ holds for every $x
\in S$. On the other hand, each finite idempotent commutative monoid
$(M,+,0)$ determines a unique $\BB$-module structure in $M$ since we
are forced to have $1x = x$ and $0x = 0$ for every $x \in M$, and the
arrows turn out to be the same mappings.
\qed

Now, for each category X $\in \{$ FICM, FINS, FBM $\}$, we define
another category Xg by taking objects of the form $(M,E)$, where $M$
is an object of X and $E \subseteq M \setminus \{ 0 \} $ a generating
set for $M$. For arrows
$\varphi:(M,E) \to (M',E')$, we
require also $E\p \subseteq E' \cup \{ 0 \}$. With straightforward
adaptations, Proposition \ref{isocat} yields

\bc
\label{isocatg}
The categories {\rm FLg}, {\rm FICMg}, {\rm FINSg} and {\rm FBMg} are
isomorphic. 
\ec

The following result, stated for FL and ideals, which is after all our
basic viewpoint in this 
paper, is inspired by standard concepts in semigroup theory \cite{CP},
and therefore by the viewpoint FICM. We say that $I \subseteq L$ is an
{\em ideal} (or downset) if 
$x \leq y \in I$ implies $x \in I$ for all $x,y \in L$.  There is an obvious
dual version of the Rees quotient when we consider the dual notion of
upset. 

Given $(L,E) \in \FLg$ and an ideal $I$ of $L$, the {\em Rees
  quotient} $L/I$ is the 
quotient of $L$ by the congruence $\sim_I$ defined on $L$ by
$$x \sim_I y  \hspace{.5cm} \mbox{if} \hspace{.5cm} x = y \mbox{ or
}x,y \in I.$$
The elements of $L/I$ are the equivalence class $B' = I$ (the bottom
element) and the 
singular equivalence classes $\{ x \}$ $(x \in L \setminus I)$, which
we identify with $x$. The partial ordering of $L$ translates naturally to $L/I$.
% Note that upsets convert into ideals in the isomorphic category FICMg
% and so $L/U$ is the lattice version of the classical Rees quotient
% from semigroup theory (see \cite{CP}).

\bp
\label{rees}
Let $L \in$ {\rm FL} and let $I$ be
an ideal of $L \setminus \{ T \}$. Then 
$L/I \in$ {\rm FLg}.
\ep

\proof
Clearly, $L/I$ inherits a naural $\wedge$-semilattice structure, and
then becomes a lattice with the determined join.
\qed

It is also possible to import the notion of {\em quotient submodule}
(see \cite[Section 9.1]{RS}) from FBM to FL. Let $S$ be a
$\vee$-subsemilattice of a finite lattice $L$. We define a
$\vee$-congruence $\equiv_S$ on $L$ by $x \equiv_S y$ if $(x \vee s) =
(y \vee s')$ for some $s, s' \in S$. Then $L/S$ denotes the quotient
$L/\equiv_S$.

Given a finite lattice $L$, we say that $\xi:L \to L$ is a {\em
  closure operator} if the following axioms hold for all $x,y \in L$:
\bi
\item[(C1)] $x \leq x\xi$;
\item[(C2)] $x \leq y \hspace{.5cm} \Rw \hspace{.5cm} x\xi \leq y\xi$;
\item[(C3)] $x\xi = (x\xi)\xi$.
\ei

The next proposition summarizes some of the properties of closure
operators (see \cite[Subsection I.3.12]{Gra}):

\bp
\label{cla}
Let $L$ be a lattice, let $\xi:L \to L$ be a closure operator and
let $S$ be a $\wedge$-subsemilattice of $L$. Then:
\bi
\item[(i)] $(x \vee y)\xi = (x\xi \vee y\xi)\xi$ for all $x,y \in L$.
\item[(ii)] $L\xi$ is a $\wedge$-subsemilattice of $L$ and 
  constitutes a lattice under the determined join $(x\xi \vee' y\xi) =
  (x\xi \vee y\xi)\xi$.
\item[(iii)] $S$ induces a
  closure operator $\xi_S:L \to L$ defined by $x\xi_S = \wedge\{ y \in
  S \mid y \geq x \}$.
\item[(iv)] $\xi_{L\xi} = \xi$ and $L\xi_S = S$.
\ei
\ep

Next we associate closure operators and $\vee$-congruences, making
explicit a construction suggested in \cite[Theorem 
6.3.7]{RS}.  

\bp
\label{clocon}
Let $L$ be a lattice, let $\xi:L \to L$ be a closure operator and
let $\rho$ be a $\vee$-congruence on $L$. Then:
\bi
\item[(i)] $\ker\xi$ is a $\vee$-congruence on $L$.
\item[(ii)] $\rho$ induces a
  closure operator $\eta_{\rho}:L \to L$ defined by $x\eta_{\rho} =
  \vee (x\rho) =$ {\rm max}$\, (x\rho)$.
\item[(iii)] $\eta_{\ker\xi} = \xi$ and $\ker\eta_{\rho} = \rho$.
\ei
\ep

\proof
(i) Clearly, $\ker\xi$ is an equivalence relation. Hence we must show
that $x\xi = y\xi$ implies $(x \vee z)\xi  = (y  \vee z)\xi$ for all
$x,y,z \in L$. Now $(x \vee z)\xi  = (x\xi \vee z\xi)\xi = (y\xi \vee
z\xi)\xi = (y  \vee z)\xi$
by Proposition \ref{cla}(i).

(ii) Axioms (C1) and (C3) follow immediately from $x \in x\rho$ and
$x\eta_{\rho}\rho = x\rho$. 

Assume that $x \leq y$ in $L$. Then $y = (x \vee y)$, hence 
$$y\rho = (x \vee y)\rho = (x\rho \vee y\rho) = (x\eta_{\rho}\rho \vee
y\rho) = (x\eta_{\rho} \vee y)\rho$$ and so 
$x\eta_{\rho} \leq (x\eta_{\rho} \vee y) \leq y\eta_{\rho}$ and (C2)
holds as well.

(iii) Let $x,y \in L$. Suppose that $y\xi = x\xi$. By (C1), we get
$y \leq y\xi = x\xi$ and so $x\eta_{\ker\xi} = \max\, (x(\ker\xi)) =
x\xi$. Thus $\eta_{\ker\xi} = \xi$.

On the other hand, 
$$(x,y) \in \ker\eta_{\rho} \iff x\eta_{\rho} = y\eta_{\rho} \iff
\max\, (x\rho) = \max\, (y\rho) \iff (x,y) \in \rho,$$
therefore $\ker\eta_{\rho} = \rho$.
\qed

Since $\vee$-congruences are nothing but kernels of $\vee$-maps,
Proposition \ref{clocon} establishes a correspondence between kernels
of $\vee$-maps and closure operators.

We can combine Propositions
\ref{cla} and \ref{clocon} with the
representation of lattices by flats. Given $(L,E) \in \FLg$, let
$\isflats(L,E)$ denote the set of all $\wedge$-subsemilattices of
$\flats(L,E)$, equivalently described as subsets of $\flats(L,E)$ closed
under intersection. 

\bt
\label{clofla}
Let $(L,E) \in$ {\rm FLg}, %. Then there is a bijective correspondence, 
let $\rho$ be a $\vee$-congruence on $L$
and let $F \in$ {\rm ISFl}$\,(L,E)$. Then:
\bi
\item[(i)] $F_{\rho} = \{ Z_{\max\, (x\rho)} \mid x \in L \} \in$ {\rm
    ISFl}$\,(L,E)$. 
\item[(ii)] The relation $\rho_F$ on $L$ defined by $x \rho_F y$ if
$$\cap \{ Z \in F \mid Z_x \subseteq Z \} = \cap \{ Z \in F \mid Z_y
\subseteq Z \}$$ 
is a $\vee$-congruence.
\item[(iii)] $\rho_{F_{\rho}} = \rho$ and $F_{\rho_F} = F$.
\ei
\et

\proof
We combine the correspondences in Propositions
\ref{cla} and \ref{clocon} with the lattice isomorphism $\p: (L,\leq) \to
(\flats(L,E),\subseteq): x \mapsto Z_x$ from Proposition
\ref{ltom}. Note that a $\wedge$-subsemilattice of $L$ corresponds to
a subset of $\flats(L,E)$ closed under intersection. Thus we only have
to check that the three correspondences mentioned above yield the
claimed ones.

Starting from $\rho$, we get $\eta_{\rho}$ and $L\eta_{\rho}$ by Propositions
\ref{clocon} and \ref{cla}, respectively, and then $\{ Z_m \mid m \in
L\eta_{\rho} \} = F_{\rho}$ by applying $\p$.

Starting from $F$, application of $\p\inv$ gives us $F\p\inv$. Then Propositions
\ref{cla} and \ref{clocon} give us successively $\xi_{F\p\inv}$ and
$\ker\xi_{F\p\inv}$. Now, for all $x,y \in L$, we get
$$\begin{array}{lll}
x\xi_{F\p\inv} = y\xi_{F\p\inv}&\iff&\wedge\{ p \in F\p\inv \mid p
\geq x \} = \wedge\{ p \in F\p\inv \mid p
\geq y \}\\
&\iff&\cap\{ Z \in F \mid Z \supseteq Z_x \} = \cap\{ Z \in F \mid Z
\supseteq Z_y \}\\
&\iff&x \rho_F y
\end{array}$$
and we are done.
\qed

% \addcontentsline{toc}{subsection}{B. Decomposing $\vee$-maps}
\subsection{Decomposing $\vee$-maps}

Once again, we import to the context of finite lattices a concept
originated in semigroup theory. Following \cite[Section 5.2]{RS}, we
call an onto $\vee$-map a $\vee$-surmorphism and say that a $\vee$-surmorphism
$\p:L \to L'$ is a {\em maximal proper 
  $\vee$-surmorphism} (MPS) of lattices if $\ker\p$
is a minimal nontrivial $\vee$-congruence on $L$. This amounts to saying
that $\p$ cannot be factorized as the composition of two proper
$\vee$-surmorphisms. 

Given $a,b \in L$, let $\rho_{a,b}$ denote the equivalence relation
on $L$ defined by
$$x\rho_{a,b} = \left\{
\begin{array}{ll}
\{ a,b\}&\mbox{ if $x = a$ or $x = b$}\\
\{ x \}&\mbox{ otherwise}
\end{array}
\right.$$

\bp
\label{mps}
Let $\p:L \to L'$ be a $\vee$-surmorphism of lattices. Then: 
\bi
\item[(i)] If $\p$ is not one-to-one, then $\p$ factorizes as a
  composition of MPSs.
\item[(ii)] If $a$ covers $b$ and $b$ is smi, then $\rho_{a,b}$ is a
  minimal nontrivial $\vee$-congruence on $L$. 
\item[(iii)] $\p$ is an MPS if and only if $\ker \p = \rho_{a,b}$ for
  some $a,b \in L$ such that $a$ covers $b$ and $b$ is smi.
\ei
\ep

\proof
(i) Since $L$ is finite, there exists a minimal nontrivial
$\vee$-congruence $\rho_1 \subseteq \ker \p$ and we can factor $\p$ as
a composition $L \to L/\rho_1 \to L'$. Now we apply the same argument
to $L/\rho_1 \to L'$ and successively.

(ii) Let $x \in L$. We must prove that $(x\vee a,x\vee b) \in
\rho_{a,b}$. Since $b$ is
smi, $a$ is the unique element of $L$ covering $b$. Hence either $(x
\vee b) = b$ or $(x \vee b) \geq a$. In the first case, we get $x \leq b$
and so $(x \vee a) = a$; in the latter case, we get $(x \vee b) = (x \vee
(x \vee b)) \geq (x \vee a) \geq (x \vee b)$ and so $(x \vee b) = (x \vee
a)$. Hence  $(x\vee a,x\vee b) \in \rho_{a,b}$ and so $\rho_{a,b}$ is a
(nontrivial) $\vee$-congruence on $L$. Minimality is obvious.

(iii) Assume that $\p$ is an MPS and let $a \in L$ be maximal among the
elements of $L$ which belong to a nonsingular $\ker\p$ class. Then
there exists some $x \in L \setminus \{ a \}$ such that $x\p =
a\p$. It follows that $(x \vee a)\p = (x\p \vee a\p) = a\p$ and so by
maximality of $a$ we get $(x \vee a) = a$ and so $x < a$.
Then there exists some $b \geq x$ such that $a$ covers
$b$. Since every $\vee$-map preserves order, we get $a\p = x\p \leq
b\p \leq a\p$ and so $a\p = b\p$. 

Suppose that $b$ is not smi. Then $b$ is covered by some other element
$c \neq a$, hence $b = (a \wedge c)$ and $a,c < (a \vee c)$. It follows
that $(a \vee c)\p = (a\p \vee c\p) = (b\p \vee c\p) = c\p$. Since $c \neq
(a \vee c)$ and $a < (a \vee c)$, this contradicts the maximality of
$a$. Thus $b$ is smi. Since $\rho_{a,b} \subseteq \ker\p$, it follows
from (ii) that $\ker\p = \rho_{a,b}$.

The converse implication is immediate.
\qed
 
We prove next the dual of Proposition \ref{mps} for injective
$\vee$-maps.  We
say that a $\vee$-map $\p:L \to L'$ is a {\em maximal proper
  injective $\vee$-map} (MPI) of lattices if $\p$ is injective and $L\p$
is a maximal proper $\vee$-subsemilattice of $L'$. This amounts to saying
that $\p$ cannot be  factorized as the composition of two proper
injective $\vee$-maps.

\bp
\label{mpi}
Let $\p:L \to L'$ be an injective $\vee$-map of lattices. Then: 
\bi
\item[(i)] If $\p$ is not onto, then $\p$ factorizes as a
  composition of MPIs.
\item[(ii)] If $a \in L'\setminus \{ B \}$ is sji, then the inclusion
  $\iota: L' \setminus \{ a\} \to L'$ is an MPI of lattices. 
\item[(iii)] $\p$ is an MPI if and only if $L\p = L' \setminus \{ a\}$ for
  some sji $a \in L\setminus \{ B \}$.
\ei
\ep

\proof
(i) Immediate since $L'$ is finite and each proper injective
$\vee$-map increases the number of elements.

(ii) Let $x,y \in L' \setminus \{ a\}$. Since $a$ is sji, the join of
$x$ and $y$ in $L'$ is also the join of
$x$ and $y$ in $L' \setminus \{ a\}$. Hence $L' \setminus \{ a\}$ is a
$\vee$-semilattice and therefore a lattice with the determined
meet. Since $(x\iota \vee y\iota) = (x \vee y) = (x \vee y)\iota$,
then $\iota$ is a $\vee$-map. Since $|L' \setminus \im\iota| = 1$, it
must be an MPI.

(iii)  
Assume that $\p$ is an MPI. Let $a$ be a minimal element of $L'
\setminus L\p$. We claim that $a$ is an sji in $L'$. Otherwise, by
minimality of $a$, we would have $a = (x\p \vee y\p)$ for some $x,y
\in L$. Since $\p$ is a $\vee$-map, this would imply $a = (x \vee
y)\p$, contradicting  $a \in L' \setminus L\p$. 

Thus $a$ is an sji in $L'$ and we can factor $\p:L \to L'$ as the composition
of $\p:L \to L' \setminus \{ a\}$ with the inclusion $\iota: L'
\setminus \{ a\} \to L'$. Since $\p$ is an MPI, then $\p:L \to L'
\setminus \{ a\}$ must be onto as required.

The converse implication is immediate.
\qed

\bt
\label{csi}
Let $\p:L \to L'$ be a $\vee$-map of lattices. Then $\p$ factorizes as a
  composition of MPSs followed by a composition of MPIs.
\et

\proof
In view of Propositions \ref{mps} and \ref{mpi}, it suffices to
note that $\p$ can always be factorized as $\p = \p_1\p_2$ with $\p_1$
a $\vee$-surmorphism and $\p_2$ an injective $\vee$-map. This can be
easily achieved taking $\p_1:L \to L\p$ defined like $\p$, and $\p_2:L\p
\to L'$ to be the inclusion.
\qed

We can produce a partial version of this result for the category FLg:

\bc
\label{csic}
Let $\p:(L,E) \to (L',E')$ be a $\vee$-surmorphism in {\rm FLg}. Then
the decomposition of $\p$ as a composition of MPSs constitutes a
composition of maps in {\rm FLg}.
\ec

\proof
Clearly, the $\vee$-generating set $E$ has a canonical correspondent
$E\rho_{a,b}$ in the construction $L/\rho_{a,b}$, and the restriction
$\p|_E E \to E' \cup \{ B \}$ factors adequately if $\p$ is onto. 
\qed

The analogous result fails for injective $\vee$-maps. For instance,
it is easy to see that the chain of inclusions
$$
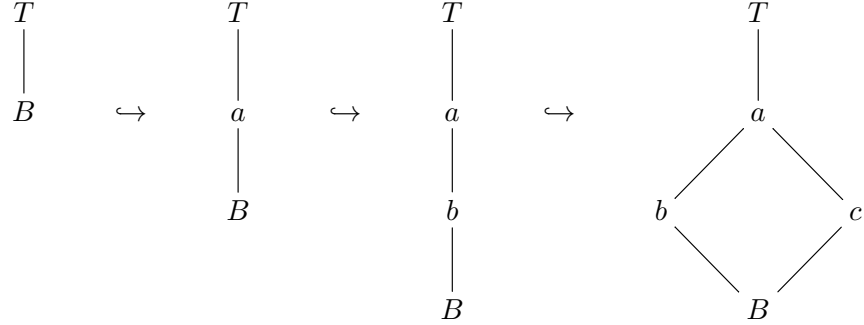
\xymatrix{
T \ar@{-}[d] && T \ar@{-}[d] && T \ar@{-}[d] &&& T \ar@{-}[d] &\\
B &\hookrightarrow & a \ar@{-}[d] &\hookrightarrow & a \ar@{-}[d]
&\hookrightarrow &&a  \ar@{-}[dl] \ar@{-}[dr] &\\
&&B && b \ar@{-}[d] && b && c \\
&&&&B &&& B  \ar@{-}[ul]
\ar@{-}[ur] & 
}$$
can induce no chain of inclusions between $\vee$-generating sets when
we consider $E = \{ T \}$ and $E' = \{ T,b,c\}$.

% \addcontentsline{toc}{subsection}{C. Geometry}
\subsection{Geometry}
\label{sgeom}

Let $P$ be a finite nonempty
set and let $\L$ be a nonempty subset of $2^P$. We shall always assume
that $P \cap 2^P = \emptyset$. We say that $(P,\L)$
is a {\em partial euclidean geometry} (abbreviated to PEG) if
the following axioms are satisfied:
\bi
% \item[(G1)] $P \subseteq \cup \L$;
\item[(G1)] if $L,L' \in \L$ are distinct, then $|L \cap L'| \leq 1$;
\item[(G2)] $|L| \geq 2$ for every $L \in \L$.
\ei 
The elements of $P$ are called {\em points} and the elements of $\L$
are called {\em lines}.  Given $p \in P$, we denote by $\L(p)$ the set
of all lines 
containing $p$.

The concept of PEG is  an abstract
combinatorial generalization of the following geometric situation: 

Consider a finite set of lines $\L$ in the euclidean space
$\RR^n$. Consider also a finite subset $P$ of $\cup\L 
\subset \RR^n$ such that:
\bi
% \item if $L,L' \in \L$ are distinct, then $|L \cap L'| \leq 1$;
\item if $L,L' \in \L$ and $L \cap L' = \{ p \}$, then $p \in P$;
\item $|L \cap P| \geq 2$ for every $L \in \L$.
\ei
Representing each $L\in \L$ by $L\cap P$, it follows that $(\L,P)$ constitutes a
PEG. It is well known that not all PEG's can be represented over an
euclidean space (nor 
any field) (see \cite{Cox} and \cite[Section 2.6]{Gru}).

In view of Proposition \ref{rtoh}, if $\het L = 2$, the subsets of $L
\setminus \{ B \}$ with at most two elements are the only c-independent
subsets of $L$. What about the case $\het L = 3$? This is the starting
point for a digression into some interesting connections between
c-independence and geometry.

% \subsection{Height 3}
% \label{subs3}

Given a lattice $L$ and $\ell \in \mathring{L} = L \setminus \{
T,B\}$, we define 
$$\ell\iota = \ell\dw \cap \mathring{L}.$$
If $(L,E) \in \FLg$, we define
$$\lin (L,E) = \{ \ell\iota \cap E : \ell \in \mathring{L},\;
|\ell\iota \cap E| \geq 2
\},$$
$$\poi (L,E) = E,\quad \geo (L,E) = (\poi (L,E),\lin
(L,E)).$$  

\bt
\label{geolat}
Let $(L,E) \in$ {\rm FLg} with {\rm ht}$\, L = 3$. Then {\rm
  Geo}$\,(L,E)$ is a PEG. 
\et

\proof
Axiom (G2) holds trivially, it remains (G1) to be
checked. Let $k, \ell \in \mathring{L}$ be such that $|k\iota \cap
\ell\iota \cap E| \geq 2$. We must show that $k = \ell$.

Suppose that $k \neq \ell$. Without loss of generality, we may assume
that $k > (k\wedge \ell)$. On the other hand, if $e_1 ,e_2$ are distinct
elements of $k\iota \cap
\ell\iota \cap E$, then we have $e_1,e_2 \leq (k\wedge \ell)$ and we may
assume that $e_1 < (k\wedge \ell)$. Thus we obtain a chain
$$B < e_1 < (k\wedge \ell) < k < T$$
in $\L$, contradicting $\het L = 3$. Therefore $k = \ell$ and we are done.
\qed

% It is not difficult to see that we can ensure $\poi (L,E) =
% E \setminus\{ 1 \}$ if $L$ has no
% maximal chain of length 2, which amounts to say that $L$ is modular!

Next we associate a matroid to every lattice $L$ of heigth 3:  we
define $\matroid_0 L$ to contain: 
\bi
\item
all the $i$-subsets of $L \setminus \{ B\}$ for $i \leq 2$;
\item
all the 3-subsets $X$ of $L \setminus \{ B\}$ such that $\vee X = T$.
\ei
Note that the latter condition is
equivalent to saying that $X \not\subseteq \ell\iota$ for every $\ell \in
\mathring{L}$. Finally, write $\matro L = (L \setminus \{ B\},
\matroid_0 L)$.

\bt
\label{assmat}
Let $L$ be a lattice of height 3. Then {\rm
  Mat}$\, L$ is a matroid. 
\et

\proof
It is immediate that $\matroid_0 L$ is a hereditary
collection. Let $\{ x,y \}$ be a 2-subset of $L \setminus \{ B\}$ and
let $\{ \ell_1,\ell_2,\ell_3 \}$ be a 3-subset of $L \setminus \{ B\}$
satisfying $(\ell_1 \vee \ell_2 \vee \ell_3) = T$. We must show that
$(x \vee y \vee \ell_i) = T$ for some $i \in \hat{3}$. Suppose
not. Then $(x \vee y \vee \ell_i) = k_i < T$ for every $i \in
\hat{3}$. 

Suppose first that $k_1 = k_2 = k_3$. Then $\ell_i \leq k_1 < T$ for
$i = 1,2,3$, contradicting  $(\ell_1 \vee \ell_2 \vee \ell_3) =
T$. Hence we have $(k_1 \wedge k_2 \wedge k_3) < k_j$ for some $j \in
\hat{3}$. Since $x,y \leq k_i$ for $i = 1,2,3$, we get $x,y \leq (k_1
\wedge k_2 \wedge k_3)$. Since $x$ and $y$ are distinct, we may assume
that $x < (k_1 \wedge k_2 \wedge k_3)$ and so we get a chain in $L$ of
the form
$$B < x < (k_1 \wedge k_2 \wedge k_3) < k_j < T,$$
contradicting $\het L = 3$. Hence $(x \vee y \vee \ell_i) = T$ for some
$i \in \hat{3}$ and so $\{ x,y,\ell_i\} \in \matro L$. It follows that 
$\matroid_0 L$ is a matroid.
\qed

Next we use $\matro L$ to characterize the c-independent subsets of
$L$. To do so, we introduce one more concept: a 3-subset $X \subseteq
L \setminus \{ B\}$ is called a {\em potential line} (of $\geo L$) if $|X \cap
\ell\iota| \leq 1$ for every $\ell \in \mathring{L}$. This is
equivalent to saying that $(x \vee y) = T$ for any distinct $x,y \in X$.

\bt
\label{chacind}
Let $L$ be a lattice of height 3 and let $X \subseteq L \setminus \{
B\}$. Then the following conditions are equivalent: 
\bi
\item[(i)] $X$ is c-independent;
\item[(ii)] $X \in$ {\rm Mat}$_0L$ and is not a potential line;
\item[(iii)] $|X| \leq 2$ or ($|X| = 3$, $\vee X = T$ and $(x\vee y) < T$ for
  some distinct $x,y \in X$).
\ei
\et

\proof
(i) $\Rw$ (ii) Suppose that $X$ is c-independent. We may assume that
$|X| = 3$. 
By Proposition \ref{cindch}, we may write $X = \{ x,y,z\}$ to get a
chain $(x \vee y \vee z) > (x \vee y) > x > 0$. Since $\het L = 3$, it
follows that $\vee X = T$ and so $X \in \matroid_0 L$. Since $(x \vee
y) < T$, $X$ is not a potential line.

(ii) $\Rw$ (iii). Immediate.

(iii) $\Rw$ (i). The case $|X| = 2$ following from Proposition
\ref{onetwo}, assume that $|X| = 3$, $\vee X = T$ and $x\vee y < T$ for some
distinct $x,y \in X$. Since $x \neq y$, we may assume that $(x \vee y) >
x$ and so we get a chain $\vee X > (y \vee x) > x$. By
Proposition \ref{cindch}, $X$ is c-independent.
\qed

Next we associate a $\vee$-generated lattice of height 3 to every PEG
$\G = (P,\L)$ with at least two lines: 
let $\lat \G = P \cup \L \cup \{ B,T\}$, where $x \leq y$ if and only
if
$$x = B \hspace{.5cm}\mbox{or}\hspace{.5cm} y = T
\hspace{.5cm}\mbox{or}\hspace{.5cm} (x\in P \mbox{ and } y \in \L \mbox{
  and } x \in y).$$ 
It is immediate that $\lat \G$ is a lattice of height 3. Moreover, if $W = \{
p_1,\ldots,p_k\} \in \L$, then $W = (p_1\vee \ldots \vee p_k)$, and we
can also get the top T as the join of two lines. Thus $(\lat \G,P) \in
\FLg$. 

\bp
\label{geolatpeg}
Let $\G = (P,\L)$ with $|\L| \geq 2$. Then {\rm Geo}$\,(${\rm Lat}$\,
\G,P) = \G$.
\ep

\proof
It follows from the definitions that <$\geo(\lat\G,P)$ is
of the form $(P,\L')$. If $p \in P$, then $p\iota = \{ p
\}$ in $\lat G$. If  $W = \{
p_1,\ldots,p_k\} \in \L$, then $W\iota = \{ W,
p_1,\ldots,p_k\}$. Thus, by definition of the construction $\geo$, the
elements of $\L'$ are of the form $W\iota \cap P$ for $W \in
\L$. Since $W\iota \cap P = W$, we get $\L' = \L$ and so
$\geo(\lat\G,P) = \G$.
\qed

% \subsection{Height 4}

We say that $h \in \flats(L,E)$ is a {\em hyperplane} of $(L,E)$ if
$h$ is maximal in $\flats(L,E) \setminus \{ E \}$. For height 4, we
can prove the following result:

\bp
\label{fourind}
Let $(L,E) \in$ {\rm FLg} have height 4 and let $X$ be a 4-subset of
$E$. Then the 
following conditions are equivalent: 
\bi
\item[(i)] $X$ is c-independent;
\item[(ii)] every 3-subset of $X$ is c-independent and $|X \cap h| =
  3$ for some hyperplane $h$ of $L$.
\ei
\ep

\proof
(i) $\Rw$ (ii). The first claim follows from c-independent sets being
closed under inclusion. On the other hand, by Proposition
\ref{cindch}, $X$ admits an 
enumeration $x_1, x_2, x_3, x_4$ such that 
$$\clos_LX \supset \clos_L(x_2,x_3,x_4) \supset \clos_L(x_3,x_4)
\supset \clos_L(x_4) \supset \emptyset.$$
Let $h = \clos_L(x_2,x_3,x_4)$. Since $\het \flats(L,E) = \het L = 4$
by Proposition \ref{ltom}, it follows that $h$ is a hyperplane of
$(L,E)$. Clearly, $x_1 \notin \clos_L(x_2,x_3,x_4)$ and so $|X \cap h| = 3$.

(ii) $\Rw$ (i). Write $X \setminus h = \{ x_1 \}$. Since $X \cap h$ is
c-independent, it follows from Proposition \ref{threeind} that
$X\cap h$ admits an enumeration $x_2, x_3, x_4$ such that
$(x_2 \vee x_3 \vee x_4) > (x_3 \vee x_4) > x_4$. Since $\{
x_2,x_3,x_4\} \subseteq h$ implies $\cl_L(x_2,x_3,x_4) \subseteq
\cl_L(h) = h$, we get $x_1 \notin \cl_L(x_2,x_3,x_4) = Z_{x_2 \vee x_3
  \vee x_4}$ (by (\ref{clz})) and so $(x_1 \vee
x_2 \vee x_3 \vee x_4) > (x_2 \vee x_3 \vee x_4)$. Thus $X$ is
c-independent by Proposition
\ref{cindch}.
\qed

% \seection{Arbitrary height}

We can generalize to higher dimensions the concept of PEG
to get generalizations of some results obtained in the height 3 case, namely
Proposition \ref{geolatpeg}. For technical reasons, we include the
full space of points as the highest dimension subspace, but it could
as well be omitted.

For $m \geq 3$, we say that $(P_1,P_2,\ldots,P_m)$ is an $m$-PEG over
a finite set $E$ 
if:
\bi
\item[(J1)] $P_1,\ldots, P_m$ are mutually disjoint 
subsets of $2^E$ and $P_m = \{ E \}$;
\item[(J2)] $\forall p \in P_1 \; |p| = 1$;
\item[(J3)] $\forall i \in \{ 2, \ldots,m\} \; \cup P_i \subseteq \cup P_{1}$;
\item[(J4)] $\forall i \in \{ 2, \ldots,m\} \; \forall p \in P_i \; \exists
  q \in P_{i-1}: q \subset p$;
\item[(J5)] $\forall i,j \in \{ 2, \ldots,m\} \; \forall p \in P_i \;
  \forall q \in P_j$, one of the following five conditions holds:
\bi
\item[(J5a)]
$p \cap q = \emptyset$;
\item[(J5b)]
$p \cap q \in P_r$ for some $r < i,j$;
\item[(J5c)]
$i < j$ and $p \subset q$;
\item[(J5d)]
$i > j$ and $p \supset q$;
\item[(J5e)]
$p = q$.
\ei
\ei

The 3-PEG case corresponds to our original concept of PEG, replacing
each point $p$ by $\{ p \}$ and adding the full subspace $P$. A natural
example for the general case is given by a (finite) collection of affine
subspaces of various ranks in an euclidean space of arbitrary
dimension, where the subspaces are defined through collections of
points from a finite set $E$. 

Two $m$-PEGs $(P_1,P_2,\ldots,P_m)$ (over $E$) and
$(P'_1,P'_2,\ldots,P'_m)$ (over $E'$) are said to be isomorphic if
there exists a bijection $\p:E \to E'$ such that 
$$\{ e_1,\ldots, e_k \} \in P_i \iff \{ e'_1,\ldots, e'_k \} \in P_i$$
holds for all $i \in \hat{m}$ and $e_1,\ldots, e_k \in E$. A
particularly important case arises with the canonical bijections $X
\to \wt{X}$, where $\wt{X} = \{ \{ x \} \mid x \in X \}$.

Given a lattice $L$ and $\ell \in \L$, we define $\hei(\ell)$ to be
the maximum length $n$ of a chain $\ell = \ell_0 > \ell_1 > \ldots >
\ell_n$ in $\L$. Obviously, $\het L = \hei(1)$.

We recall
now the notions of atom and atomic lattice. If we denote
by $A_L$ the set of atoms of $L$, then $L$ is atomic if and only if
$(L,A_L) \in \FLg$.

Given an atomic lattice $L$ of height $m$, we define $\geo(L,A_L) =
(P_1,\ldots,P_m)$ by
$$P_i = \{ \ell\iota \cap A_L \mid \ell \in L,\; \hei(\ell) =
i \}\hspace{1cm}(i = 1,\ldots,m).$$
We claim that $\geo(L,A_L)$ is an $m$-PEG over $A_L$. Axiom (J1) follows from
(\ref{ltom1}) and $T\iota = \mathring{L}$, and (J2) is
immediate. 
Since the elements of $P_1$ are of the form $\{ a \}$ for $a \in
\A_L$, (J3) holds. 
Since every element of $L$ of height 
$i$ covers some element of height $i-1$, (J4) follows. Finally, let
$k, \ell \in L$ have height $i$ and $j$, respectively. Since 
$(k\iota \cap A_L)\cap(\ell\iota \cap A_L) = (k\wedge\ell)\iota \cap
A_L$, then (J5a) or (J5b) hold if $k\wedge \ell < k,\ell$. Hence we are
left with the cases $k < \ell$, $k > \ell$ and $k = \ell$ which give
us respectively (J5c), (J5d) and (J5e) in view of (\ref{ltom1}). Thus
$\geo(L,A_L)$ is an $m$-PEG over $A_L$.

Conversely, given an $m$-PEG $\G = (P_1,\ldots,P_m)$ over $E$, we
define the poset $\latti_0\G = \{ \emptyset \} \cup P_1 \cup \ldots \cup
P_m$, ordered by inclusion. By (J5), $\latti_0 \G$ is closed under
intersection and constitutes then a $\wedge$-semilattice with bottom
element $\emptyset$ and top element $E$. Hence $\latti_0 \G$ is a lattice
with $p\vee q = \cap\{ r \in \lat\G \mid p\cup q \subseteq r\}$ (note
that $E \in \latti_0\G$). We
claim that $P_1$ is the set of atoms of $\latti_0 G$. 

Indeed, it follows from (J4) that any atom is contained in $P_1$, and
the converse is a consequence of (J2). Now, by (J2) and (J3), $\latti_0
\G$ is atomic and we can define $\lat\G = (\latti_0 \G,P_1) \in \FLg$.  

\bt
\label{mpeglat}
Let $\G$ be an $m$-PEG and let $L$ be an atomic lattice. Then
\bi
\item[(i)] {\rm Geo}$\,${\rm Lat}$\,\G \cong \G$;
\item[(ii)] {\rm Lat}$\,${\rm Geo}$\,(L,A_L) \cong (L,A_L)$.
\ei
\et

\proof
(i) Let $\G = (P_1,\ldots,P_m)$ be an $m$-PEG over $E$. We claim that 
\beq
\label{mpeglat1}
p \in P_i \hspace{.5cm} \Rw \hspace{.5cm} \hei(p) = i \mbox{ in } \lat\G.
\eeq
We use induction on $i$. The case $i = 1$ follows from (J2), hence we
assume that $p \in P_i$ with $i > 1$ and (\ref{mpeglat1}) holds for
$i-1$. By (J4) and the induction hypothesis, we have $\hei(p) \geq
i$. Suppose that $\hei(p) > i$. This would imply that there would
exist distinct $q,r \in P_j$ for some $j$ such that $q \subset r$,
contradicting (J5). Hence $\hei(p) = i$ and (\ref{mpeglat1}) holds.

Hence $\lat\G$ is a lattice of height $m$ with set of atoms $P_1$
 and we may write
$\geo\lat\G = (P'_1,\ldots,P'_m)$, an $m$-PEG over $P_1 \subseteq
\wt{E}$. In view of (\ref{mpeglat1}), the elements of $P'_i$ are of
the form $p\iota \cap P_1$ for $p \in P_i$.

For all $p \in P_i$ and $e \in E$, we have $\{ e \} \in
p\iota$ in $\lat\G$ if and only if $\{ e \} \subseteq p$ in $\G$ if
and only if $e \in p$, hence $p\iota \cap P_1 = \wt{p}$ and so $e
\mapsto \wt{e}$ $(e \in E)$ induces an isomorphism between $\G$ and
$\geo\lat\G$.

(ii) Let $L$ be an atomic lattice. Then the elements of
$\latti_0\geo(L,A_L)$ are of the form $\ell\iota \cap A_L$ (note that
$\emptyset = B\iota \cap A_L$), ordered by inclusion. By
(\ref{ltom1}), $\p: L \to \latti_0\geo(L,A_L)$ defined by $\ell\p =
\ell\iota \cap A_L$ is an isomorphism of posets and therefore of
lattices. 
Since the set of atoms of a lattice is uniquely determined, $\p$
induces an isomorphism from $(L,A_L)$ onto $\lat\geo(L,A_L)$.
\qed

%\addcontentsline{toc}{subsection}{D. Strong maps}
\subsection{Strong maps}

The concept of strong map can be defined for lattices
and for hereditary collections. We start discussing the lattice case.

We say that a mapping $\p:(L,E) \to
(L',E')$ of lattices is a {\em strong map} if 
$$\forall Z \in \flats (L',E') \hspace{.5cm} (Z\cup \{ B \})\p\inv
\cap E \in \flats (L,E).$$
Now assume that $\p$ is a $\vee$-map in FLg. In particular, $E\p
\subseteq E' \cup 
\{ 0 \}$. Given $X \subseteq E$, write $\oo{X} = \clos_LX.$ 
Then $\p$ induces a map $\oo{\p}:
\flats (L,E) \to \flats (L',E')$ by $Z\oo{\p} = \oo{Z\p \setminus
\{ B \}}$. 

\bp
\label{strong}
Let $\p:(L,E) \to
(L',E')$ be a $\vee$-map in {\rm FLg}. Then:
% following conditions are equivalent:
\bi
\item[(i)] $\p$ is a strong map;
\item[(ii)] $\oo{\p}$ is a $\vee$-map and $\oo{e}\,\oo{\p} = \oo{e\p
    \setminus \{ B \}}$
  for every $e \in E$.
\ei
\ep

\proof
(i) Let $Z' \in \flats(L',E')$. Then $Z' = Z_{\ell'} = \ell'\dw \cap
E'$ for some $\ell' \in L'$. Assume that $(Z' \cup \{ B
\})\p\inv = \{ x_1,\ldots, x_k\}$, and write $\ell = (x_1 \vee \ldots
\vee x_k)$. We claim that
\beq
\label{injsm1}
(Z' \cup \{ B \})\p\inv \cap E = Z_{\ell} \in \flats(L,E).
\eeq

Indeed, let $e \in (Z' \cup \{ B \})\p\inv \cap E$. Then $e = x_i$ for
some $i \in \hat{k}$ and so $e \leq \ell$. Thus $e \in
Z_{\ell}$. Conversely, assume that $e \in
Z_{\ell}$. Since $Z_{\ell} \subseteq E$, it suffices to show that $e\p
\in Z' \cup \{ B \}$. Since $E\p \subseteq E' \cup \{ B \}$, all we
need is to prove that $e\p \leq \ell'$. Now $e \leq \ell = (x_1 \vee \ldots
\vee x_k)$ and $\p$ being a $\vee$-map yields $e\p \leq (x_1 \vee \ldots
\vee x_k)\p = (x_1\p \vee \ldots
\vee x_k\p)$. Since $x_i\p \in Z' \cup \{ B
\}$ and therefore $x_i\p \leq \ell'$ for every $i$, we get $e\p \leq
\ell'$ as required. Thus (\ref{injsm1}) holds and so $\p$ is a strong map.

(ii)  
Let $x,y \in L$. To show that $\oo{\p}$ is a $\vee$-map, we need to
show that $(Z_x \vee Z_y)\oo{\p} = 
(Z_x\oo{\p} \vee Z_y\oo{\p})$, i.e. 
$$(\oo{Z_x \cup Z_y})\oo{\p} =
\oo{Z_x\oo{\p} \cup Z_y\oo{\p}}.$$
Since $\oo{\p}$ is order-preserving, we get $Z_x\oo{\p} \cup
Z_y\oo{\p} \subseteq (\oo{Z_x \cup Z_y})\oo{\p}$ and so $\oo{Z_x\oo{\p} \cup
Z_y\oo{\p}} \subseteq (\oo{Z_x \cup Z_y})\oo{\p}$ since $(\oo{Z_x \cup
Z_y})\oo{\p}$ is closed. The inclusion
\beq
\label{strong1}
(\oo{Z_x \cup Z_y})\oo{\p} \subseteq
\oo{Z_x\oo{\p} \cup Z_y\oo{\p}}.
\eeq
remains to be proved. We have 
\beq
\label{strong2}
(\oo{Z_x \cup
Z_y})\oo{\p} = \oo{(\oo{Z_x \cup
Z_y})\p \setminus \{ B \}}.
\eeq
Let $Z = \oo{Z_x\oo{\p} \cup
Z_y\oo{\p}} \in \flats(L',E')$. Then $Z_x \cup Z_y \subseteq (Z \cup
\{ 0 \})\p\inv$.
By (i), we have $(Z \cup
\{ B \})\p\inv \cap E \in
\flats(L,E)$, hence  $\oo{Z_x \cup Z_y} \subseteq (Z \cup
\{ B \})\p\inv$ and so $(\oo{Z_x \cup Z_y})\p \setminus \{ B \} 
\subseteq Z$. Since $Z$ is closed, it follows that
$$(\oo{Z_x \cup Z_y})\oo{\p} \subseteq \oo{(\oo{Z_x \cup Z_y})\p
  \setminus \{ B \}} \subseteq Z,$$
hence (\ref{strong1}) holds and so $\oo{\p}$ is a $\vee$-map.

Let $e \in E$. We must show that $\oo{\oo{e}\p \setminus \{ B \}} =
\oo{e\p \setminus \{ B \}}$. The opposite inclusion being immediate,
we set $Z = \oo{e\p \setminus \{ B \}} \in \flats(L',E')$ and show
that $\oo{\oo{e}\p \setminus \{ B \}} \subseteq Z$. 

Clearly, $e \in (Z \cup
\{ B \})\p\inv$. By (i), $(Z \cup
\{ B \})\p\inv$ is closed and so $\oo{e} \subseteq (Z \cup
\{ B \})\p\inv$. Hence $\oo{e}\p \subseteq Z \cup
\{ B \}$ and so $\oo{e}\p \setminus \{ B \} \subseteq Z$. Since $Z$ is
closed, we get $\oo{\oo{e}\p \setminus \{ B \}} \subseteq Z$ as
required.
\qed

The following examples show that a strong map is not necessarily a
$\vee$-map, even if we assume injectivity or surjectivity:

\be
\label{cestrong}
Consider the inclusion $\iota:(L,E) \to (L',E')$ for the following
lattices, where the elements of 
the $\vee$-generating sets are marked with an asterisk:
$$
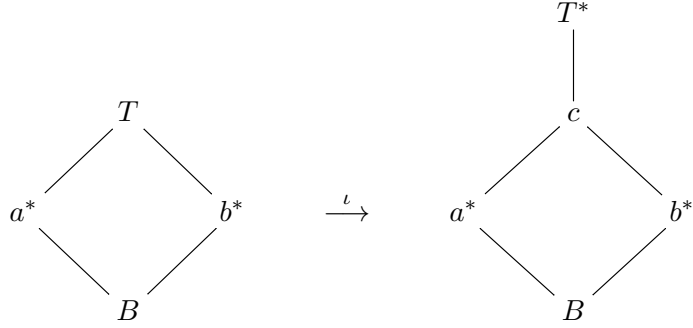
\xymatrix{
&&&&& T^* \ar@{-}[d] &\\
&T \ar@{-}[dl] \ar@{-}[dr] &&&&c \ar@{-}[dl] \ar@{-}[dr] &\\
a^* \ar@{-}[dr] &&b^* \ar@{-}[dl] &{\mapright{\iota}}&a^* \ar@{-}[dr]
&&b^* \ar@{-}[dl] \\
&B&&&&B&
}$$
Then $\iota$ is strong but not a $\vee$-map.
\ee

Indeed, since $\flats(L,E) = 2^E$, $\iota$ is strong. Since $(a\vee
b)\iota = T \neq c = (a\iota \vee b\iota)$, $\iota$ is not a $\vee$-map. 

\be
\label{cestrongtwo}
Consider the onto mapping $\p:(L,E) \to (L',E')$ defined by
$$
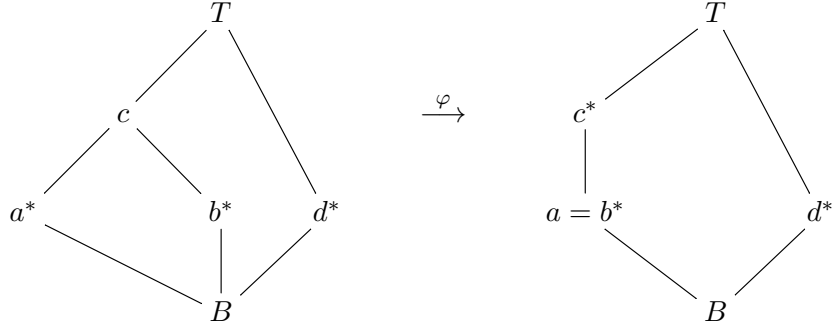
\xymatrix{
&&T \ar@{-}[dl] \ar@{-}[ddr] &&&&T \ar@{-}[dl] \ar@{-}[ddr] &\\
&c \ar@{-}[dl] \ar@{-}[dr] &&&{\mapright{\p}}&c^* \ar@{-}[d] &&\\
a^*\ar@{-}[drr] &&b^*\ar@{-}[d] &d^* \ar@{-}[dl] &&a=b^* \ar@{-}[dr]
&&d^* \ar@{-}[dl] \\
&&B&&&&B&
}$$
where the elements of the $\vee$-generating sets are marked with an asterisk.
Then $\p$ is strong but not a $\vee$-map.
\ee

Indeed, we have $\flats(L,E) = 2^E \setminus \{ ad, bd \}$, but the
inverse image of a flat from $(L',E')$ contains $a$ if and only if it
contains $b$. Thus $\p$ is strong. Since $(a\vee
b)\p = c \neq (a\p \vee b\p)$, then $\p$ is not a $\vee$-map.

We discuss now the concepts of strong and weak maps
for hereditary 
collections. Let $(E,H), (E',H')$ be hereditary
collections and let $\p:E \to E'$ be a mapping. We say that $\p$ is
a {\em weak map} (with respect to $(E,H), (E',H')$) if
$$\p|_X \mbox{ injective and } X\p \in H' \hspace{.5cm} \Rw
\hspace{.5cm} X \in H$$
holds for every $X \subseteq E$. Assume now that $(E,H), (E',H')$ are
boolean representable. We say that $\p$ is
a {\em strong map} (with respect to $(E,H), (E',H')$) if
$$X \in \flats(E',H') \Rw X\p\inv \in \flats(E,H)$$
holds for every $X \subseteq E'$.

\bp
\label{stwe}
Let $(E,H), (E',H')$ be boolean representable simple hereditary
collections and let $\p:E \to E'$ be a strong map. Then $\p$ is a
weak map.
\ep

\proof
Let $X \subseteq E$ and assume that $\p|_X$ is injective and $X\p \in
H'$. By Proposition \ref{hcl}, and since $\p|_X$ is injective, $X$
admits an enumeration 
$x_1,\ldots, x_k$ such that 
$$\clos(x_1\p,\ldots, x_k\p) \supset \clos(x_2\p,\ldots, x_k\p)
\supset \ldots \supset \clos(x_k\p).$$
We claim that
\beq
\label{stwe1}
\clos(x_1,\ldots, x_k) \supset \clos(x_2,\ldots, x_k)
\supset \ldots \supset \clos(x_k).
\eeq
Indeed, suppose that $\clos(x_i,\ldots, x_k) = \clos(x_{i+1},\ldots,
x_k)$ for some $i \in \{1,\ldots, k-1\}$. Then $x_i \in \clos(x_{i+1},\ldots,
x_k)$. Since $\{ x_{i+1},\ldots,
x_k \} \subseteq (\clos(x_{i+1}\p,\ldots, x_k\p))\p\inv$ and the
latter is closed since $\p$ is a strong map, we get $x_i \in
\clos(x_{i+1},\ldots, x_k) \subseteq (\clos(x_{i+1}\p,\ldots, x_k\p))\p\inv$
and so $x_i\p \in \clos(x_{i+1}\p,\ldots, x_k\p)$, contradicting
$\clos(x_{i}\p,\ldots, x_k\p) \supset \clos(x_{i+1}\p,\ldots,
x_k\p)$. Hence $\clos(x_i,\ldots, x_k) \supset \clos(x_{i+1},\ldots,
x_k)$ for every $i$ and so (\ref{stwe1}) holds. Now Proposition
\ref{hcl} yields $X \in H$ and so $\p$ is a weak map.
\qed

\section*{Acknowledgments}

The second author acknowledges support from the European Regional
Development Fund through the programme COMPETE
and by the Portuguese Government through FCT (Funda\c c\~ao para a Ci\^encia e a
Tecnologia) under the project PEst-C/MAT/UI0144/2011.

\end{document}